\documentclass[10pt, twoside]{article}
\usepackage{verbatim, version}
\usepackage{amsmath,amssymb,latexsym, amsthm}
\usepackage{graphicx, subfigure, epsfig, psfrag}
\usepackage{amsfonts, amssymb, dsfont, mathrsfs, color} 
\usepackage{diagbox}

\graphicspath{{../figs/}}
\newlength{\tmarg}
\newlength{\bmarg}
\newlength{\lmarg}
\newlength{\rmarg}

\newcommand{\beal}[1]{\begin{subequations}\label{#1}\begin{eqnarray}}
\newcommand{\bea}{\begin{subequations}\begin{eqnarray}}
\newcommand{\eea}{\end{eqnarray}\end{subequations}}

\newcommand{\bel}[1]{\begin{equation}\label{#1}}

\newcommand{\bmath}{\begin{displaymath}}
\newcommand{\emath}{\end{displaymath}}

\newcommand{\bear}{\begin{eqnarray}}
\newcommand{\eear}{\end{eqnarray}}

\newcommand{\bM}{\begin{displaymath}}
\newcommand{\eM}{\end{displaymath}}

\newcommand{\be}{\begin{equation}}
\newcommand{\ee}{\end{equation}}

\newcommand{\bfg}{\begin{figure}}
\newcommand{\efg}{\end{figure}}

\newcommand{\bc}{\begin{center}}
\newcommand{\ec}{\end{center}}
\newcommand{\vect}[1]{{\bf #1}}

\definecolor{gray}{rgb}{0.8,0.8,0.8}
\definecolor{black}{rgb}{0.0,0.0,0.0}
\definecolor{darkgray}{rgb}{0.4,0.4,0.4}
\definecolor{darkblue}{rgb}{0.0,0.0,0.5}

\definecolor{golden}{rgb}{0.9,0.8,0.0}
\definecolor{pink}{rgb}{1.0,0.0,1.0}

\def\been{\begin{enumerate}}
\def\enen{\end{enumerate}}
\def\beit{\begin{itemize}}
\def\enit{\end{itemize}}

\include{macros}

\begin{document}

\title{ Kernel-free boundary integral method for two-phase Stokes equations with discontinuous viscosity on staggered grids
%\thanks{This work was supported by National Science Foundation
%Grant DMS-9873275 and VIGRE Grant DMS-9983320. } 
}

%\thanks{ This work was in part supported by National Science Foundation
%Grant DMS-9873275 and VIGRE Grant DMS-9983320. The software part of the
%project was mostly performed under the auspices of the U.S. Department 
%of Energy by University of California Lawrence Livermore National 
%Laboratory under Contract No. B522274. } }

%\author{Wenjun Ying
%\thanks{School of Mathematical Sciences and Institute of Natural Sciences, 
%Shanghai Jiao Tong University, Minhang, Shanghai 200240, P.~R.~China 
%(wying@sjtu.edu.cn). } 
%}

\author{
  Haixia Dong\thanks{MOE-LCSM, School of Mathematics and Statistics, Hunan Normal University, Changsha, Hunan 410081, P. R. China
    ({\tt hxdong@csrc.ac.cn}).}
  \and
  Shuwang Li\thanks{Department of Applied Mathematics, Illinois Institute of Technology, Rettaliata Engineering Center, Room 11B, 10 W. 32nd Street, Chicago, IL60616, USA ({\tt sli@math.iit.edu}).}
  \and
  Wenjun Ying\thanks{Corresponding author. School of Mathematical Sciences, MOE-LSC, and Institute of Natural Sciences, Shanghai Jiao Tong University, Minhang,
Shanghai, 200240, P. R. China ({\tt wying@sjtu.edu.cn}).}
    \and
   Zhongshu Zhao\thanks{School of Mathematical Sciences, and Institute of Natural Sciences, Shanghai Jiao Tong University, Minhang, Shanghai, 200240, P. R. China 
   ({\tt zhaozs@sjtu.edu.cn}).}
 }
   
% \\ {\sl \small School of Mathematical Sciences, 
% Shanghai Jiao Tong University} 
%  \\ {\sl \small Minhang, Shanghai 200240, P.R.C.}
%  \\ {\sl \small wying@sjtu.edu.cn} }

%\date{\today}
\date{}

\maketitle

%\bc
%{\Large \bf An efficient Cartesian grid-based boundary integral method for scattering problems around irregular obstacles}
%
%{\bf \large Wenjun Ying}
%
%{\sl \small School of Mathematical Sciences and Institute of Natural Sciences, \\
%Shanghai Jiao Tong University, Minhang, Shanghai 200240, P.~R.~China}
%\ec

%------------------------------------------------------------------------------ 

\begin{abstract}
 A discontinuous viscosity coefficient makes the jump conditions of the velocity and normal stress coupled together, which 
brings great challenges to some commonly used numerical methods to obtain accurate solutions. To overcome the difficulties, a kernel free boundary integral (KFBI) method combined with a modified marker-and-cell (MAC) scheme is developed to solve the two-phase Stokes problems with discontinuous viscosity. The main idea is to reformulate the two-phase Stokes problem into a single-fluid Stokes problem by using boundary integral equations and then evaluate the boundary integrals indirectly through a Cartesian grid-based method. Since the jump conditions of the single-fluid Stokes problems can be easily decoupled, the modified MAC scheme is adopted here and the existing fast solver can be applicable for the resulting linear saddle system. 
The computed numerical solutions are second order accurate in discrete $\ell^2$-norm for velocity and pressure as well as the gradient of velocity, and also second order accurate in maximum norm for both velocity and its gradient, even in the case of high contrast viscosity coefficient, which is demonstrated in numerical tests.

\end{abstract}


%============================================================================== 
\section{Introduction}
Incompressible Stokes equations are used intensively for flows with small to modest Reynolds numbers.  Applications for Stoke interface problems with moving interface in computational fluid dynamics include multi-phase incompressible flows \cite{chang1996level,cogan2005modeling,gross2011numerical}, fluid structure interaction (FSI) problems \cite{hou2012numerical,mokbel2018phase,kim2019immersed} and so on.  

In the past decades, efficient and accurate numerical approaches to approximating interface problems  have received wide attention,  which can be particularly classified into two categories: methods with interface-fitted meshes and methods with interface-unfitted meshes according to the discretization of the physical domain. The former approach does not allow the interface to cut across any element, and the jump conditions across the interface can be incorporated into a standard numerical formulation, such as finite element method\cite{gross2007finite}, hybridizable discontinuous Galerkin method \cite{wang2013hybridizable}. It tends to capture discontinuities of the solution more accurately. However, generating an interface-fitted mesh of relatively high quality is challenging. Especially when the interface evolves with time, the generation of fitted mesh consumes much time and needs large memory. The latter approach is more desirable as it allows much simpler meshes independent of the location of the interface, such as Cartesian grids and quasi-uniform meshes. The success of the interface-unfitted method relies on how to effectively handle the jump conditions and many efforts have been made to it.

First, among the interface-unfitted method to solve the Stokes interface problem, the finite element method (FEM) is a popular choice.   Hansbo et al. \cite{hansbo2014cut} proposed a cut FEM to solve the Stokes interface problem, which weakly enforces the jumping condition on the interface with a weighted coefficient in the Nitsche’s numerical flux. Adjerid et al. \cite{adjerid2015immersed} presented an immersed $Q_1/Q_0$ discontinuous Galerkin FEM for solving the Stokes interface problem, which achieves optimal convergence in both the velocity and the pressure. This discontinuous immersed finite element space is designed according to the location of the interface and pertinent interface jump conditions. Then, this idea was applied to Stokes interface problem with moving interfaces \cite{adjerid2019immersed}. 
Later, many approaches based on immersed finite element method (IFEM) have been further developed, such as nonconforming IFEM \cite{jones2021class}, partially penalized IFEM \cite{chen2021p2}, Immersed $CR-P_0$ element \cite{ji2022immersed}. 
 In addition, there also exist other FEMs to solve Stoke interface problems including Nitsche's Extended FEM \cite{lehrenfeld2012nitsche,wang2015new,wang2019nonconforming,he2019stabilized}, XFEM \cite{chessa2003extended,gross2007extended,kirchhart2016analysis}, fictitious domain FEM \cite{lundberg2019distributed,sun2019fictitious,sun2020distributed},  corrected FEM \cite{laymuns2022corrected} and so on.

Second, another widely used numerical approach to solving the Stokes interface problem is immersed interface method (IIM) within the finite difference framework, which is originally proposed by LeVeque and Li \cite{leveque1994immersed} for solving elliptic interface problems and motivated by Peskin's immersed boundary (IB) method \cite{peskin1977numerical,peskin2002immersed} to improve the accuracy to at least second order, particularly near the interface \cite{li2001maximum,li2017accurate}.  IIM incorporates the correction terms computed from the jump conditions into the finite difference scheme,  thus it can capture the solution and its derivative jumps sharply. Later, the IIM has been widely used in interface problems, such as acoustic wave equations \cite{zhang1997immersed}, Stokes flow with elastic boundaries or surface tension \cite{leveque1997immersed}, Navier-Stokes problems \cite{li2001immersed,lee2003immersed,tan2009immersed}, fluid-solid interaction \cite{xu20083d}.  A detailed IIM overview can be found in the book by Li and Ito \cite{li2006immersed}. Especially, to solve incompressible 2D Stokes flow with discontinuous viscosity,  Li et al. developed an augmented approach using IIM \cite{li2007augmented}, which decomposes the incompressible Stokes equations into three Poisson equations. The bi-periodic boundary condition is assumed and a numerical boundary condition for the pressure should be designed \cite{li2007augmented}.
Then, Tan, Lim, and Khoo \cite{tan2011implementation} combined the augmented IIM with the MAC scheme to design an efficient algorithm for two-phase incompressible Stokes equations, which numerically produces second-order accurate solution for both velocity and pressure. Since the introduction of augmented variables destroyed the original nice matrix structure, a direct IIM \cite{chen2018direct} based on the MAC scheme is proposed to solve two-phase incompressible Stokes equations. In this approach,  the resulting linear system can be solved by the regular Uzawa iterative method.

Third, boundary integral methods (BIMs) are the third attractive computational technique for Stokes interface problems as they reformulate PDEs into integral equations on the domain boundary which reduces the dimensionality of the problem.
The major advantage of this method is that the unknown stress and velocity fields of the flow equations are only solved at the domain boundaries or the fluid interface, thus it avoids the generation of unstructured mesh and requires less computer memory. 
Especially for two-fluid Stokes equations in the free space, Layton \cite{layton2008efficient} reduced the problem into a single-fluid problem by solving a Fredholm integral equation of the second kind. The solution to the reduced problem is then computed on a finite domain using IIM. Such an approach can be applied to problems with a variety of boundary conditions or three dimensions by using appropriate Green's functions in boundary integrals. However, in many cases, it is difficult to obtain the analytical expression of Green's functions, for example, when Green's function is defined on a bounded domain and subject to a non-periodic boundary condition, or when it is associated with the variable coefficients differential operators. Ying and Henriquez \cite{ying2007kernel} proposed the kernel-free boundary integral  (KFBI) method for elliptic boundary value problems, which was then applied to various problems \cite{ying2013kernel,ying2013fast,ying2014kernel,dong2018hybridizable,xie2019fourth,xie2019high}.  The most significant merits of the KFBI method are  that it solves boundary value or interface problems in the framework of BIMs but does not need to know the analytical expression of Green's function, and it does not have singularity issues associated with traditional BIMs.

Very recently, a modified MAC scheme for the Stokes interface problem with constant viscosity in the framework of finite difference was presented and the corresponding rigorous error analysis was also provided in \cite{dong2022second}. %However, the modified MAC scheme only solves Stokes interface problem with constant viscosity. 
This work further extends the modified MAC scheme to solve two-phase Stokes equations with discontinuous viscosity in the framework of the KFBI method.  The proposed approach has a variety of advantages and is efficient as well as accurate even for large contrast viscosity across the interface.  More precisely, ingredients of the algorithm and contributions of this work include
\begin{itemize}
\item[(i)]By solving a Fredholm integral equation of the second kind, the original two-phase Stokes problem is reduced into a single-fluid case, which can be solved by the previous modified MAC scheme in \cite{dong2022second}. Unlike the aforementioned augmented methods\cite{li2007augmented,tan2011implementation} , which introduce additional variables and equations to decouple the jump conditions, the original problem with discontinuous viscosity is not discretized directly in this approach. Thus the coupled difficulty between the velocity and the normal traction due to the jump conditions, caused by the discontinuity of the fluid viscosity across the boundary, is successfully overcome. Furthermore, this method can also be applied to problems with a variety of boundary conditions even in three dimensions.  
\item[(ii)]
All the integrals encountered in the boundary integral equations (BIEs)  are evaluated indirectly by solving a single-fluid Stokes interface problem on staggered grid with a procedure of polynomial interpolation, thus the requirement of the analytical expressions
of the Green’s functions is successfully avoided. It is noteworthy that this technique is essentially different from that in \cite{layton2008efficient}. The integrals in \cite{layton2008efficient} are approximated directly using the trapezoid rule, so it may not be computationally efficient due to the fact that the integral kernel becomes singular in some circumstances. 
\item[(iii)]
 Another contribution of this work is the use of fast efficient solvers for the two linear systems resulting from BIEs and a modified MAC scheme respectively. The former one is solved iteratively by the GMRES method with a relatively small iteration number. And the presented numerical examples show that the number of GMRES iterations is actually independent of the system dimension and mesh parameter $h$. As for the latter linear system, since the jump conditions for the pressure $p$ and velocity $\vect u$ are decoupled for single-fluid Stokes interface problems, the coefficient matrix is identical to the case without interface, so that a conjugate gradient method incorporating the FFT-based solvers similar to that mentioned in \cite{dong2022second} for the velocity components can be adopted directly.
 
\item[(iv)]  Two types of numerical examples with different interface shapes are provided to show the efficiency and accuracy of the proposed method.   Examples with exact solution illustrate that second-order accuracy for velocity $\vect u$ and the pressure $p$ as well as the gradient of velocity in $\ell^2$-norm; for velocity $\vect u$ and its gradient in maximum norm can be achieved even with a relatively large ratio $\mu^+/\mu^-$. Examples with moving interface are also provided to demonstrate the stability and efficiency of the proposed method, where explicit time-stepping algorithms are developed for the motion of the interface.
\end{itemize}

The remaining part of the paper is organized as follows. In Section 2, the model of the steady incompressible two-fluid Stokes equations with interfaces is described, and the corresponding boundary integral formulation is presented in Section 3. Details of the evaluation of boundary or volume integral are presented in Section 4.
The numerical algorithm is summarized in the ensuing Section 5. In Section 6,  some numerical results are given to verify the accuracy and computation efficiency of the proposed approach. Some concluding remarks are made in Section 7. The equivalence between the integrals and the single-fluid Stokes interface problem is presented in Appendix.

%============================================================================== 
\section{The model problem}
\begin{figure}[ht!]
\centering
\includegraphics[width=0.5\textwidth]{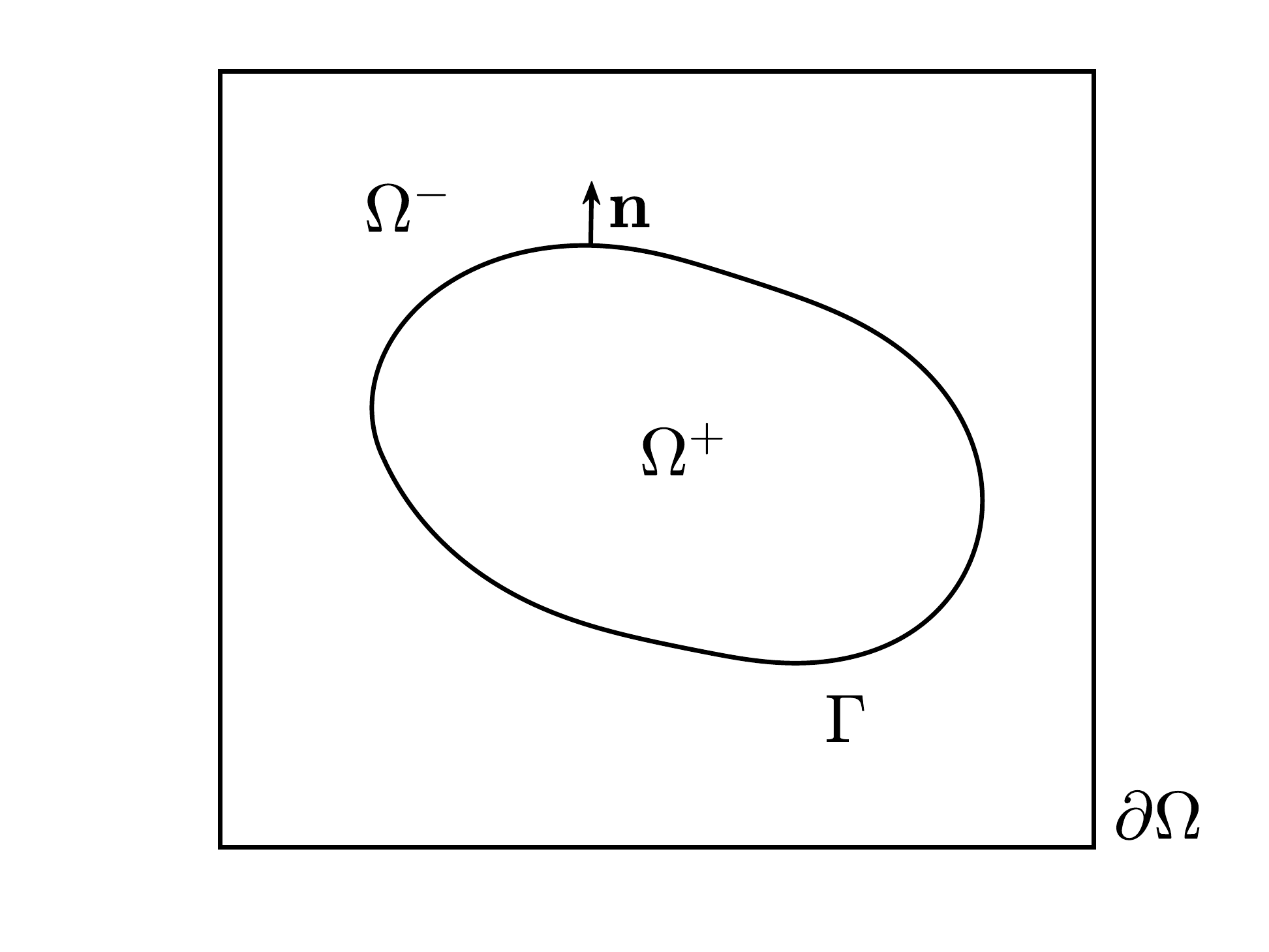}
\setlength{\abovecaptionskip}{-0.1cm}
 \caption{A sketch map for a rectangle domain $\Omega =\Omega^+\cup\Omega^-$ with interface $\Gamma$. }
 \label{OmegaDomain}
\end{figure}
Let $\Omega$ be a rectangle in two space dimensions, and $\Omega^+\subset\subset\Omega$ be a simply connected domain enclosed by a smooth interface $\Gamma$, illustrated in Fig. \ref{OmegaDomain}. Denote by $\Omega^-$ the complement of $\Omega^+$ in $\Omega$, satisfying $\Omega^-=\Omega\backslash\bar{\Omega}^+$. Consider the following Stokes interface problem
\begin{subequations}
\label{interfaceP}
\begin{align}
-\mu \Delta\vect  u+\nabla \widetilde p&=\widetilde{\vect  f},  \quad\;\;\hbox{in}\; \Omega^+\cup\Omega^-,\label{one}\\
\nabla\cdot \vect  u&=0,\quad\;\, \hbox{in}\; \Omega^+\cup\Omega^-,\label{two}\\
[\![ \vect  u]\!]&=\vect 0, \quad\;\,\hbox{on}\; \Gamma,\label{jump1}\\
[\![ \widetilde{\pmb \sigma}(\vect  u, \widetilde{p})\vect  n]\!]&=\vect  g, \,\quad\hbox{on}\; \Gamma,\label{jump2}\\
\vect  u&=\vect  u_b, \quad\hbox{on}\; \partial \Omega.\label{three}
\end{align}
\end{subequations}
Here, $\vect  u=(u^{(1)},u^{(2)})^T$, $\widetilde p$ and $\widetilde{\vect  f}=(\widetilde{f}^{(1)},\widetilde{f}^{(2)})^T$ represent the velocity, pressure and external force, respectively. The stress tensor is defined by 
$\widetilde{\pmb\sigma}(\vect  u, \widetilde p)
=-\widetilde p\vect  I+\mu(\nabla \vect  u+(\nabla \vect  u)^T)$, and $\vect  n$ is the unit outward normal vector on $\Gamma$ pointing from $\Omega^+$ to $\Omega^-$. The viscosity $\mu$ is assumed to be piecewise positive constant across the interface over the whole domain given by
\begin{equation*}
\mu=
\begin{cases}
\mu^+, \;\;\hbox{in}\; \Omega^+,\\
\mu^-,  \;\;\hbox{in}\;\Omega^-,
\end{cases}
\end{equation*}
where $\mu^+$ and $\mu^-$ are two positive constants.
The jump notation along the interface $\Gamma$ is denoted by 
$$[\![ \vect  v]\!](\vect x)\equiv\lim_{\epsilon\rightarrow 0}(\vect  v(\vect x-\epsilon\vect n) - \vect  v(\vect x+\epsilon\vect n)), \qquad \vect x\in \Gamma.$$
It is known that due to the incompressibility constraint,  the boundary data $\vect  u_b$ should satisfy the following compatibility condition
\begin{equation*}
\int_{\partial \Omega} \vect  u_b\cdot \vect  n_bds =0,
\end{equation*}
where $\vect  n_b$ is the outer unit  normal vector on $\partial\Omega$.

In order to facilitate the derivation of the integral equation, the model problem \eqref{interfaceP} is rewritten as
\begin{subequations}
\label{interfaceP2}
\begin{align}
-\Delta \vect  u+\nabla p&=\vect  f,  \,\;\;\;\hbox{in}\; \Omega^+\cup\Omega^-,\label{interfaceP2-1}\\
\nabla\cdot \vect  u&=0,\;\,\;\, \hbox{in}\; \Omega^+\cup\Omega^-,\label{interfaceP2-2}\\
[\![ \vect  u ]\!]&=\vect 0, \;\,\;\,\hbox{on}\; \Gamma,\label{interfaceP2-3}\\
[\![ \mu\pmb\sigma(\vect u, p)\vect  n]\!] &=\vect  g,\;\;\hbox{on}\; \Gamma,\label{interfaceP2-4}\\
\vect  u&=\vect  u_b, \;\, \hbox{on}\; \partial \Omega,\label{interfaceP2-5}
\end{align}
\end{subequations}
by introducing the scaled variables $p=\widetilde p/\mu, \; \vect  f=\widetilde {\vect  f}/\mu$  and 
$\pmb\sigma =\pmb\sigma(\vect  u, p)=-p\vect  I+(\nabla \vect  u+(\nabla \vect  u)^T)$.

Without loss of generality, this work assumes homogeneous boundary condition $\vect u = \vect 0$ in \eqref{interfaceP2-5}. Otherwise, the solution can be split into two components: $\vect u = \vect u_0 +\vect u_d, p = p_0 + p_d$. Here, $(\vect u_d, p_d)$ is the solution to the following Stokes problem
\begin{subequations}
\label{interfaceP-d}
\begin{align}
-\Delta\vect  u_d+\nabla p_d&= 0,  \,\quad\;\;\hbox{in}\; \Omega^+\cup\Omega^-,\label{one-d}\\
\nabla\cdot \vect  u_d&=0,\quad\;\;\, \hbox{in}\; \Omega^+\cup\Omega^-,\label{two-d}\\
\vect  u_d&=\vect  u_b, \quad\; \hbox{on}\; \partial \Omega,\label{three-d}
\end{align}
\end{subequations}
which can be solved by a classical numerical method since there is no discontinuity across the interface $\Gamma$. $(\vect u_0, p_0)$ is the solution to the Stokes interface problem with homogeneous Dirichlet boundary condition $\vect u = \vect 0$ in \eqref{interfaceP2-5} and the jump condition $[\![ \mu\pmb\sigma(\vect u, p)\vect  n]\!] = \vect g - [\![ \mu\pmb\sigma(\vect u_d, p_d)\vect  n]\!]$ instead of \eqref{interfaceP2-4}. Once $(\vect u_d, p_d)$ has been solved from \eqref{interfaceP-d}, the above jump condition is exactly known. So it is enough to consider problem $\eqref{interfaceP2}$ with $\vect u = \vect 0$ in \eqref{interfaceP2-5} . 
%%%%%%%%%%%%%%%%%%%%%%%%%%%%%%%%%%%%%
\section{The boundary integral formulation}
\label{sec;BIE}
In this section, the model problem \eqref{interfaceP2} will be reformulated into a boundary integral equation. 
Firstly, introduce Green's function pairs $(\vect  G_{\vect  v}(\vect  x, \vect  y), G_q(\vect  x, \vect  y))$ of the standard Stokes equations, satisfying 
\begin{equation}
\label{Greendefine}
\begin{split}
-\Delta\vect  G_{\vect  v}(\vect  x, \vect  y)+\nabla G_q(\vect  x, \vect  y) &= \vect  I\delta(\vect  x-\vect  y),\;\;\,\hbox{in}\;\Omega,\\
\nabla\cdot\vect  G_{\vect  v} (\vect  x, \vect  y)& = 0, \qquad\qquad\;\hbox{in}\;\Omega, \\
\vect  G_{\vect  v}(\vect  x, \vect  y)&=\vect  0,\qquad\qquad\;\hbox{on}\;\partial \Omega,
\end{split}
\end{equation}
 where $\vect  I$ is the unit matrix in $\mathbb{R}^d$ and $\delta(\vect  x-\vect  y) $ is the Dirac delta function.
 All differentiations are carried out with respect to the variable $\vect  x$. It is noted that the Green's pairs defined here are different from that in the free space.  Its expression is in general not analytically known, but its  existence is guaranteed.
 
To obtain the boundary integral equations,  rewrite function $\vect  u$, $p$ and $\vect  f$ respectively as 
\begin{equation*}
\vect  u=\begin{cases}
\vect  u^+,\;\hbox{in} \;\Omega^+,\\
\vect  u^-,\;\hbox{in}\;\Omega^-,
\end{cases}\quad\;
p=\begin{cases}
p^+,\;\hbox{in} \;\Omega^+,\\
p^-,\;\hbox{in}\;\Omega^-,
\end{cases}\quad\;
\vect  f=\begin{cases}
\vect  f^+,\;\hbox{in} \;\Omega^+,\\
\vect  f^-,\;\hbox{in}\;\Omega^-.
\end{cases}
\end{equation*}
 
 Applying the Green's second formula to $\vect  u^+$ on $\Omega^+$, one gets
\begin{equation}
\label{BIE-1}
\begin{split}
\int_{\Omega^+}\vect  G_{\vect  v}(\vect  x, \vect  y)\vect  f^+(\vect  y)d\vect  y&+\int_{\Gamma}\vect  G_{\vect  v}(\vect  x, \vect  y)T(\vect  u^+, p^+)ds_{\vect  y} \\
&- \int_{\Gamma}T(\vect  G_{\vect  v}, G_q) \vect  u^+(\vect  y)ds_{\vect  y} =\begin{cases}
\vect  u^+,\;\;\hbox{if}\; \vect  x\in\Omega^+,\\
\vect  0,\;\;\;\;\;\hbox{if}\;\vect  x\in\Omega^-.
\end{cases}
\end{split}
\end{equation}
Similarly,  applying the Green's second formula to $\vect  u^-$ on $\Omega^-$, one derives
\begin{equation}
\label{BIE-2}
\begin{split}
\int_{\Omega^-}\vect  G_{\vect  v}(\vect  x, \vect  y)\vect  f^-(\vect  y)d\vect  y &-\int_{\Gamma}\vect  G_{\vect  v}(\vect  x, \vect  y)T(\vect  u^-, p^-)ds_{\vect  y} \\
&+ \int_{\Gamma}T(\vect  G_{\vect  v}, G_q) \vect  u^-(\vect  y)ds_{\vect  y} 
= \begin{cases}
\vect  0,\;\;\;\;\,\hbox{if}\;\vect  x\in\Omega^+,\\
\vect  u^-,\;\;\hbox{if}\; \vect  x\in\Omega^-.
\end{cases}
\end{split}
\end{equation}
Adding \eqref{BIE-1} and \eqref{BIE-2}, one arrives at
\begin{equation}
\label{BIE-3}
\begin{split}
\int_{\Omega}\vect  G_{\vect  v}(\vect  x, \vect  y)\vect  f(\vect  y)d\vect  y &+\int_{\Gamma}\vect  G_{\vect  v}(\vect  x,\vect  y)\pmb \psi(\vect  y) ds_{\vect  y}\\
&-\int_{\Gamma}T(\vect  G_{\vect  v}, G_q) \pmb \varphi(\vect  y) ds_{\vect  y}=
\begin{cases}
\vect  u^+(\vect  x),\;\;\hbox{if}\;\vect  x\in\Omega^+,\\
\vect  u^-(\vect  x),\;\;\hbox{if}\; \vect  x\in\Omega^-,
\end{cases}
\end{split}
\end{equation}
with $\pmb \psi=T(\vect  u^+, p^+) -  T(\vect  u^-, p^-)=\pmb\sigma(\vect  u^+, p^+)\vect  n - \pmb\sigma(\vect  u^-, p^-)\vect  n$ and $\pmb \varphi = \vect  u^+-\vect  u^-$.
Substituting the relation $\nabla p=\Delta \vect  u + \vect  f$ into \eqref{BIE-1}-\eqref{BIE-2} and using the following identity 
\begin{equation}
\label{dd}
\begin{split}
\Delta_{\vect  x} T_{i,j}(\vect  G_{\vect  v}, G_q) &= -\delta_{i,j}\Delta_{\vect  x} G_q + \dfrac{\partial}{\partial y_i} \Delta_{\vect  x} (G_v)_j + \dfrac{\partial}{\partial y_j} \Delta_{\vect  x} (G_v)_i\\
&= -\delta_{i,j}\Delta_{\vect  x} G_q - \dfrac{\partial}{\partial x_i} \Delta_{\vect  x} (G_v)_j - \dfrac{\partial}{\partial x_j} \Delta_{\vect  x} (G_v)_i\\
&=- \frac{\partial^2 G_q}{\partial x_i\partial x_j} - \frac{\partial^2 G_q }{\partial x_j\partial x_i} \\
&= -2\frac{\partial^2G_q }{\partial x_i\partial x_j} 
=2\frac{\partial^2G_q }{\partial x_i\partial y_j},
\end{split}
\end{equation}   
one  gets
\begin{equation}
\label{BIE-4}
\begin{split}
\int_{\Omega^+} G_q(\vect  x, \vect  y) \cdot \vect  f^+ d\vect  y &+\int_{\Gamma} G_q(\vect  x, \vect  y)T(\vect  u^+, p^+) ds_{\vect  y} \\
&-2\int_{\Gamma} \dfrac{\partial G_q}{\partial \vect  n} \vect  u^+ ds_{\vect  y}=
\begin{cases}
p^+,\;\;\hbox{if}\;\vect  x\in\Omega^+,\\
0,\,\;\;\;\;\hbox{if}\; \vect  x\in\Omega^-,
\end{cases}
\end{split}
\end{equation}
and 
\begin{equation}
\label{BIE-5}
\begin{split}
\int_{\Omega^-} G_q(\vect  x, \vect  y) \cdot \vect  f^- d\vect  y &-\int_{\Gamma} G_q(\vect  x, \vect  y)T(\vect  u^-, p^-) ds_{\vect  y} \\
&+2\int_{\Gamma} \dfrac{\partial G_q}{\partial \vect  n} \vect  u^- ds_{\vect  y}=
\begin{cases}
0,\,\;\;\;\;\hbox{if}\;\vect  x\in\Omega^+,\\
p^-,\;\;\hbox{if}\; \vect  x\in\Omega^-.
\end{cases}
\end{split}
\end{equation}
Here $\delta_{i,j}$ is the Kronecker delta.
Adding \eqref{BIE-4} and \eqref{BIE-5}, one arrives
\begin{equation}
\label{BIE-6}
\begin{split}
\int_{\Omega}G_q(\vect  x, \vect  y)\cdot\vect  f(\vect  y)d\vect  y &+\int_{\Gamma}G_q(\vect  x,\vect  y)\cdot\pmb \psi(\vect  y) ds_{\vect  y}\\
&-2\int_{\Gamma}\dfrac{\partial G_q}{\partial \vect  n} \pmb \varphi(\vect  y) ds_{\vect  y}=
\begin{cases}
p^+(\vect  x),\;\;\hbox{if}\;\vect  x\in\Omega^+,\\
p^-(\vect  x),\;\;\hbox{if}\; \vect  x\in\Omega^-.
\end{cases}
\end{split}
\end{equation}
For the density function $\pmb\varphi, \pmb\psi$ and $\vect  f$, introduce the double layer boundary integrals $\mathcal{M}_{\vect v}\pmb\varphi$, $\mathcal{M}_q\pmb\varphi$, the single layer boundary integrals $\mathcal{L}_{\vect v}\pmb\psi$, $\mathcal{L}_q\pmb\psi$ and the volume integrals $\mathcal{G}_{\vect v}\vect  f$, $\mathcal{G}_{q}\vect  f$, which are given respectively by
\begin{equation*}
\begin{split}
(\mathcal{M}_{\vect  v}\pmb\varphi)(\vect  x) &= \int_{\Gamma} T(\vect  G_{\vect  v}, G_q) \pmb\varphi(\vect  y) ds_{\vect  y}, \quad\;\;
(\mathcal{M}_q\pmb\varphi)(\vect  x) = 2\int_{\Gamma} \dfrac{\partial G_q(\vect  x, \vect  y)}{\partial \vect  n_{\vect  y}} \pmb\varphi(\vect  y)ds_{\vect  y},\\
(\mathcal{L}_{\vect  v}\pmb\psi)(\vect  x) &= \int_{\Gamma} \vect  G_{\vect  v}(\vect  x, \vect  y) \pmb\psi(\vect  y) ds_{\vect  y},\qquad\;\;\,
(\mathcal{L}_q\pmb\psi)(\vect  x) = \int_{\Gamma} G_q(\vect  x, \vect  y)\cdot \pmb\psi(\vect  y) ds_{\vect  y},\\
(\mathcal{G}_{\vect  v}\vect  f)(\vect  x) &= \int_{\Omega}\vect  G_{\vect  v}(\vect  x, \vect  y)\vect  f(\vect  y)d\vect  y, \qquad\quad\;\;
(\mathcal{G}_q\vect  f)(\vect  x) = \int_{\Omega} G_q(\vect  x, \vect  y)\cdot\vect  f(\vect  y)d\vect  y,
\end{split}
\end{equation*}
and the adjoint double layer $\mathcal{D}\pmb\varphi$ and hyper-singular operator $\mathcal{M}^*\pmb\psi$, which are  defined respectively by
\begin{equation*}
\begin{split}
(\mathcal{D}\pmb\varphi)(\vect  x) &= -T_x(\mathcal{M}_{\vect  v}\pmb\varphi, \mathcal{M}_q\pmb\varphi) 
= -\Big( -\mathcal{M}_q\pmb\varphi  +\nabla (\mathcal{M}_{\vect  v}\pmb\varphi )
+ \big(\nabla (\mathcal{M}_{\vect  v}\pmb\varphi) \big)^T\Big)\vect  n,\\
(\mathcal{M}^*\pmb\psi)(\vect  x) &= \int_{\Gamma} T_y(\vect  G_{\vect  v}, G_q)\pmb\psi(\vect  y)ds_{\vect  y}.
\end{split}
\end{equation*}
It is remarked that the subscript $\vect  x$ or $\vect  y$ in operator $T$ represents that the differentiations are with respect to the corresponding  variable $\vect  x$ or $\vect  y$. 

By the continuity properties of the hyper-singular boundary and volume integrals, and the discontinuity properties of the adjoint double layer potential, one has
\begin{equation}
\label{BIE-7}
\begin{split}
\pmb\sigma(\vect  u^+, p^+)\vect n\big|_{\Gamma} &= \dfrac{1}{2} \pmb\psi + \mathcal{M}^*\pmb\psi + T(\mathcal{G}_{\vect  v}\vect  f, \mathcal{G}_q\vect  f) +\mathcal{D}\pmb\varphi,\\
\pmb\sigma(\vect  u^-, p^-)\vect n\big|_{\Gamma} &= -\dfrac{1}{2} \pmb\psi + \mathcal{M}^*\pmb\psi + T(\mathcal{G}_{\vect  v}\vect  f, \mathcal{G}_q\vect  f) +\mathcal{D}\pmb\varphi.
\end{split}
\end{equation}

Plugging \eqref{BIE-7} into the interface conditions \eqref{interfaceP2-4}, one can obtain the following BIE
\begin{equation}
\label{BIEinterface}
\frac{1}{2}\pmb\psi+\gamma\mathcal{M}^*\pmb\psi = \hat{\vect  g}
-\gamma T(\mathcal{G}_{\vect  v}\vect  f, \mathcal{G}_q\vect  f)-\gamma\mathcal{D}\pmb\varphi,  \; \hbox{on}\;\Gamma,
\end{equation}
with $\gamma=(\mu^+-\mu^-)/(\mu^++\mu^-)$ and $\hat{\vect  g}=\vect  g/(\mu^++\mu^-)$. 
The integral equation \eqref{BIEinterface} is a Fredholm integral equation of the second kind, which is uniquely solvable for the unknown density $\pmb\psi$ \cite{kress1989linear}.  Once a numerical solution $\pmb \psi$ is solved from BIE \eqref{BIEinterface}, the approximation of $(\vect u, p)$ to the interface problem \eqref{interfaceP2} can be further computed by the following representation formulas
\begin{equation}
\label{Ss}
\begin{split}
\vect u &=  \mathcal{G}_{\vect v}\vect f(\vect x) +\mathcal{L}_{\vect v}\pmb\psi(\vect x) -\mathcal{M}_{\vect v}\pmb\varphi(\vect x),\\
p &=  \mathcal{G}_{q}\vect f(\vect x) +\mathcal{L}_{q}\pmb\psi(\vect x) -\mathcal{M}_{q}\pmb\varphi(\vect x).
\end{split}
\end{equation}

At the end of this chapter, two remarks about the BIEs are given: 
\begin{itemize}
\item[1).]
As the boundary integral operator $\mathcal{M}^*$ is a compact operator and its spectrum lies on the interval $[-1/2, 1/2)$ \cite{kress1989linear, Beale2004grid}, the spectral radius of the operator $\mathbb{A} = \frac{1}{2}\mathbb{I} + \gamma \mathcal{M}^*$ falls on the interval $(0,1]$. Thus the minimum eigenvalue of the operator $\mathbb{A}$ is strictly greater than a positive number and the corresponding discrete linear system is well conditioned. Therefore, the integral equation \eqref{BIEinterface} can be solved iteratively using a Krylov subspace method. 
In this work, the generalized minimal residual (GMRES) method \cite{saad1986gmres,Saad1993GMRES} is adopted, which converges to a prescribed tolerance in a fixed number of steps with any initial guess $\pmb\psi_0$ in the solution space.   It is remarkable that the condition number of the system grows in the case of $\mu^+\gg\mu^-$, because the minimum eigenvalue of the operator $\mathbb{A}$ becomes close to $0$, causing the fact that the number of GMRES iteration increases as illustrated in Section \ref{sec;num}. 
\item[2).]
Owing to the conceivable unavailability of Green's function pairs $({\vect G}_{\vect v}, G_q)$ defined in the bounded domain $\Omega$, the main challenge in solving the BIE \eqref{BIEinterface} is the evaluation of the boundary and volume integrals.  This work does not discretize the corresponding integrals directly by numerical quadratures. Instead, it follows the lines of the KFBI method \cite{ying2007kernel} and makes the evaluation indirect by a Cartesian grid based method in two steps: solving interface problem and extracting boundary data, which will be discussed in details in the next section. 
\end{itemize}
%%%%%%%%%%%%%%%%%%%%%%%%%%%%%%%%%%%%%

%%%%%%%%%%%%%%%%%%%%%
\section{Evaluation of boundary or volume integral}
As analytical expressions of Green's pairs $(G_{\vect v}, G_q)$ are not easily available, it is difficult to directly evaluate the integrals encountered in this work. Interestingly, the KFBI method makes the evaluation indirectly by a Cartesian grid based approach and analytical expressions of Green's functions are no longer needed,  which is a major difference from the traditional boundary integral methods.
As discussed in  \cite{ying2007kernel,ying2013kernel,ying2013fast,ying2014kernel}, each boundary or volume integral appeared in the BIE has an equivalent but simple interface problem in the sense that  the partial differential equations (PDEs) only involve uniformly continuous coefficients.  Therefore, evaluating a boundary or volume integral primarily consists of two essential components:
\begin{itemize}
\item [i).] Solve an equivalent but simple interface problem in $\Omega$, including discretizing the corresponding interface problem with an appropriate scheme, modifying the established linear system and solving the resulting linear system on uniform Cartesian grids with a fast Fourier transform (FFT)-based conjugate gradient (CG) method. 
\item [ii).]  Extract the boundary values of the integrals on $\Gamma$ by quadratic polynomial interpolation with the discrete numerical solution.
\end{itemize}

\subsection{Equivalent interface problems}
This subsection lists the equivalent interface problems for the volume, single layer boundary and double layer boundary integrals, respectively. It is assumed that any vector function $\vect v(\vect x)$ or scalar function $q(\vect x)$ and their partial derivatives involved in the following are at least piecewise smooth, with its potential discontinuity only existing on the interface $\Gamma$.

The volume integrals $\vect v(\vect x) = \mathcal{G}_{\vect v}\vect f(\vect x)$ and $q=\mathcal{G}_q\vect f(\vect x)$ are computed by solving the following simple interface problem with discontinuous inhomogeneous source
\begin{equation}
\label{v-i}
\begin{split}
-\Delta \vect  v+\nabla q&=
\begin{cases}
\vect  f^+, \;\hbox{in}\; \Omega^+,\\
\vect  f^-,  \;\hbox{in}\; \Omega^-,\\
\end{cases}\\
\nabla\cdot \vect  v&=0,\;\; \hbox{in}\; \Omega^+\cup\Omega^-,\\
[\![ \vect  v ]\!]&=\vect 0, \;\;\hbox{on}\; \Gamma,\\
[\![ \pmb \sigma(\vect v, q)\vect  n]\!]&=\vect 0, \;\; \hbox{on}\; \Gamma,\\
\vect  v&=\vect 0,\;\;\hbox{on}\;\partial \Omega.
\end{split}
\end{equation}
 It is noted that the interface conditions above indicate the continuous property of the volume potential $\vect v$ 
 as well as  its traction $\pmb \sigma(\vect v, q)\vect  n$.
 
The double layer boundary integrals $\vect v(\vect x) = -\mathcal{M}_{\vect v}\pmb \varphi(\vect x)$ and 
$q = -\mathcal{M}_q\pmb\varphi(\vect x)$ are calculated by solving a  simple homogeneous interface problem
\begin{equation}
\label{d-i}
\begin{split}
-\Delta \vect  v+\nabla q&=
\vect 0,\;\; \hbox{in}\; \Omega^+\cup\Omega^-,\\
\nabla\cdot \vect  v&=0,\;\; \hbox{in}\; \Omega^+\cup\Omega^-,\\
[\![ \vect  v ]\!]&=\pmb\varphi, \,\;\hbox{on}\; \Gamma,\\
[\![ \pmb \sigma(\vect v, q)\vect  n]\!]&=\vect 0, \;\; \hbox{on}\; \Gamma,\\
\vect  v&=\vect 0,\;\;\hbox{on}\;\partial \Omega.
\end{split}
\end{equation} 
The discontinuity properties of the double layer potential imply that 
\begin{equation*}
\begin{split}
\vect  v^+ &= \frac{1}{2}\pmb\varphi - \mathcal{M}_{\vect v}\pmb\varphi, \;\;\; \;{\rm on}\;\Gamma,\\
\vect  v^- &= -\frac{1}{2}\pmb\varphi - \mathcal{M}_{\vect v}\pmb\varphi, \;\; {\rm on}\;\Gamma.
\end{split}
\end{equation*}

The single layer boundary integrals $\vect v(\vect x) = \mathcal{L}_{\vect v}\pmb \psi(\vect x)$ and 
$q = \mathcal{L}_q\pmb\psi(\vect x)$  are  also evaluated by solving a simple homogeneous interface problem
\begin{equation}
\label{s-i}
\begin{split}
-\Delta \vect  v+\nabla q&=\vect 0, \;\,\;\hbox{in}\; \Omega^+\cup\Omega^-,\\
\nabla\cdot \vect  v&=0,\;\; \,\hbox{in}\; \Omega^+\cup\Omega^-,\\
[\![ \vect  v ]\!]&=\vect 0, \,\;\;\hbox{on}\; \Gamma,\\
[\![ \pmb \sigma(\vect v, q)\vect  n]\!]&=\pmb\psi, \;\; \hbox{on}\; \Gamma,\\
\vect  v&=\vect 0,\,\;\;\hbox{on}\;\partial \Omega.
\end{split}
\end{equation} 
Moreover, the jump relations above can be rewritten as 
\begin{equation*}
\begin{split}
\pmb\sigma(\vect  v^+, q^+)\vect  n &= \frac{1}{2}\pmb\psi+ \mathcal{M}^*\pmb\psi, \;\; \;\;{\rm on}\;\Gamma,\\
\pmb\sigma(\vect  v^-, q^-)\vect  n&= -\frac{1}{2}\pmb\psi+ \mathcal{M}^*\pmb\psi, \;\; {\rm on}\;\Gamma.
\end{split}
\end{equation*}

Note that different integrals correspond to different terms on the right hand side. Moreover, based on the assumptions on the domain $\Omega$ and interface $\Gamma$, each interface problem above has a unique solution for sufficiently smooth functions $\vect f, \pmb\varphi$ and $\pmb\psi$. By the linearity of the problems, the solution \eqref{Ss} to the interface problem \eqref{interfaceP2} is the sum of the solutions to the previous three interface problems \eqref{v-i}-\eqref{s-i}, which can be presented and solved in a unified framework,
\begin{equation}
\label{i}
\begin{split}
-\Delta \vect  v+\nabla q&=\vect f, \;\,\;\hbox{in}\; \Omega^+\cup\Omega^-,\\
\nabla\cdot \vect  v&=0,\;\; \,\hbox{in}\; \Omega^+\cup\Omega^-,\\
[\![ \vect  v ]\!]&=\pmb\varphi, \,\;\;\hbox{on}\; \Gamma,\\
[\![ \pmb \sigma(\vect v, q)\vect  n]\!]&=\pmb\psi, \;\; \hbox{on}\; \Gamma,\\
\vect  v&=\vect 0,\,\;\;\hbox{on}\;\partial \Omega.
\end{split}
\end{equation} 
Therefore, to evaluate the volume or boundary integrals, one can turn to solving the interface problem \eqref{i} with the source term $\vect f$ or the jumps $\pmb\varphi, \pmb\psi$, which are determined by the corresponding value given in equations \eqref{v-i}-\eqref{s-i}, respectively.
It is worth pointing out that the interface problem \eqref{i} does not have any discontinuous coefficient at all, which is much simpler to solve than the interface problem \eqref{interfaceP2}.

Further more, proof of the equivalences between the interface problems with continuous viscosity  and the  integrals here is similar to  
that for elliptic problem \cite{ying2007kernel, ying2014kernel}, but need some extra techniques, which will be addressed in the Appendix.

\subsection{A Cartesian grid-based MAC scheme for simple interface problem}
\label{sub;MAC}
Since the interface problem \eqref{i} is a simple interface problem in the sense that the PDEs involve only constant viscosity, it can be solved with various existing methods in the literature, one can refer to \cite{cortez2001method,fogelson1986numerical,leveque1997immersed,mayo1992implicit,tu1992stability,dong2022second} and the references therein.  Very recently, a modified finite difference MAC scheme with second-order accuracy is proposed in \cite{dong2022second}. This subsection will give a brief review of the implementation. 

To this end, the computational domain $\Omega$ is partitioned into $N\times N$ small rectangles of the same shape. 
The discretization is carried out on a standard MAC staggered grid with mesh size $h = x_{i+1} - x_i = y_{j+1} - y_j$. With the MAC mesh,  
the velocity component $v^{(1)}$ is located at the vertical edges of a cell $(ih, (j-1/2)h)$, with $i\in\{1, 2, \cdots, N-1\}$ and $j\in\{1,2,\cdots, N\}$;  the velocity component $v^{(2)}$ is located at the horizontal edges of a cell $((i-1/2)h,jh)$, with $i\in\{1, 2, \cdots, N\}$ and $j\in\{1,2,\cdots, N-1\}$; the pressure field $q$ is defined at the cell center $((i-1/2)h, (j-1/2)h)$, with $i\in\{1, 2, \cdots, N\}$ and $j\in\{1,2,\cdots, N\}$.  The modified MAC scheme reads as 
\begin{equation}
\begin{split}
-\Delta_h v^{(1)}_{i,j-\frac{1}{2}} +\delta_{h,1}^+ q_{i-\frac{1}{2},j-\frac{1}{2}}
 &=f^{(1)}_{i,j-\frac{1}{2}} + C\{\Delta v^{(1)}\}_{i,j-\frac{1}{2}} +C\{q_x\}_{i,j-\frac{1}{2}},\\[8pt]
-\Delta_h v^{(2)}_{i-\frac{1}{2},j} +\delta_{h,2}^+ q_{i-\frac{1}{2},j-\frac{1}{2}}
 &=f^{(2)}_{i-\frac{1}{2},j} + C\{\Delta v^{(2)}\}_{i-\frac{1}{2},j} +C\{q_y\}_{i-\frac{1}{2},j},\\[8pt]
 \delta_{h,1}^-v^{(1)}_{i,j-\frac{1}{2}} +\delta_{h,2}^-v^{(2)}_{i-\frac{1}{2},j}& = C\{v_x^{(1)}\}_{i-\frac{1}{2},j-\frac{1}{2}} + C\{v_y^{(2)}\}_{i-\frac{1}{2},j-\frac{1}{2}}
\end{split}
\end{equation}
where 
\begin{equation*}
\begin{split}
\delta_{h,1}^+\,v_{l,m}&=h^{-1}\left(v_{l+1,m}
-v_{l,m}\right), \;\;\quad
\delta_{h,1}^-\,v_{l,m}=h^{-1}\left(v_{l,m}
-v_{l-1,m}\right), \\[4pt]
\delta_{h,2}^+\,v_{l,m}&=h^{-1}\left(v_{l,m+1}
-v_{l,m}\right),\;\;\quad
\delta_{h,2}^-\,v_{l,m}=h^{-1}\left(v_{l,m}
-v_{l,m-1}\right), \\
\Delta_{h}v_{l,m} &= \delta_{h, 1}^+\,\delta_{h, 1}^-\,v_{l,m}
+\delta_{h, 2}^+\,\delta_{h, 2}^-\,v_{l,m},
\end{split}
\end{equation*}
with $l$ 
taking values $i, i-\frac{1}{2}$ for integer $i$ and $m$ taking values 
$j, j-\frac{1}{2}$ for integer $j$. 

It is noted that the above MAC scheme has been modified due to the jump conditions across the interface $\Gamma$. The correction terms are non-zero only at irregular grid nodes near the interface and only appear on the right hand of the linear system. Therefore, the coefficient matrix of the discrete system is the same as that resulting from the discretization of the Stokes problems without an interface. Thus the CG method together with an FFT-based Poisson solver can be applied directly.  Moreover, these correction terms will improve the truncation errors near the interface to at least first-order accuracy. 
Interestingly, one lower order of truncation error at the interface will not affect overall second order accuracy of the solution, which is verified simultaneously by the numerical experiment and theoretical analysis in \cite{dong2022second}.  
In addition, the correction terms need the information about the jumps of the solution and their derivatives, which can be evaluated from the original interface conditions \eqref{interfaceP2-4}.
For more details about the calculation of jump conditions, the derivation of corrections terms, and the rigorous proof of second-order accuracy, the reader is referred to \cite{dong2022second}.

\subsection{Interpolation for integral values on the interface}
As we can see,  the approximation solution $(\vect v_h, q_h)$ to the interface problem \eqref{i}, the equivalent boundary integral or volume integral, is obtained at the staggered grid,  while the approximation of the corresponding boundary or volume integral needed in \eqref{BIEinterface} should be evaluated at discretization points of the interface. 
Thus, a polynomial interpolation should be designed to extract the boundary value $\vect v_h$ and its flux $\pmb\sigma(\vect v_h, q_h)$ at any given discretization points on the interface.  Here, the interpolation technique is completely the same as the treatment of elliptic interface problem in \cite{ying2007kernel,ying2014kernel} except that the velocity $\vect v$ and pressure $q$ are treated separately on the respective grid positions.  

\begin{figure}[ht!]
\centering
\includegraphics[width=0.23\textwidth]{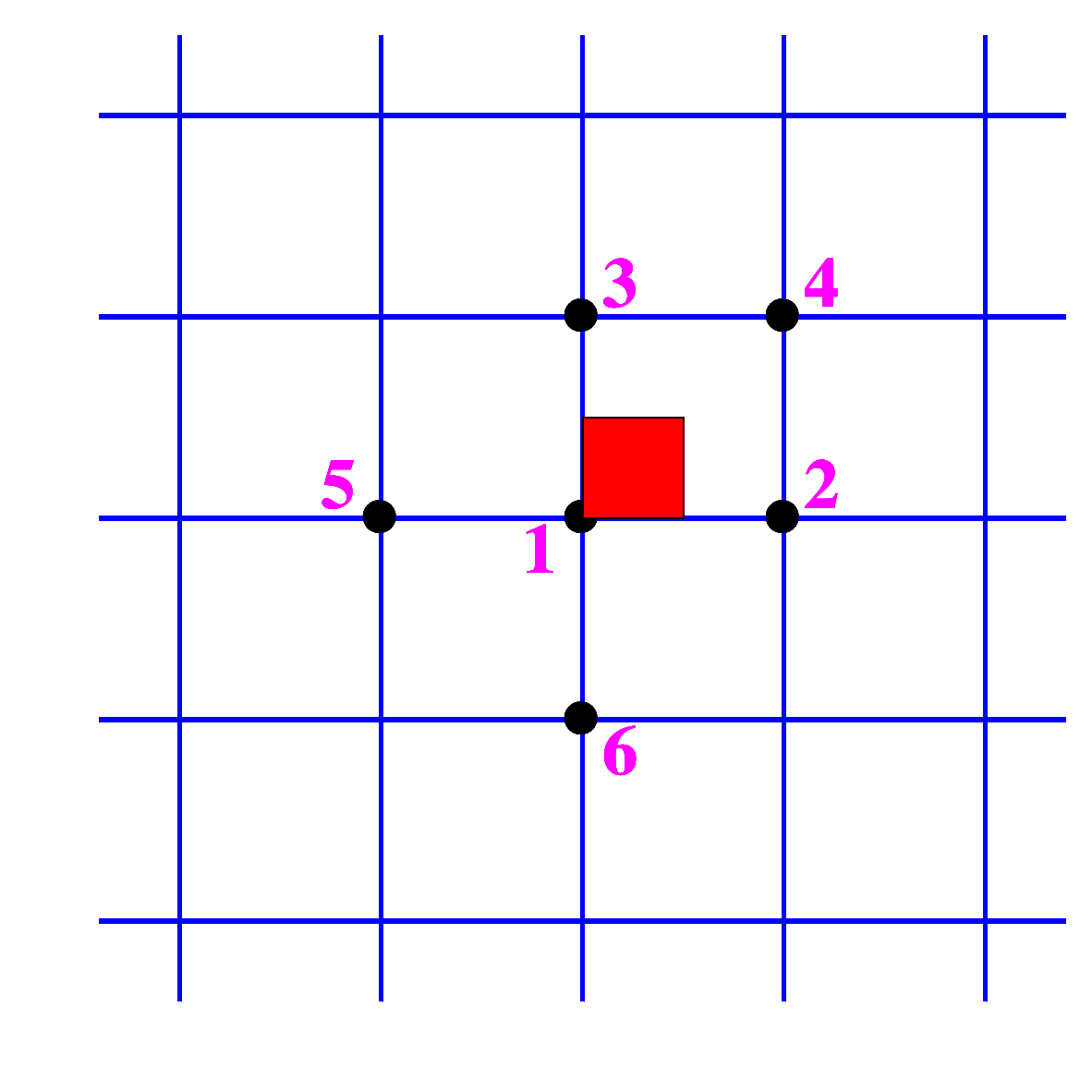}
\includegraphics[width=0.23\textwidth]{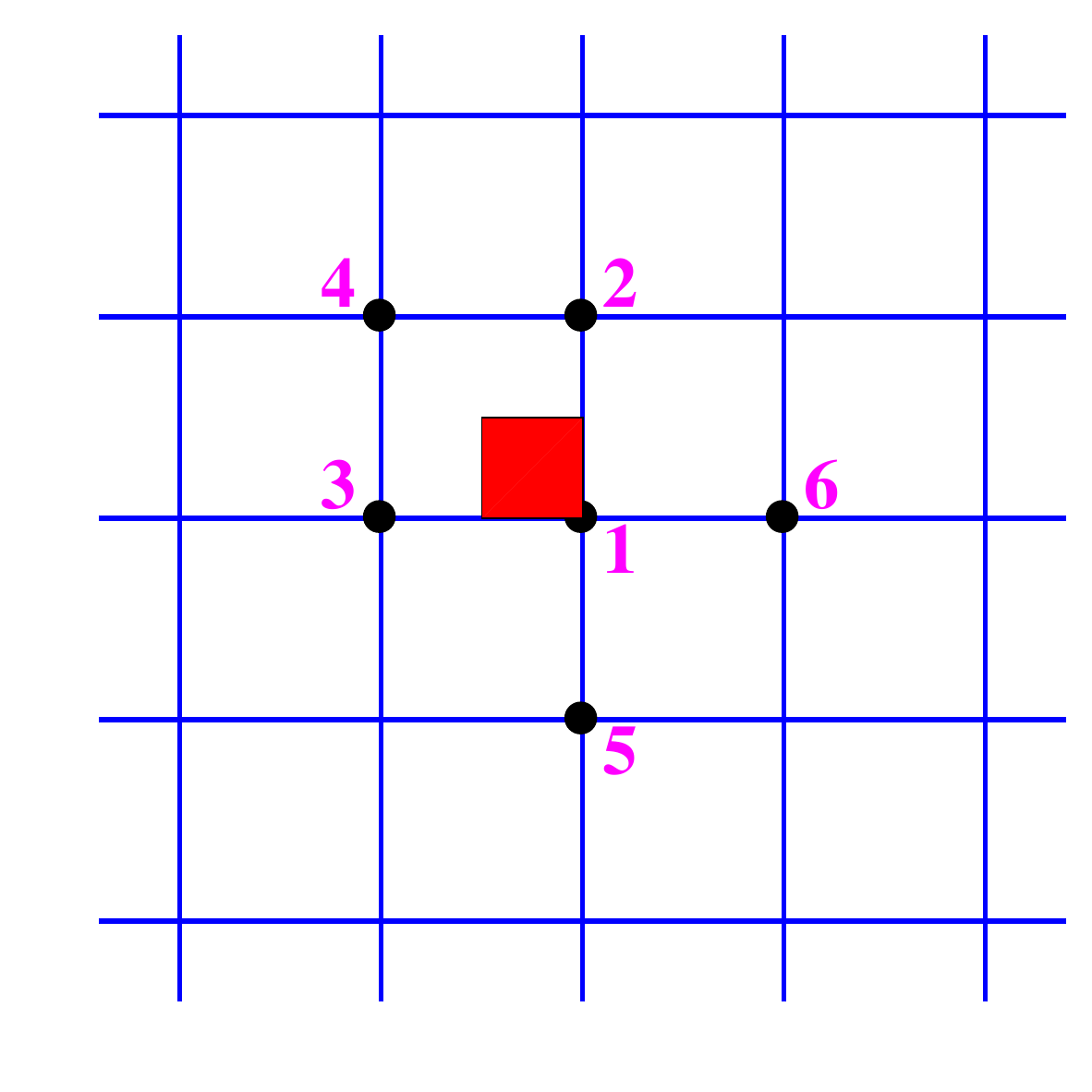}
\includegraphics[width=0.23\textwidth]{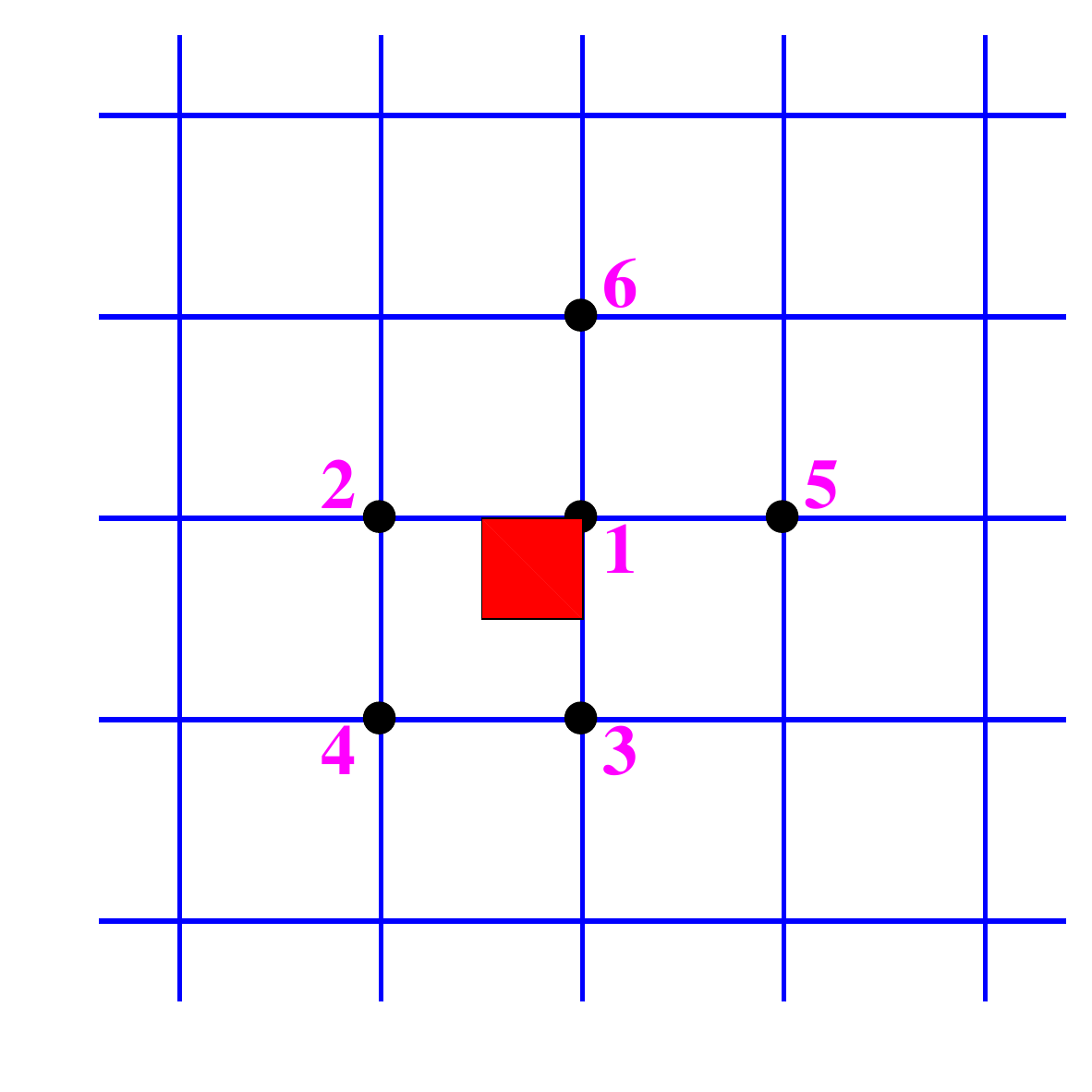}
\includegraphics[width=0.23\textwidth]{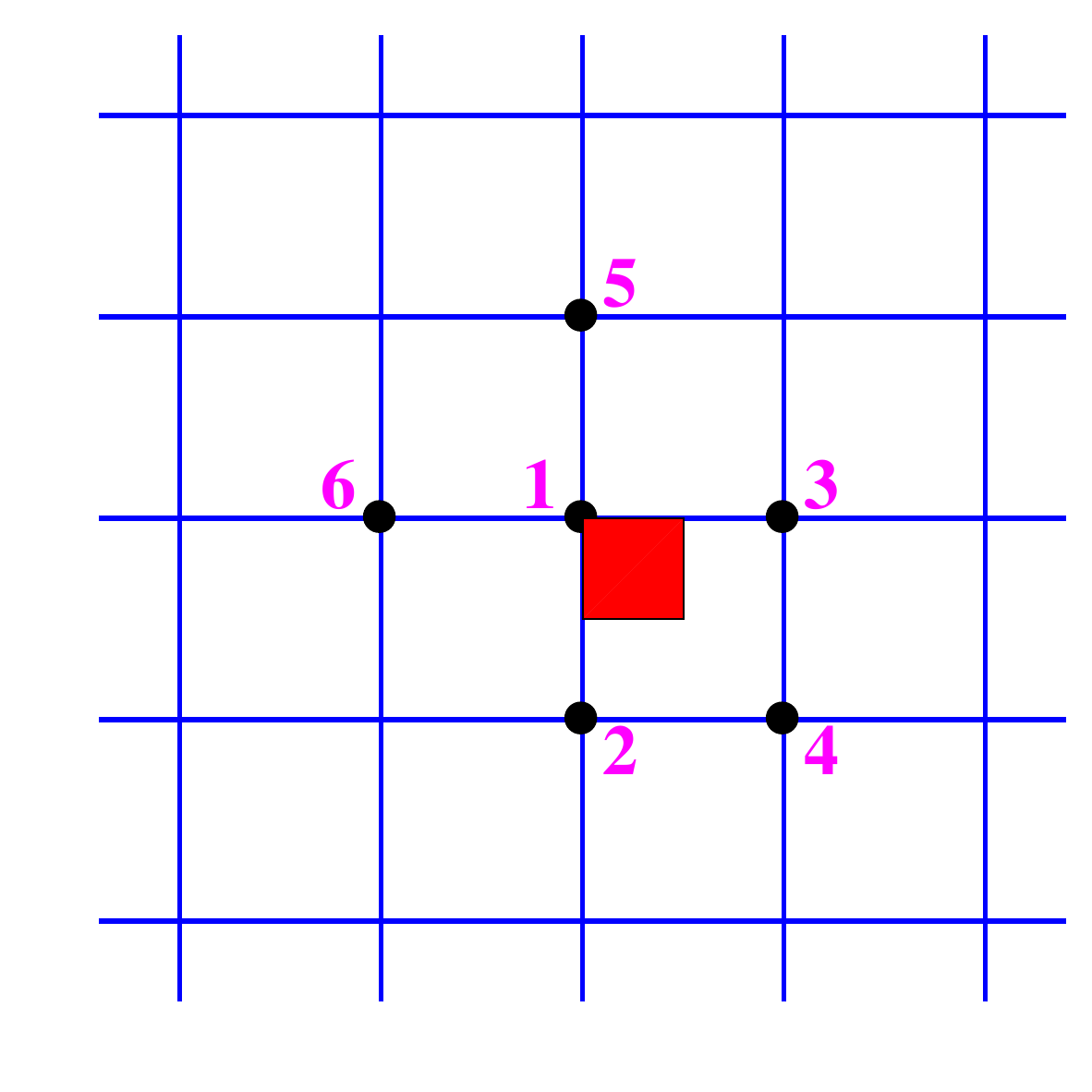}
\setlength{\abovecaptionskip}{-0.1cm}
 \caption{Six grid nodes for the quadratic Lagrange interpolation at points in the red region, which is located at four different corner of a grid cell respectively.  }
 \label{sixgrid}
\end{figure}

Assume that the second-order finite difference solution $\vect v_h$ to the interface problem \eqref{i} maintains the same piecewise smooth property as $\vect v$ does. For a point $\vect  x$ on $\Gamma$, which is just located in the red square region,  Fig. \ref{sixgrid} shows the interpolation stencil, where six grid nodes $\vect z_k$ are involved. 
Taylor expansion of the approximation solution $\vect v_h$  around $\vect  x\in \Gamma$  for each interpolation point $\vect z_k(k = 1,2,\cdots, 6)$, gives us
\begin{subequations}
\begin{align}
\nonumber
\vect v_h(\vect  z_k)&=\vect v_h^+(\vect  x)+\xi_k\frac{\partial \vect v_h^+(\vect  x)}{\partial x}
+\eta_k\frac{\partial \vect v_h^+(\vect  x)}{\partial y}
+\frac{1}{2}\xi_k^2\frac{\partial^2 \vect v_h^+(\vect  x)}{\partial x^2}\\[4pt]
\label{Taylor-1}
&+\xi_k\eta_k\frac{\partial^2 \vect v_h^+(\vect  x)}{\partial x\partial y}
+\frac{1}{2}\eta_k^2\frac{\partial^2 \vect v_h^+(\vect  x)}{\partial y^2}
+\mathcal{O}(|\vect  z_k-\vect  x|^3), \quad\;\hbox{if}\;  \vect   z_k\in \Omega^+,\\[4pt]
\nonumber
\vect v_h(\vect  z_k)&=\vect v_h^-(\vect  x)+\xi_k\frac{\partial \vect v_h^-(\vect  x)}{\partial x}
+\eta_k\frac{\partial \vect v_h^-(\vect  x)}{\partial y}
+\frac{1}{2}\xi_k^2\frac{\partial^2 \vect v_h^-(\vect  x)}{\partial x^2}\\[4pt]
\label{Taylor-2}
&+\xi_k\eta_k\frac{\partial^2 \vect v_h^-(\vect  x)}{\partial x\partial y}
+\frac{1}{2}\eta_k^2\frac{\partial^2 \vect v_h^-(\vect  x)}{\partial y^2}
+\mathcal{O}(|\vect  z_k-\vect  x|^3), \quad\;\hbox{if}\;  \vect   z_k\in \Omega^-,
\end{align}
\end{subequations}
Here, $(\xi_k, \eta_k)^T\equiv \vect  z_k-\vect  x$. It is remarked that the interpolation node $\vect z_k=(x_i,y_{j-\frac{1}{2}})$ for component $v^{(1)}_h$ and $\vect z_k=(x_{i-\frac{1}{2}}, y_j)$ for component $v^{(2)}_h$, but the above Taylor expansions  do not make distinguish between them only for the convenience of writing. Let  
\begin{equation*}
\vect J_k\equiv [\![ \vect v_h ]\!]+\xi_k\Big[\!\!\Big[ \dfrac{\partial \vect v_h}{\partial x}\Big]\!\!\Big]
+\eta_k\Big[\!\!\Big[ \dfrac{\partial \vect v_h}{\partial y}\Big]\!\!\Big]+\frac{1}{2}\xi_k^2\Big[\!\!\Big[ \dfrac{\partial^2 \vect v_h}{\partial x^2} \Big]\!\!\Big]+
\xi_k\eta_k\Big[\!\!\Big[ \dfrac{\partial^2 \vect v_h}{\partial x\partial y} \Big]\!\!\Big]+\frac{1}{2}\eta_k^2\Big[\!\!\Big[ \dfrac{\partial^2 \vect v_h}{\partial y^2} \Big]\!\!\Big], \quad \hbox{if}\; \vect z_k\in \Omega^-,
\end{equation*}
where $ [\![ \vect v_h ]\!]$ and its partial derivatives are computable according to subsection \ref{sub;MAC} with the assumption that $\vect v_h$ satisfies the same jump conditions as $\vect v$ does. Then expansion formula \eqref{Taylor-1}-\eqref{Taylor-2} can be rewritten as 
\begin{subequations}
\begin{align}
\label{inter-1}
\vect v_h(\vect  z_k)&=\vect v_h^+(\vect  x)+\Big(\xi_k\frac{\partial }{\partial x}
+\eta_k\frac{\partial }{\partial y}\Big) \vect v_h^+(\vect  x)
+\frac{1}{2}\Big(\xi_k^2\frac{\partial^2}{\partial x^2}
+2\xi_k\eta_k\frac{\partial^2 }{\partial x\partial y}
+\eta_k^2\frac{\partial^2}{\partial y^2}\Big) \vect v_h^+(\vect  x), 
\quad\hbox{if}\;  \vect   z_k\in \Omega^+,\\[4pt]
\label{inter-2}
\vect v_h(\vect  z_k)+\vect J_k&=\vect v_h^+(\vect  x)+\Big(\xi_k\frac{\partial }{\partial x}
+\eta_k\frac{\partial }{\partial y}\Big) \vect v_h^+(\vect  x)
+\frac{1}{2}\Big(\xi_k^2\frac{\partial^2}{\partial x^2}
+2\xi_k\eta_k\frac{\partial^2 }{\partial x\partial y}
+\eta_k^2\frac{\partial^2}{\partial y^2}\Big) \vect v_h^+(\vect  x), 
\quad\hbox{if}\;  \vect   z_k\in \Omega^-,
\end{align}
\end{subequations}
with the third-order term $\mathcal{O}(|\vect  z_k-\vect  x|^3)$ omitted. 
Note that the coefficient matrix of \eqref{inter-1} and \eqref{inter-2} is dependent of the mesh parameter $h$ due to the dependency of $\xi_k$ and $\eta_k$ on $h$, thus solution of this linear system may involve large computer round-off errors.   In the practical calculation, the above equations are transformed  into the following formula
 \begin{equation}
 \label{Tplus1}
\begin{split}
\vect v_h(\vect z_k)&=\vect {\bar v} ^+ +\alpha_k \vect {\bar v}^+_x + \beta_k \vect {\bar v}^+_y
+\dfrac{1}{2}\alpha_k^2  {\bar v}^+_{xx} + \alpha_k\beta_k \vect {\bar v}^+_{xy} 
+\dfrac{1}{2}\beta_k^2 \vect {\bar v}^+_{yy}, \; \;\;\hbox{if}\;  \vect   z_k\in \Omega^+,\\
\vect v_h(\vect z_k)+\vect J_k&=\vect {\bar v} ^+ +\alpha_k \vect {\bar v}^+_x + \beta_k \vect {\bar v}^+_y
+\dfrac{1}{2}\alpha_k^2  {\bar v}^+_{xx} + \alpha_k\beta_k \vect {\bar v}^+_{xy} 
+\dfrac{1}{2}\beta_k^2 \vect {\bar v}^+_{yy}, \; \;\;\hbox{if}\;  \vect   z_k\in \Omega^-,
\end{split}
\end{equation}
with $\alpha_k=\xi_k/h, \beta_k=\eta_k/h$ being new coefficients and 
\begin{equation*}
\begin{split}
&\vect {\bar v}^+=\vect v_h^+(\vect x), 
\qquad\qquad\;\; \vect {\bar v}^+_x=h\frac{\partial \vect v_h^+(\vect x)}{\partial x}, 
\qquad\quad \vect {\bar v}^+_y=h\frac{\partial \vect v_h^+(\vect x)}{\partial y}, \\
&\vect {\bar v}^+_{xx}=h^2\frac{\partial^2 \vect v_h^+(\vect x)}{\partial x^2},
\qquad \vect {\bar v}^+_{xy}=h^2\frac{\partial^2 \vect v_h^+(\vect x)}{\partial x\partial y}, 
\qquad \vect {\bar v}^+_{yy}=h^2\frac{\partial^2 \vect v_h^+(\vect x)}{\partial y^2},
\end{split}
\end{equation*}
being new quantities. 
The coefficient matrix of the rescaled system \eqref{Tplus1} is independent of $h$ and
 the choice of the six interpolation nodes ${\vect  z_k}$ selected as stated in Fig. \ref{sixgrid} guarantees that the coefficient matrix is always invertible. The polynomial interpolation scheme is stable for a general shape of interface $\Gamma$. This implies that the quadratic interpolation for  the limit values of the approximate solution  $\vect  v_h$ on $\Gamma$ is uniquely determined by \eqref{Tplus1}.
 
Since the stress tensor $\pmb\sigma(\vect v_h, q_h)$ defined on the interface $\Gamma$ does not involve the derivatives of $q_h$, linear interpolation is enough in the calculation of boundary value. As done before,  one can get
\begin{equation}
\label{Tplus2}
\begin{split}
q_h(\vect z_k)&={\bar q} ^+ +\alpha_k {\bar q}^+_x + \beta_k {\bar q}^+_y, \; \;\;\hbox{if}\;  \vect   z_k\in \Omega^+,\\[4pt]
q_h(\vect z_k)+J_k&={\bar q} ^+ +\alpha_k {\bar q}^+_x + \beta_k {\bar q}^+_y, \; \;\;\hbox{if}\;  \vect   z_k\in \Omega^-,
\end{split}
\end{equation} 
with 
\begin{equation*}
\begin{split}
&{\bar q}^+=q_h^+(\vect x), 
\qquad\bar q^+_x=h\frac{\partial q_h^+(\vect x)}{\partial x}, 
\qquad\bar q^+_y=h\frac{\partial q_h^+(\vect x)}{\partial y},
\end{split}
\end{equation*}
and
\begin{equation*}
J_k\equiv [\![ q_h ]\!]+\xi_k\Big[\!\!\Big[ \dfrac{\partial q_h}{\partial x}\Big]\!\!\Big]
+\eta_k\Big[\!\!\Big[ \dfrac{\partial q_h}{\partial y}\Big]\!\!\Big], \quad \hbox{if}\; \vect z_k\in \Omega^-,
\end{equation*}
where $ [\![ q_h ]\!]$ and its partial derivatives are also computable. Here, $(\xi_k, \eta_k)^T\equiv \vect  z_k-\vect  x$, with the interpolation node $\vect z_k=(x_{i-\frac{1}{2}},y_{j-\frac{1}{2}})$ for component $q_h$. The  linear polynomial interpolation defined in \eqref{Tplus2} is also guaranteed to uniquely exist.
After the linear systems \eqref{Tplus1} and \eqref{Tplus2} are solved, simple manipulation leads to the stress tensor as desired, which are then applied to discretize the corresponding boundary  or volume  integral encountered in section \ref{sec;BIE}.

\section{Algorithm summary}
A KFBI method is presented to solve the Stokes interface problem with piecewise constant viscosity. The major difference between this method and classical BIE method is that the associated volume and boundary integrals are computed as limit values of the Cartesian grid-based approximate solutions. In this section, the algorithm is summarized as follows:
\begin{itemize}
\item[(1)] Some preparatory work, such as partitioning the interface $\Gamma$  into a set of quasi-uniformly spaced nodes $\mathcal{N}_{\Gamma}$ and computing normals, tangents, curvatures of the interface $\Gamma$; discretizing the computational domain $\Omega$ by a staggered grid; identitying the interior and irregular grid points; finding the intersection points of the interface with the staggered grid lines.
\item[(2)] Evaluation of the boundary or volume integral on the interface $\Gamma$, which is concretely translated into the following two points:
  \begin{description}
   \item[$\bullet$] solve the interface problem \eqref{i}: discretize the interface problem with second-order accuracy MAC scheme; compute jumps of partial derivatives at intersection points;  compute the correction terms and modify the right hand side of the discrete interface problem at irregular gird nodes; solve the modified linear system with a  CG method incorporated with an FFT-based Poisson solvers. 
   \item[$\bullet$] extract the boundary data:
   compute jumps of partial derivatives at  the discrete points $\mathcal{N}_{\Gamma}$ of the interface $\Gamma$; 
   use the computed jumps to calculate the limiting values of the numerical solution or its stress tensor at the interface points to obtain the boundary integral or volume integral by a polynomial  interpolation.
\end{description}
\item[(3)] The GMRES iteration:
\begin{itemize}
   \item[(3.1).]  evaluate the volume potential and the double layer potential boundary data  with steps $(2)$;
   \item[(3.2).]  choose an initial guess $\pmb \psi_0$  to start the GMRES iteration and set up a stopping criterion $\epsilon$.
   \item[(3.3).]  evaluate the single layer potential boundary data $\pmb\sigma(\vect v_h, q_h)$ with steps $(2)$;
   \item[(3.4).]   update the unknown discrete density $\pmb \psi_{\nu+1}$ by the GMRES iteration.
   \item[(3.5).] go to step $(3.3)$ until the residual is small enough in some norm.
\end{itemize}
\item[(4)] Superposition of the solutions: once the GMRES iteration converges, add the volume potential, the double layer potential and the single layer potential together.
\end{itemize}

\section{Numerical examples}
\label{sec;num}
In this section, some numerical examples with different coefficients and interface geometries are presented to validate the accuracy and efficiency of the proposed method. 
To do so, the scaled discrete $l^2$-norms are defined  respectively by
$$\|e_{\vect  u}\|=\frac{\|\vect  u-\vect  u_h\|}{\|\vect  u\|},\qquad 
\| e_{\vect u}\|_1 = \frac{\|\vect  u-\vect  u_h\|_1}{\|\vect  u\|_1},\qquad
\|e_{p}\|=\frac{\|p-p_h\|}{\|p\|},$$
where 
\begin{equation*}
\begin{split}
\| \vect u \|^2 \equiv  h^2 \sum\limits_{i=1}^{N-1}\sum\limits_{j=1}^N \Big(u_{i,j-\frac{1}{2}}^{(1)}\Big)^2 +h^2 \sum\limits_{i=1}^{N}\sum\limits_{j=1}^{N-1} \Big(u_{i-\frac{1}{2},j}^{(2)}\Big)^2, \qquad
\|p\|^2 = h^2 \sum\limits_{i=1}^N\sum\limits_{j=1}^N \big(p_{i-\frac{1}{2},j-\frac{1}{2}}\big)^2,
\end{split}
\end{equation*}
and
\begin{equation*}
\begin{split}
\| \vect u \|_1^2 &\equiv 
h^2 \sum\limits_{i=1}^{N} \sum\limits_{j=1}^{N} \Big(\delta_{h,1}^-v^{(1)}_{i-\frac{1}{2}, j-\frac{1}{2}}\Big)^2 + 
h^2\sum_{i=1}^{N-1}\sum_{j=0}^{N}\rho_{j}^y \Big(\delta_{h,2}^-v^{(1)}_{i, j}\Big)^2 \\ 
&\quad+
h^2\sum_{i=0}^{N}\sum_{j=1}^{N-1}\rho_{i}^x \Big(\delta_{h,1}^-v^{(2)}_{i, j}\Big)^2 + 
h^2 \sum\limits_{i=1}^{N} \sum\limits_{j=1}^{N} \Big(\delta_{h,2}^-v^{(2)}_{i-\frac{1}{2}, j-\frac{1}{2}}\Big)^2,
\end{split}
\end{equation*}
with $\rho_{0}^x=\rho_{N}^x=\dfrac{1}{2}, \rho_{i}^x=1$ when 
$i=1, 2, \cdots, N-1$, and  $\rho_{0}^y=\rho_{N}^y=\dfrac{1}{2},\rho_{j}^y=1$ 
when $j=1, 2, \cdots, N-1$. 

The scaled discrete maximum norm  respectively by  
$$\|e_{\vect  u}\|_{\infty}=\frac{\|\vect  u-\vect  u_h\|_{\infty}}{\|\vect  u\|_{\infty}},\qquad
\|e_{\vect  u}\|_{1,\infty}=\frac{\|\vect  u-\vect  u_h\|_{1,\infty}}{\|\vect  u\|_{1,\infty}},\qquad
\|e_{p}\|_{\infty}=\frac{\|p-p_h\|_{\infty}}{\|p\|_{\infty}},$$
where
\begin{equation*}
\|\vect u\|_{\infty} \equiv \frac{1}{2}\Big(\max\limits_{1\leq i\leq N-1, 1\leq j \leq N} u^{(1)}_{i,j-\frac{1}{2}} + \max\limits_{1\leq i\leq N, 1\leq j \leq N-1} u^{(2)}_{i-\frac{1}{2},j}\Big),  \qquad
\|p\|_{\infty} \equiv \max\limits_{1\leq i\leq N-1, 1\leq j \leq N-1}p_{i-\frac{1}{2},j-\frac{1}{2}}, 
\end{equation*}
and
\begin{equation*}
\begin{split}
\|\vect u\|_{1,\infty} \equiv \frac{1}{4}\Big(&\max\limits_{1\leq i\leq N, 1\leq j \leq N} \delta_{h,1}^-u^{(1)}_{i-\frac{1}{2},j-\frac{1}{2}} + \max\limits_{1\leq i\leq N-1, 0\leq j \leq N} \delta_{h,2}^-u^{(1)}_{i,j}  \\[6pt]
+ &\max\limits_{0\leq i\leq N, 1\leq j \leq N -1} \delta_{h,1}^-u^{(2)}_{i,j}  +  \max\limits_{1\leq i\leq N, 1\leq j \leq N} \delta_{h,2}^-u^{(2)}_{i-\frac{1}{2},j-\frac{1}{2}} \Big).
\end{split}
\end{equation*}
In all numerical experiments, the standard GMRES iterative method is used to solve the discrete linear system corresponding to the BIE \eqref{BIEinterface}. 
Moreover,  the iterated unknown density $\pmb\psi$ is initialized  with an initial guess whose entries equal the values of $\vect g$ at the corresponding discretization points of the interface and the GMRES iteration stops when the maximum norm of the residual is less than the tolerance $\epsilon=10^{-8}$.

The proposed algorithm was implemented in custom codes written in the C++ computer language. All calculation reported below were performed on an iMAC with $3.8$GHz Inter Core $i7$.
%:::::::::::::::::::::::::::::::::::::::::::::::::::::::::::::::::::::::::::::::::::::::::::::::::::::::::::::::::::::::::::::::::::::::::::::::::::::::

\subsection{Test examples with exact solution} 
As we know, it is quite challenging to construct the exact solutions to incompressible Stokes flow with an interface. 
Here, only two examples with the exact solution known are performed to demonstrate the accuracy 
for the velocity $\vect  u$, its gradient and the pressure $p$.  In order to numerically illustrate the robustness of our algorithm to large jumps across the interface $\Gamma$,  six different cases with large viscosity  will be considered:
\begin{itemize}
\item[I)]  \; $\mu^+ = 1,  \;\qquad\mu^-=10$;
\item[II)] \; $\mu^+ = 1,  \;\qquad\mu^-=100$;  
\item[III)] \; $\mu^+ = 1, \;\qquad\mu^-=1000$;  \\
\item[IV)]\,\;$\mu^+ = 10, \quad\quad\mu^-=1$; 
\item[V)]\,\;$\mu^+ = 100, \,\;\quad\mu^-=1$; 
\item[VI)] \,\;$\mu^+ = 1000, \quad\mu^-=1$.\\ 
\end{itemize}
In the following two examples, a grid refinement analysis is performed in Tables \eqref{tabl2-e1-1} -\eqref{tabmax-e2-2}. The second column is the number of grid lines in both $x-$ and $y-$ directions. The third column is the maximum error (or $l^2$-error) of the velocity $\vect u$ while the fourth column is the approximate convergence order. The fifth column is the  maximum error (or $l^2$-error) of  the gradient of  velocity $\vect u$ while the sixth column is the approximate convergence order.  The seventh column is the maximum error (or $l^2$-error) of  the pressure $p$ and the eighth column is the corresponding approximate convergence order.

%====================================================================
\begin{table}[h]
\caption{The GMRES iteration number of solving BIEs in {\em Example} 1 and {\em Example} 2.}
\begin{center}
\begin{tabular}{|c|c|c|c|c|c|c|c|}
\hline
&\diagbox{N}{case} &  $I)$  &  $II)$  & $III)$  & $IV)$ & $V)$ & $VI)$\\
\hline
{\em Example} 1     & 128 & 9 & 10 & 10 & 12 & 12 & 22\\ 
                   & 256 & 8 & 9 & 9 & 11 & 11 & 22\\ 
                   & 512 & 7 & 8 & 8 & 9 & 9 &14\\
                   &1024 & 6 & 6 & 7 & 8 & 8 & 14\\
                   & 2048 & 6 & 6 & 6 & 7 & 7 & 13\\
\hline
{\em Example} 2     & 128 & 10 & 12 & 12 & 14 & 19 & 23\\ 
                   & 256 & 10 & 11 & 11 & 12 & 19 & 24\\ 
                   & 512 & 9 & 11 & 11 & 11 & 14 &16\\
                   &1024 & 9 & 9 & 9 & 10 & 13 & 16\\
                   & 2048 &9 & 9 & 9 & 10 & 13 & 15\\
\hline
\end{tabular}
\end{center}
\label{tab-GMRES}
\end{table}

%========================================================================
{\em Example} 1. In this example,  the interface is a circle $\{(x, y)|x^2+y^2=1\}$, which is located at the center of the box $\Omega=(-2.0, 2.0)^2$. The exact solution and the pressure are given by 
\begin{equation*}
\begin{split}
u^{(1)}(x, y)&=\begin{cases} 
\dfrac{y}{r}-\dfrac{3y}{4}, \qquad\;\;\; x^2+y^2>1,\\[8pt]
\dfrac{y}{4}(x^2+y^2),\quad\;x^2+y^2\leq 1,
\end{cases}\\[4pt]
u^{(2)}(x, y)&=\begin{cases} 
-\dfrac{x}{r}+\dfrac{x}{4}(3+x^2), \;\;\; x^2+y^2>1,\\[8pt]
-\dfrac{xy^2}{4},\qquad\qquad\;\;\;\;\; \,x^2+y^2\leq 1,
\end{cases}\\[4pt]
p(x, y)&=\begin{cases} 
(-\dfrac{3}{4}x^3+\dfrac{3}{8}x)y, \quad\; x^2+y^2>1,\\[8pt]
5.0,\qquad\qquad\quad\;\;\;\; \,x^2+y^2\leq 1.
\end{cases}
\end{split}
\end{equation*}
Here $r = \sqrt{x^2+y^2}$.  It is easy to check that $\nabla\cdot\vect u =0$.  The pressure and the gradient of velocity are discontinuous across the interface. The external force and boundary data can be evaluated from the exact solution.

\begin{figure}[h!]
\centering
\subfigure[$u^{(1)}$]{
\includegraphics[width=0.31\textwidth]{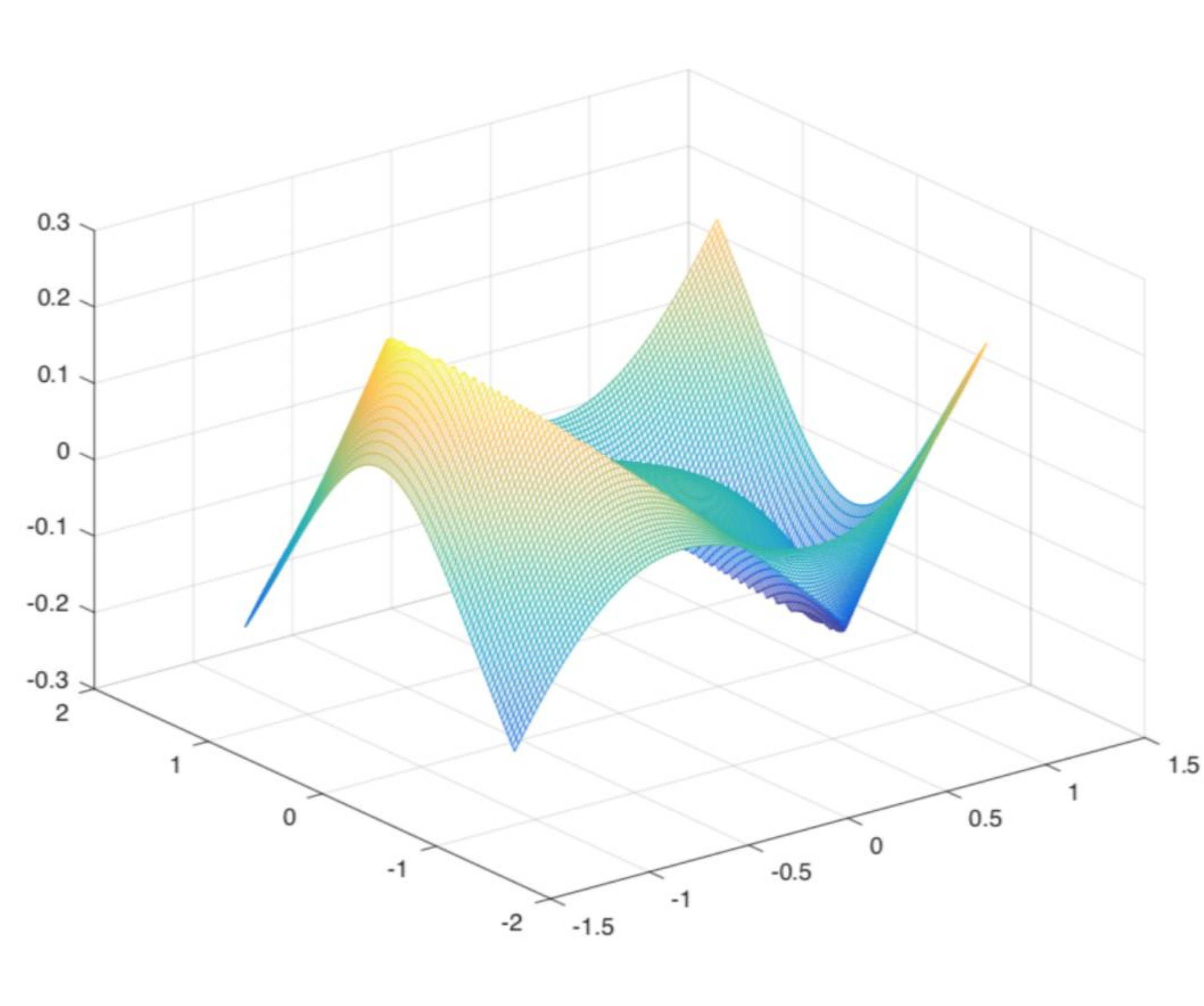}
}
\subfigure[$u^{(2)}$]{
\includegraphics[width=0.31\textwidth]{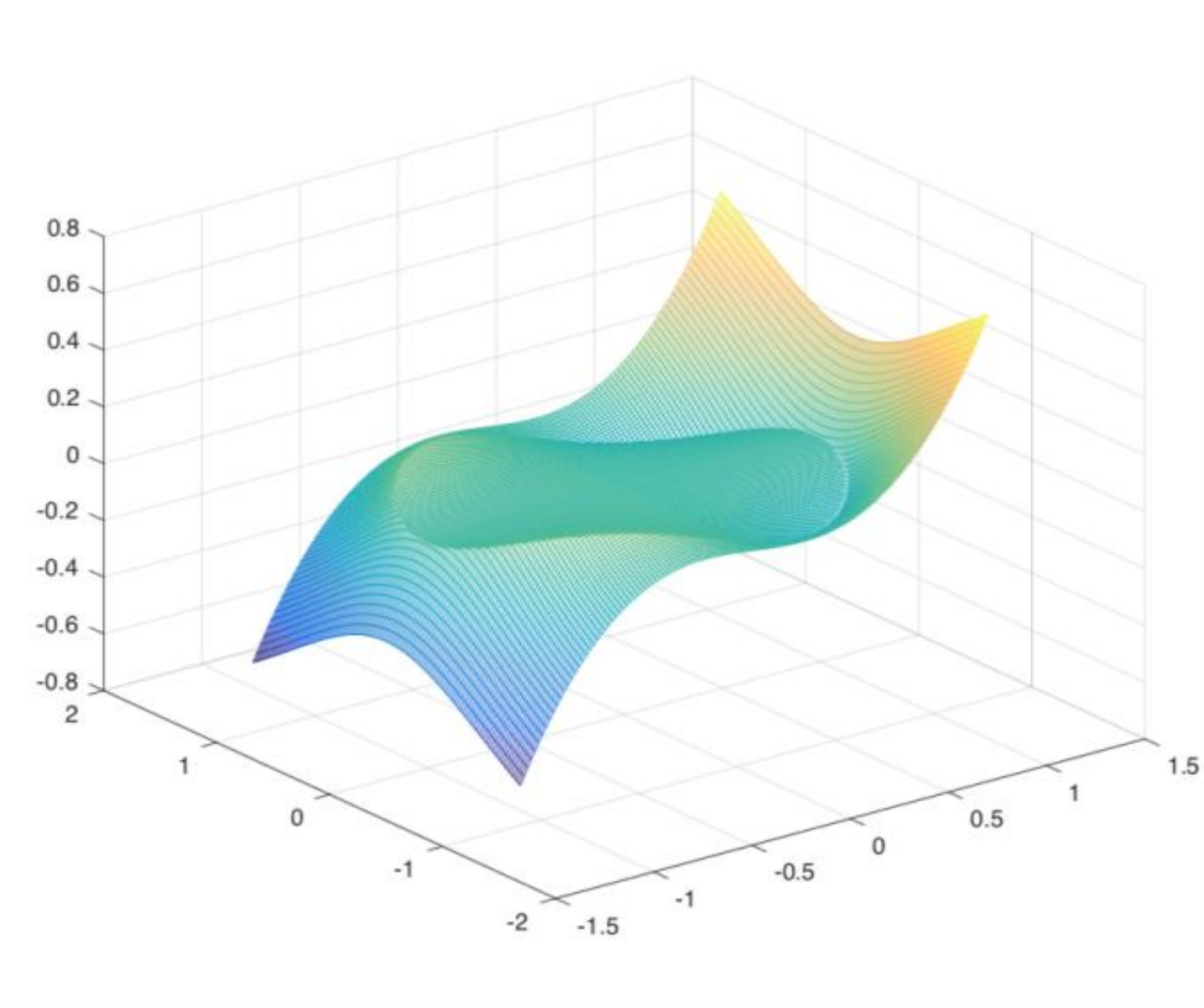}
}
\subfigure[$p$]{
\includegraphics[width=0.31\textwidth]{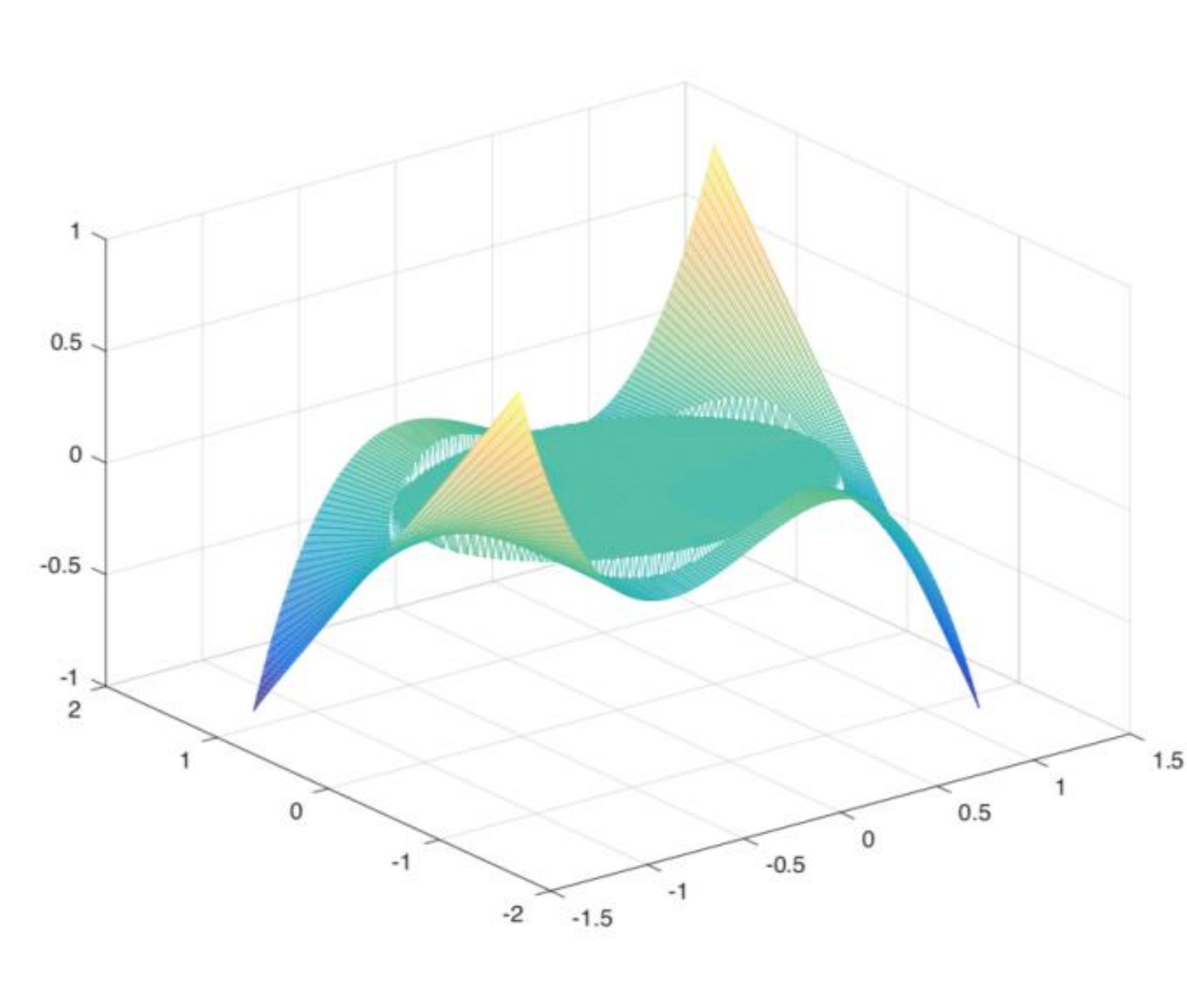}
}
\setlength{\abovecaptionskip}{-0.0cm}
\setlength{\belowcaptionskip}{-0.0cm}
 \caption{The solution plots of the $x$-component of velocity field $u^{(1)}$, the $y$-component of velocity field $u^{(2)}$ and the pressure $p$ in {\em Example} 1 with viscosity $\mu^+ = 10$ and $\mu^-=1$ on a $128\times 128$ 
 grid.}
 \label{v-p-e1}
\end{figure}

Normalized errors for the velocity $\vect u$, the gradient of the velocity $\nabla \vect u$ and the pressure $p$ in the discrete $l^2$-norms are shown in Table \ref{tabl2-e1-1} and \ref{tabl2-e1-2}.  Normalized errors in maximum norms are shown in Table \ref{tabmax-e1-1} and \ref{tabmax-e1-2}.  It can be seen that  the velocity and its derivatives are all second order accurate in both the discrete  $l^2$-norm and the discrete maximum norm, 
 and the pressure is second accurate in  $l^2$-norm but first-order accurate in maximum norm for different coefficient ratios.  Furthermore, the convergence rate of the proposed method is independent of the jump in viscosity.
The GMRES iteration number is shown in Table \ref{tab-GMRES}. One can observe that more iterations are needed when $\mu^+/\mu^-\gg 1$ than the other case. 
However, a limited number of iterations are needed and the number of iterations is almost independent of the mesh size and is insensitive to the ratio of coefficients.
The solution plots are depicted in Fig. \ref{v-p-e1}. 

\begin{table}[h]
\caption{$l^2$-error and its convergence rates of {\em Example} 1 with viscosity $\mu^+\ll\mu^-$.}
\begin{center}
\vspace{-0.1cm}
\begin{tabular}{|c|c|c|c|c|c|c|c|c|}
\hline
case & N   & $\|e_{\vect  u}\|$
& order & $\|e_{\vect  u}\|_1$ & order &$\|e_p\|$ &order\\
\hline 
   & $128$      &  3.9927e-4  &     -       & 1.4338e-4    &   -           & 3.9806e-3  &      -  \\
   & $256$      &  1.0002e-4  &  1.9971    & 3.5969e-5   & 1.9950     & 1.1097e-3   & 1.8428 \\
I) & $512$       &  2.4051e-5  &  2.0561   & 9.0023e-6   & 1.9984     & 3.0392e-4  & 1.8684 \\
   & $1024$     &  5.9084e-6  &  2.0253  & 2.2506e-6   & 2.0000	    & 8.1987e-5   & 1.8902\\
   & $2048$     &  1.4535e-6  &  2.0232  & 5.6253e-7    &2.0003     & 2.1913e-5   &   1.9036\\		
\hline
    & $128 $      & 4.0590e-4  &    -        & 1.4341e-4    &    -        &  3.9829e-2   &     -  \\
    & $256$      &  1.0141e-4   & 2.0009   & 3.5971e-5    & 1.9952   &  1.1098e-2    & 1.8435 \\
II) & $512$       &  2.4352e-5  & 2.0581   & 9.0026e-6    & 1.9984   &  3.0395e-3  & 1.8684 \\
    & $1024$     &  5.9744e-6  & 2.0272   & 2.2507e-6   & 2.0000   &   8.1988e-4  & 1.8903\\
    & $2048$    &  1.4682e-6   &  2.0459  & 5.6253e-7   & 2.0004  &   2.1907e-4  &   1.9040\\		
\hline
    & $128$     &  4.0659e-4 &         -    &  1.4341e-4   &    -       &  3.9832e-1   &    - \\
    & $256$    &  1.0155e-4   &  2.0014   &  3.5972e-5  & 1.9952   &  1.1098e-2   & 1.8436 \\
III) & $512$     &  2.4384e-5  &  2.0582  &  9.0026e-6  & 1.9985   &  3.0395e-2  & 1.8684 \\
    & $1024$   &  5.9809e-6  &  2.0275  &  2.2507e-6   & 2.0000  &  8.1993e-3   & 1.8903\\
    & $2048$   &  1.4698e-6  &  2.0247  &  5.6254e-7   &2.0003  	&  2.1907e-3    &  1.9041\\
\hline
\end{tabular}
\end{center}
\label{tabl2-e1-1}
\end{table}

\begin{table}[h]
\caption{Maximum error and  its  convergence rates of {\em Example} 1 with viscosity $\mu^+\ll\mu^-$. }
\begin{center}
\vspace{-0.1cm}
\begin{tabular}{|c|c|c|c|c|c|c|c|}
\hline
case & N &$\|e_{\vect  u}\|_{\infty}$ 
 & order & $\|e_{\vect  u}\|_{1, \infty}$& order& $\|e_p\|_{\infty}$& order\\
 \hline
   & $128$     & 3.2974e-4  & -          & 1.4172e-4   & -          & 2.6250e-2 &   - \\
   & $256$    & 8.2791e-5   &  1.9938  & 3.5084e-5  & 2.0142  & 1.3053e-2 & 1.0079 \\
I) & $512$     & 1.9607e-5   &  2.0781  & 8.7282e-6  & 2.0071  & 6.5088e-3 & 1.0039\\
   & $1024$   & 4.7533e-6  &  2.0444  & 2.1768e-6  & 2.0035  & 3.2499e-3 & 1.0020\\
   & $2048$   & 1.1511e-6    & 2.0459   & 5.4353e-7  & 2.0018  & 1.6238e-3 &  1.0010 \\	
\hline
    & $128$    & 3.4676e-4   & -          & 1.4172e-4  &   -         & 2.6258e-1   &   - \\
    & $256$    & 8.6864e-5  & 1.9971    & 3.5084e-5 &2.0142   & 1.3055e-1    & 1.0082 \\
II) & $512$     & 2.0560e-5  & 2.0789   & 8.7282e-6 & 2.0071   & 6.5093e-2   & 1.0040 \\
   & $1024$    & 4.9761e-6   & 2.0468   & 2.1768e-6 & 2.0035  & 3.2500e-2   & 1.0021\\
   & $2048$   & 1.2027e-6   & 2.0487   & 5.4353e-7  & 2.0018   & 1.6238e-3   &  1.0011 \\	
\hline
    & $128$     & 3.4853e-4 & -        & 1.4272e-4      &   -         & 2.6258e+0 &   - \\
    & $256$    & 8.7292e-5 & 1.9978  & 3.5084e-5     &2.0243    & 1.3055e+0 & 1.0082 \\
III) & $512$    & 2.0660e-5 & 2.0790  & 8.7282e-6     & 2.0071    & 6.5090e-1 & 1.0041 \\
    & $1024$   & 4.9995e-6 & 2.0470  & 2.1768e-6     & 2.0035   & 3.2500e-1 & 1.0020\\
    & $2048$  & 1.2080e-6   & 2.0492   & 5.4353e-7  & 2.0018 	& 1.6238e-1 &  1.0011 \\	
\hline
\end{tabular}
\end{center}
\label{tabmax-e1-1}
\end{table}

\begin{table}[h]
\caption{$l^2$-error and convergence rates of {\em Example} 1 with viscosity $\mu^+\gg\mu^-$.}
\begin{center}
\vspace{-0.1cm}
\begin{tabular}{|c|c|c|c|c|c|c|c|}
\hline
case & N  & $\|e_{\vect  u}\|$
& order & $\|e_{\vect  u}\|_1$ & order &$\|e_p\|$ &order\\
\hline
     & $128$      & 9.0937e-4  &  -      & 3.2781e-4  &-        & 5.9439e-4 &      - \\
     & $256$    & 1.8310e-4  &  2.3122 & 5.1541e-5  & 2.6691  & 1.5366e-4 & 1.9517 \\
IV) & $512 $    & 4.3072e-5  &  2.0878 & 9.1775e-6  & 2.4895  & 4.0515e-5 & 1.9232 \\
    & $1024$   & 9.8459e-6  &  2.1292 & 2.2694e-6  & 2.0158	& 1.0663e-5 & 1.9258\\
    & $2048$  &  2.3281e-6  &  2.0804  & 5.6466e-7   & 2.0067  	&2.7673e-6  &   1.9461\\		
\hline
    & $128 $     & 1.1550e-2 &   -       & 3.2017e-3  & -        & 9.8045e-4 &   -  \\
    & $256$    & 2.4445e-3 & 2.2406 & 4.8802e-4 & 2.7138 & 2.3419e-4 & 2.0658 \\
V) & $512 $    & 5.7208e-4 & 2.0950 & 8.1516e-5  & 2.5818 & 5.7877e-5 & 2.0166 \\
   & $1024$   &  1.3069e-4 & 2.1301  & 1.3263e-5  & 2.6197 & 1.4503e-5 & 1.9966\\
   & $2048$   &  3.1132e-5 & 2.0697 & 2.2433e-6 & 2.5637 & 3.6735e-6 &   1.9811\\		
\hline
     & $128$    & 1.1111e-1     & -         & 3.0202e-2 &     -      & 6.7794e-3 &   --\\
     & $256$   & 2.4720e-2  &2.1682   & 4.8231e-3  & 2.6466 & 1.5697e-3 & 2.1107 \\
VI) & $512$    & 5.8618e-3  & 2.0763 & 8.1463e-4  & 2.5657  & 3.6077e-4  & 2.1213 \\
    & $1024$   &  1.3419e-3  & 2.1271  & 1.3239e-4  &  2.6214 & 8.2855e-5 & 2.1224\\
    & $2048$  &  3.2004e-4 & 2.0680 & 2.2370e-5  & 2.5652 &1.9796e-5   &   2.0654\\		
\hline
\end{tabular}
\end{center}
\label{tabl2-e1-2}
\end{table}

\begin{table}[h]
\caption{Maximum error and convergence rates of {\em Example} 1 with viscosity $\mu^+\gg\mu^-$. }
\begin{center}
\vspace{-0.1cm}
\begin{tabular}{|c|c|c|c|c|c|c|c|}
\hline
case & N &$\|e_{\vect  u}\|_{\infty}$ 
 & order & $\|e_{\vect  u}\|_{1, \infty}$& order& $\|e_p\|_{\infty}$& order\\
 \hline
     &$128$    & 1.0493e-3  & -            & 1.4173e-4   & -          & 2.7033e-3 & - \\
     & $256$   & 2.1608e-4  &  2.2798   & 3.5083e-5  & 2.0143  & 1.3241e-3  & 1.0297 \\
IV) & $512$    & 5.0796e-5  &  2.0888  & 8.7282e-6  & 2.0070  & 6.5537e-4 & 1.0146 \\
    & $1024$   & 1.1625e-5   &  2.1275   & 2.1768e-6  & 2.0035  & 3.2614e-3  & 1.0068\\
    & $2048$  & 2.7459e-6   & 2.0819   & 5.4353e-7  & 2.0018  & 1.6267e-3  &  1.0035 \\	
\hline
    & $128$     & 1.1885e-2  &        -  & 5.8630e-4 & -        & 2.7599e-3&  - \\
    & $256$    & 2.5295e-3 & 2.2322 & 6.4863e-5 & 3.1762 & 1.3448e-3 & 1.0372 \\
V) & $512$     & 5.9335e-4 & 2.0919 & 8.7280e-6 & 2.8937 & 6.5822e-4 & 1.0307 \\
    & $1024$   & 1.3590e-4 & 2.1263  & 2.1767e-6 & 2.0035 & 3.2688e-4 & 1.0098\\
    & $2048$  & 3.2406e-5 & 2.0682 & 5.4353e-7 & 2.0017 & 1.6285e-4 &  1.0052 \\	
\hline
    & $128$     &1.1253e-1  &   -       & 4.9656e-3 &  -        & 4.4718e-3&    - \\
    & $256$    & 2.5108e-2 & 2.1641  & 5.5985e-4 & 3.1489 & 1.4363e-3 & 1.6385\\
VI) & $512$    & 5.9629e-3 & 2.0741 & 6.6932e-5 & 3.0643 & 6.5929e-4 & 1.1234 \\
    & $1024$   & 1.3668e-3 & 2.1252 & 7.7331e-6  & 3.1136  & 3.2760e-4 & 1.0090\\
    & $2048$  & 3.2599e-4 & 2.0679 & 9.2280e-7 & 3.0670 & 1.6292e-4 & 1.0078  \\	
\hline
\end{tabular}
\end{center}
\label{tabmax-e1-2}
\end{table}

\begin{figure}[h!]
\centering
\subfigure[$u^{(1)}$]{
\includegraphics[width=0.31\textwidth]{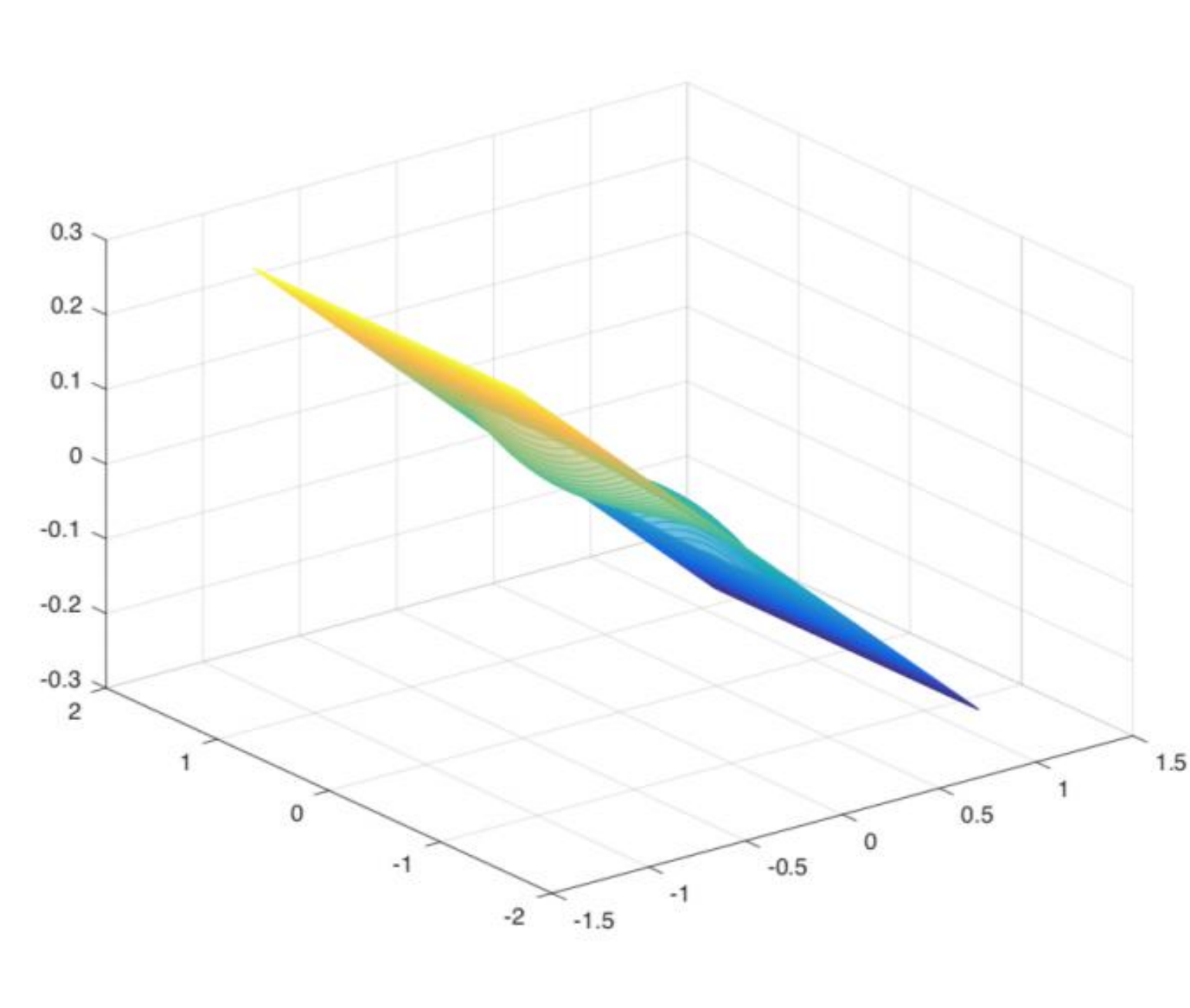}
}
\subfigure[$u^{(2)}$]{
\includegraphics[width=0.31\textwidth]{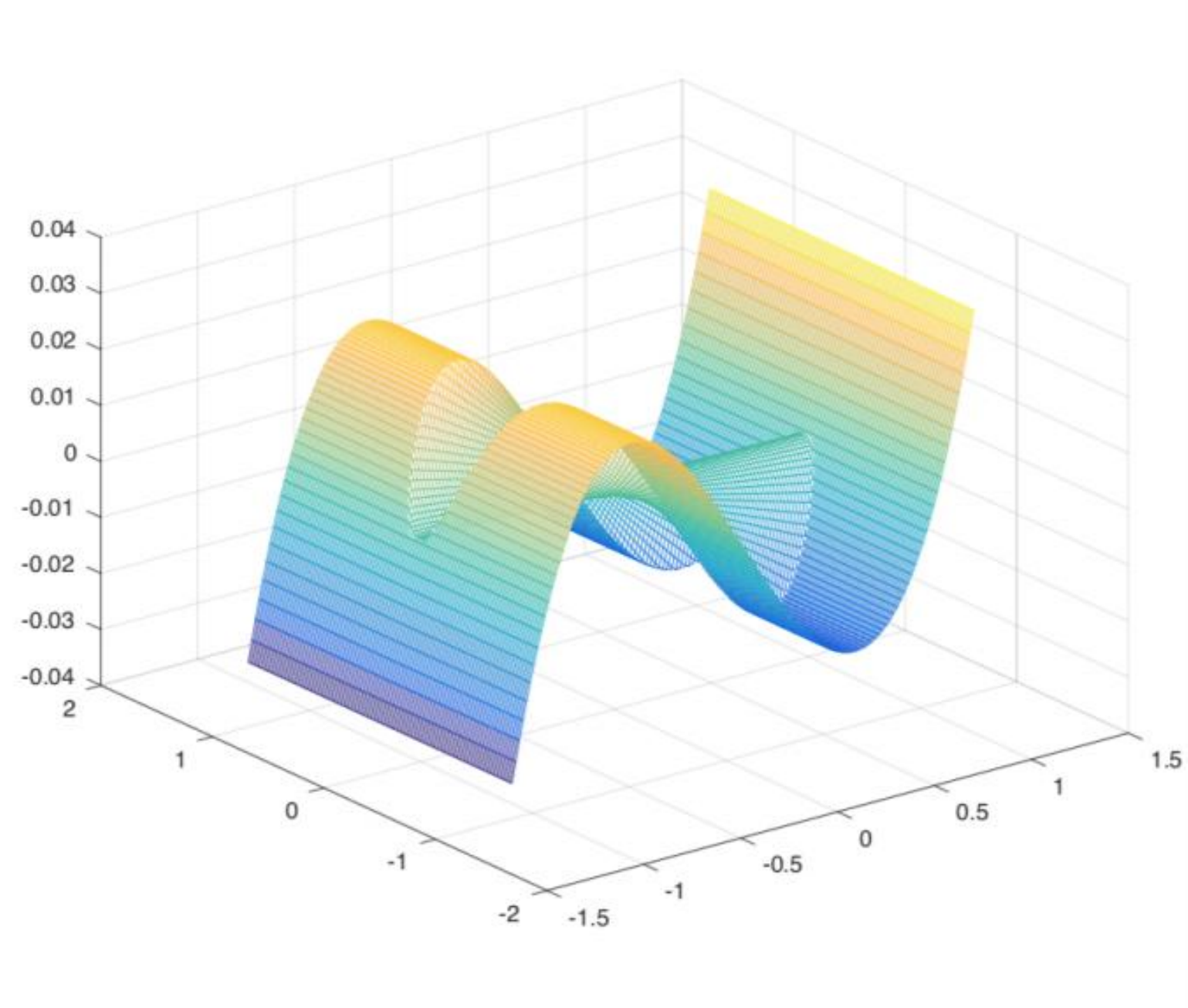}
}
\subfigure[$p$]{
\includegraphics[width=0.31\textwidth]{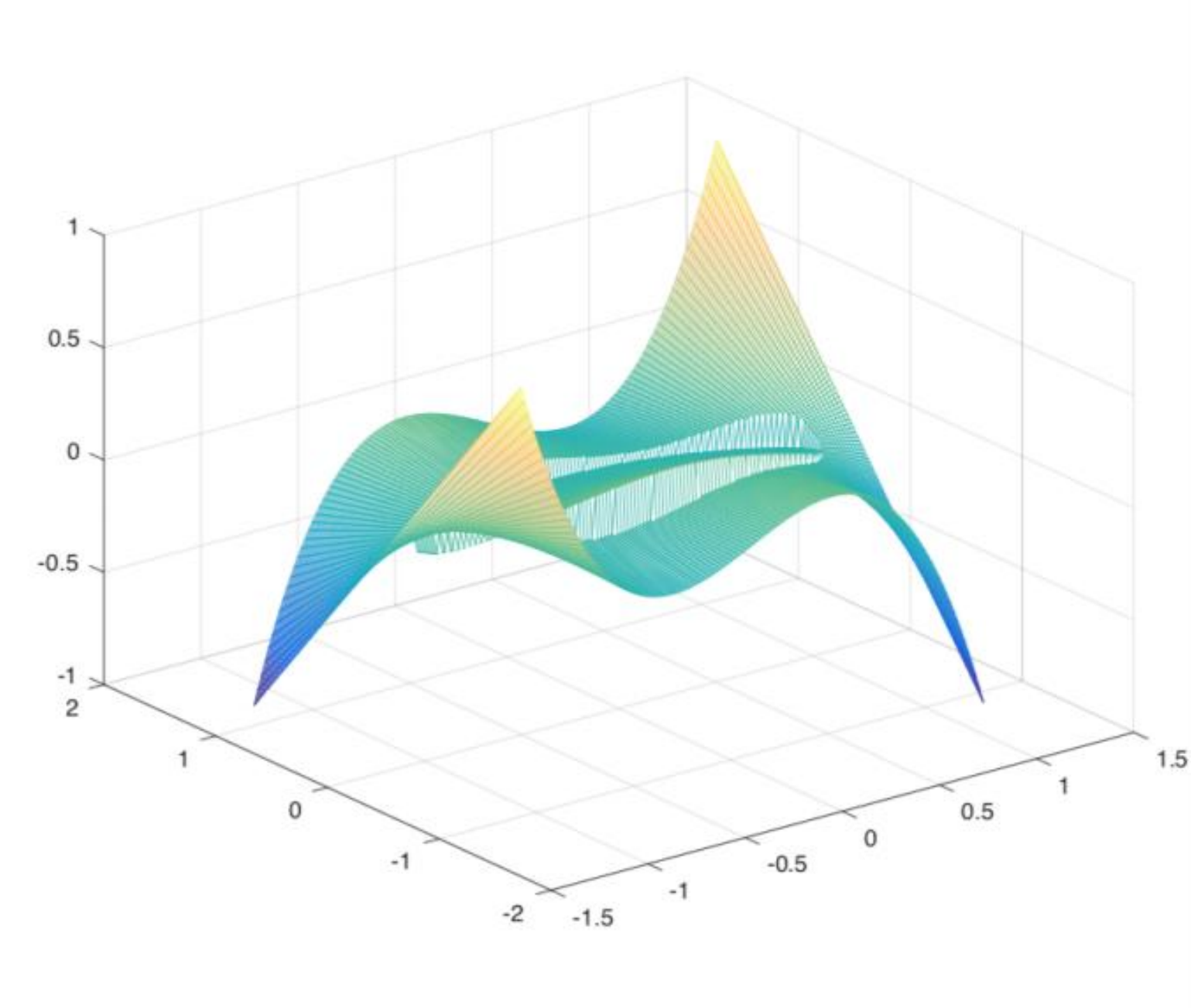}
}
\setlength{\abovecaptionskip}{-0.0cm}
\setlength{\belowcaptionskip}{-0.0cm}
 \caption{The solution plots of the $x$-component of velocity field $u^{(1)}$, the $y$-component of velocity field $u^{(2)}$ and the pressure $p$ in {\em Example} 2 with viscosity $\mu^+ = 10$ and $\mu^-=1$ on a $128\times 128$ 
 grid.}
 \label{v-p-e2}
\end{figure}

%%%%%%%%%%%%%%%%%%%%%%%%%%%%%%%%%%%%%%%%%%%%%%%%%%%%%%
{\em Example 2.} In this example,  the interface is an ellipse $\{(x, y)|x^2+4y^2=1\}$, which is located at the center of the box $\Omega=(-2.0, 2.0)^2$. The constructed exact velocity and the pressure are given by 
\begin{equation*}
\begin{split}
u^{(1)}(x, y)&=\begin{cases} 
\dfrac{y}{4}, \qquad\qquad\;\;\; x^2+4y^2>1,\\[8pt]
\dfrac{y}{4}(x^2+4y^2),\;x^2+4y^2\leq 1,
\end{cases}\\[4pt]
u^{(2)}(x, y)&=\begin{cases} 
-\dfrac{x}{16}(1-x^2), \; x^2+4y^2>1,\\[8pt]
-\dfrac{xy^2}{4},\qquad\;\;\;\;x^2+4y^2\leq 1,
\end{cases}\\[4pt]
p(x, y)&=\begin{cases} 
(-\dfrac{3}{4}x^3+\dfrac{3}{8}x)y, \quad\; x^2+4y^2>1,\\[8pt]
\exp(\sin y +\cos x),\;x^2+4y^2\leq 1.
\end{cases}
\end{split}
\end{equation*}
The external force $\vect f$ and the surface tension can be easily obtained by making the exact solution satisfy the model problem \eqref{interfaceP}. Obviously, it also has a finite jump across the interface. 
The discrete $l^2$-errors  and their corresponding convergence rates for the velocity $\vect u$, the gradient of the velocity $\nabla \vect u$ and the pressure $p$ are summarized in Tables \ref{tabl2-e2-1} and \ref{tabl2-e2-2}.  As expected, the measured rates of second order accuracy are obtained  for the velocity $\vect u$, the pressure $p$ as well as the gradient  of velocity.
The maximum norms and their corresponding convergence rates are listed in Tables \ref{tabmax-e2-1} and \ref{tabmax-e2-2}. It can be seen that  
 the velocity and its derivatives are second order accurate but the pressure is only of first order accuracy in different coefficient ratios cases. 
The numbers of GMRES iterations are also presented in Table \ref{tab-GMRES}. It shows that the GMRES iteration converges with a relatively small number of iterations, and  that the number of GMRES iterations for each case does not depend on the size of the mesh and is insensitive to the ratio of coefficients.
Finally, one can see Fig. \ref{v-p-e2} for the numerical solution $u^{(1)}, u^{(2)},p$ on a $128\times 128$ grid. 

\begin{table}[h]
\caption{$l^2$-error and convergence rates of {\em Example} 2 with viscosity $\mu^+\ll\mu^-$.}
\begin{center}
\vspace{-0.1cm}
\begin{tabular}{|c|c|c|c|c|c|c|c|}
\hline
case & N  & $\|e_{\vect  u}\|$
& order & $\|e_{\vect  u}\|_1$ & order &$\|e_p\|$ &order\\
\hline
   & $128$     & 2.8613e-4  &   -       & 2.4885e-4  &     -     & 3.4594e-4&   - \\
   & $256$    &  6.9571e-5  & 2.0401 & 6.1457e-5   & 2.0176  & 9.5385e-5& 1.8587 \\
I) & $512$     & 1.7846e-5   & 1.9629 & 1.5293e-5   & 2.0067 & 2.6346e-5& 1.8562 \\
   & $1024$   &  4.4958e-6  & 1.9890 & 3.8083e-6  & 2.0057 & 7.1185e-6 & 1.8879\\
   & $2048$  &  1.1380e-6   & 1.9821  & 9.4990e-7   & 2.0033 & 1.9090e-6 &   1.8988\\	
\hline
   & $128$     & 2.9591e-4  &    -     & 2.5003e-4 &     -    & 3.8382e-4  &   -  \\
   & $256$    & 7.1598e-5  & 2.0472 & 6.1577e-5  & 2.0216 & 1.0473e-4  & 1.8738 \\
II) & $512$    & 1.8376e-5  & 1.9621  & 1.5312e-5  & 2.0077  & 2.8918e-5 & 1.8566 \\
   & $1024$   & 4.6265e-6 & 1.9898  & 3.8108e-6 & 2.0065 & 7.7875e-6 & 1.8927\\
   & $2048$  &  1.1730e-6  & 1.9797  & 9.5024e-7 & 2.0037 & 2.0857e-6 &   1.9006\\		
\hline
    & $128$     & 2.9693e-4  &   -      & 2.5016e-4  &     -    & 3.8578e-4 &     -  \\
    & $256$    &  7.1809e-5  & 2.0479 & 6.1590e-5 & 2.0221 & 1.0515e-4  & 1.8753 \\
III) & $512$    & 1.8429e-5   & 1.9622 & 1.5314e-5   & 2.0078 & 2.9032e-5& 1.8567 \\
    & $1024$   &  4.6396e-6 & 1.9899 & 3.8111e-6   & 2.0066 & 7.8153e-6 & 1.8933\\
    & $2048$  &  1.1765e-6   & 1.9795 & 9.5028e-7 & 2.0038 & 2.0929e-6 &   1.9008\\	
\hline
\end{tabular}
\end{center}
\label{tabl2-e2-1}
\end{table}

\begin{table}[h]
\caption{Maximum error and convergence rates of {\em Example} 2 with viscosity $\mu^+\ll\mu^-$. }
\begin{center}
\vspace{-0.1cm}
\begin{tabular}{|c|c|c|c|c|c|c|c|}
\hline
case & N &$\|e_{\vect  u}\|_{\infty}$ 
 & order & $\|e_{\vect  u}\|_{1, \infty}$& order& $\|e_p\|_{\infty}$& order\\
 \hline
   & $128$     & 2.5062e-4 &      -   & 2.4883e-4 &     -   & 1.8413e-3 &  - \\
   & $256$    & 6.1799e-5 &2.0198  & 6.1627e-5 & 2.0135 & 9.3308e-4 & 0.9807 \\
I) & $512$     & 1.5358e-5 & 2.0086 & 1.5329e-5 & 2.0073 & 4.6955e-4 & 0.9907 \\
   & $1024$   & 3.8269e-6 & 2.0047 & 3.8234e-6 & 2.0033 & 2.3556e-4 & 0.9952\\
   & $2048$  & 9.5524e-7 &2.0022 & 9.5483e-7 & 2.0015    & 1.1797e-4 &  0.9977 \\	
\hline
   & $128$    & 2.5117e-4 &  	-    & 2.4915e-4 &	-        & 1.8404e-3 &    -\\
   & $256$    & 6.1875e-5 &2.0212 & 6.1677e-5 & 2.0142  & 9.3290e-4 & 0.9802 \\
II) & $512$    & 1.5368e-5 & 2.0094& 1.5335e-5 & 2.0079 & 4.6950e-4& 0.9906\\
   & $1024$   & 3.8281e-6 & 2.0052& 3.8241e-6 & 2.0036 & 2.3555e-4 & 1.0714\\
   & $2048$  & 9.5542e-7 &2.0024 & 9.5493e-7 & 2.0017 & 1.1797e-4 &  0.9976 \\	
\hline
    &$128$     &2.5122e-4 &	-     & 2.4918e-4 &	-        & 1.8403e-3 &  - \\
    & $256$    & 6.1883e-5 &2.0213 & 6.1682e-5 & 2.0143 & 9.3288e-4 & 0.9802 \\
III) & $512$    & 1.5369e-5 &2.0095 & 1.5336e-5 &2.0079 & 4.6950e-4& 0.9906\\
    & $1024$   & 3.8282e-6 &2.0053 &3.8242e-6 &2.0037 & 2.3555e-4 & 0.9951\\
    & $2048$  & 9.5552e-7 &2.0023 &9.5494e-7 & 2.0017& 1.1797e-4 &  0.9976 \\	
\hline
\end{tabular}
\end{center}
\label{tabmax-e2-1}
\end{table}

\begin{table}[h]
\caption{$l^2$-error and convergence rates of {\em Example} 2 with viscosity $\mu^+\gg\mu^-$.}
\begin{center}
\vspace{-0.1cm}
\begin{tabular}{|c|c|c|c|c|c|c|c|}
\hline
case & N  & $\|e_{\vect  u}\|$
& order & $\|e_{\vect  u}\|_1$ & order &$\|e_p\|$ &order\\
\hline
    & $128$     & 1.1353e-3  &	-	& 5.1395e-4 &		-  & 4.9017e-4& -  \\
    & $256$    &  3.3715e-4 &1.7516 & 1.0404e-4 & 2.3045 & 1.2394e-4& 1.9836 \\
IV) & $512$    &  6.7553e-5 &2.3193 & 1.7069e-5 & 2.6077& 3.2020e-5& 1.9526 \\
    & $1024$   &   1.6800e-5 &2.0076 &4.0379e-6 &2.0797 & 8.1789e-6 & 1.9690\\
    & $2048$  &  4.0694e-6 &2.0456 & 9.7789e-7 &2.0459 & 2.0752e-6 &   1.9787\\		
\hline
    & $128$     & 1.3347e-2 &	-      &4.7427e-3 & 		-  & 8.1285e-4&  -  \\
    & $256$    &  3.9362e-3 &1.7616 & 9.8891e-4 & 2.7008 & 2.1297e-4& 1.9323\\
V) & $512$     &  7.9674e-4 &2.3046 & 1.4270e-4 & 2.2618 & 5.2058e-5& 2.0325 \\
    & $1024$   &   1.9820e-4&2.0072 & 2.5145e-5 & 2.7929 & 1.3219e-5& 1.9775\\
    & $2048$  &  4.8025e-5 & 2.0451&4.3175e-6  & 2.5420  &3.3244e-6&   1.9914\\		
\hline
    & $128$     & 1.2924e-1 & -       &  4.4696e-2 &-      & 3.7146e-3&  -  \\
    & $256$    &  3.9748e-2&1.7011 & 9.7587e-3 & 2.1954& 1.1523e-3& 1.6887\\
VI) & $512$    &  8.1888e-3 &2.2792 & 1.4241e-3 & 2.7766& 2.3640e-4& 2.2852 \\
    & $1024$   &   2.0413e-3 &2.0042 & 2.5131e-4 &2.5025 & 5.9278e-5& 1.9957\\
    & $2048$  &  4.9512e-4 &2.0436 & 4.3128e-5 &2.5428 & 1.4348e-5&   2.0466\\		
\hline
\end{tabular}
\end{center}
\label{tabl2-e2-2}
\end{table}

\begin{table}[h]
\caption{Maximum error and  convergence rates of {\em Example} 2 with viscosity $\mu^+\gg\mu^-$. }
\begin{center}
\vspace{-0.1cm}
\begin{tabular}{|c|c|c|c|c|c|c|c|}
\hline
case & N &$\|e_{\vect  u}\|_{\infty}$ 
 & order & $\|e_{\vect  u}\|_{1,\infty}$& order& $\|e_p\|_{\infty}$& order\\
 \hline
    & $128$     & 1.1260e-3 & 	-     & 2.3883e-4&		-  &7.6651e-4& -\\
    & $256$    & 3.2109e-4 & 1.8102& 6.0281e-5 &1.9557 & 3.5304e-4 & 1.1185 \\
IV) & $512$    & 6.6345e-5& 2.2749&1.5171e-5 &1.9904 & 1.6965e-4 & 1.0573 \\
    & $1024$   & 1.6480e-5 &2.0093 & 3.8037e-6 &1.9958 & 8.3009e-5 & 1.0298\\
    & $2048$  & 4.0260e-6 &2.0333 & 9.5234e-7 &1.9979 & 4.1040e-5 & 1.0176 \\	
\hline
   & $128$     & 1.3136e-2 & 	-    & 5.8057e-4 &	-       & 8.3545e-4& -\\
   & $256$    & 3.8672e-3 &1.7642 & 8.6564e-5 & 2.5152& 3.7252e-4 & 1.1652 \\
V) & $512$    & 7.7743e-4 & 2.3145& 1.5142e-5 & 1.9944&1.7407e-4& 1.0977 \\
   & $1024$   & 1.9283e-4 &2.0114 & 3.8003e-6 &1.9088 & 8.4192e-5 & 1.0479\\
   & $2048$  & 4.6631e-5&2.0480 & 9.5184e-7&1.9973 & 4.1314e-5 & 1.0271 \\	
\hline
     & $128$     & 1.2851e-1 &	-     & 5.4493e-3 &-	    & 2.5221e-3&  -\\
     & $256$    & 3.9440e-2 &1.7041 & 8.6052e-4 &2.6628 & 7.9502e-4 & 1.6656 \\
VI) & $512$    & 8.0886e-3 & 2.2857& 8.6927e-5 &3.3073 & 1.7583e-4& 2.1768 \\
    & $1024$   & 2.0139e-3 & 2.0059& 1.0857e-5 & 3.0012& 8.5604e-5& 1.0384\\
    & $2048$  & 4.8816e-4&2.0446& 1.3170e-6&3.0433& 4.1406e-5& 1.0478 \\	
\hline
\end{tabular}
\end{center}
\label{tabmax-e2-2}
\end{table}

\subsection{Test examples of moving interface}
Four moving interface problems are considered in this subsection. The interfaces are explicitly represented with a spline curve in the first two examples and implicitly represented by some control points in the last two examples.  Since the emphasis of this paper is on the new method for stationary Stokes equations, the details about the algorithm for the moving interface problems is omitted, but the simulation results are shown below.  In all the following examples, the homogeneous Dirichlet boundary condition is applied, i.e. $\vect u|_{\partial \Omega} = 0$ unless it is stated otherwise. The computational domain is  set to be$(-1.2, 1.2)\times (-1.2,1.2)$. Moreover, a $128\times 128$ grid is employed in the computations, and $100$ control points are used to present the interface.

\begin{figure}[h!]
\centering
%\subfigure[$u^{(1)}$]{
\subfigure{
\includegraphics[width=0.31\textwidth]{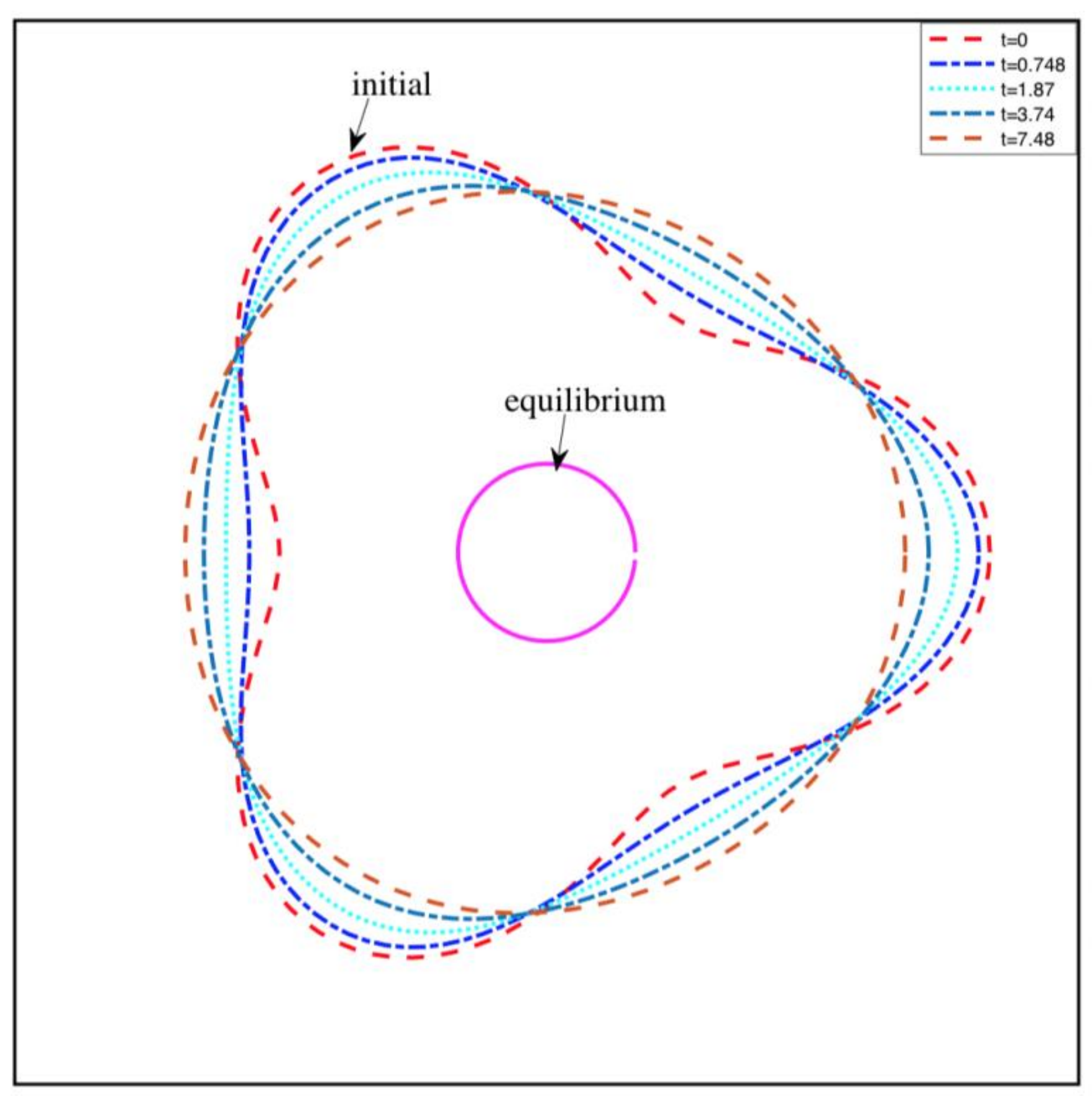}
}
%\subfigure[$u^{(2)}$]{
\subfigure{
\includegraphics[width=0.31\textwidth]{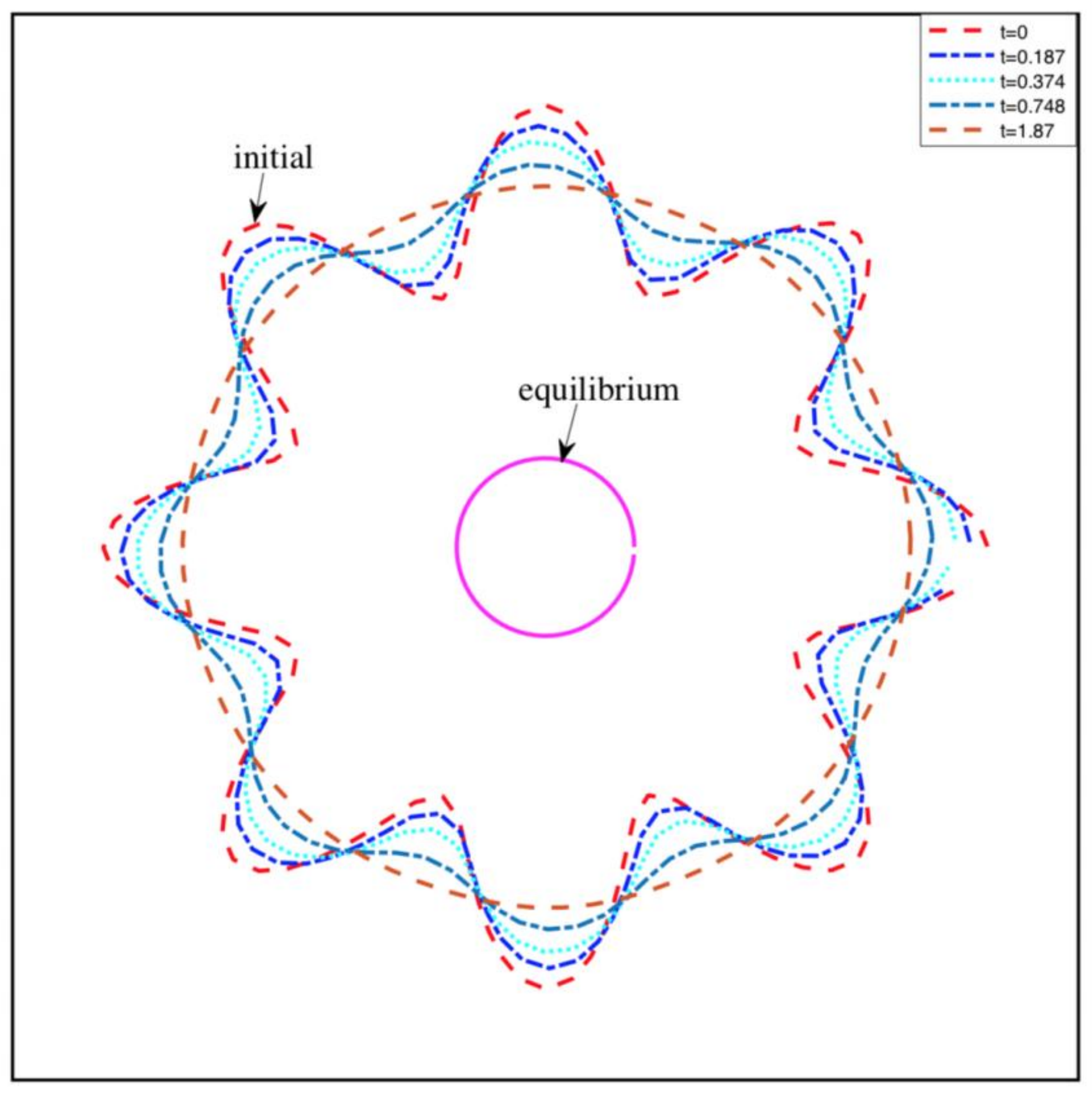}
}
\setlength{\abovecaptionskip}{-0.0cm}
\setlength{\belowcaptionskip}{-0.0cm}
 \caption{The interface configurations at different times in a square domain. (left: Three-petaled flower initial  interface; right: Eight-petaled flower initial  interface) }
\label{interface-1}
\end{figure}

\begin{figure}[h!]
\centering
\subfigure[$$]{
\includegraphics[width=0.31\textwidth]{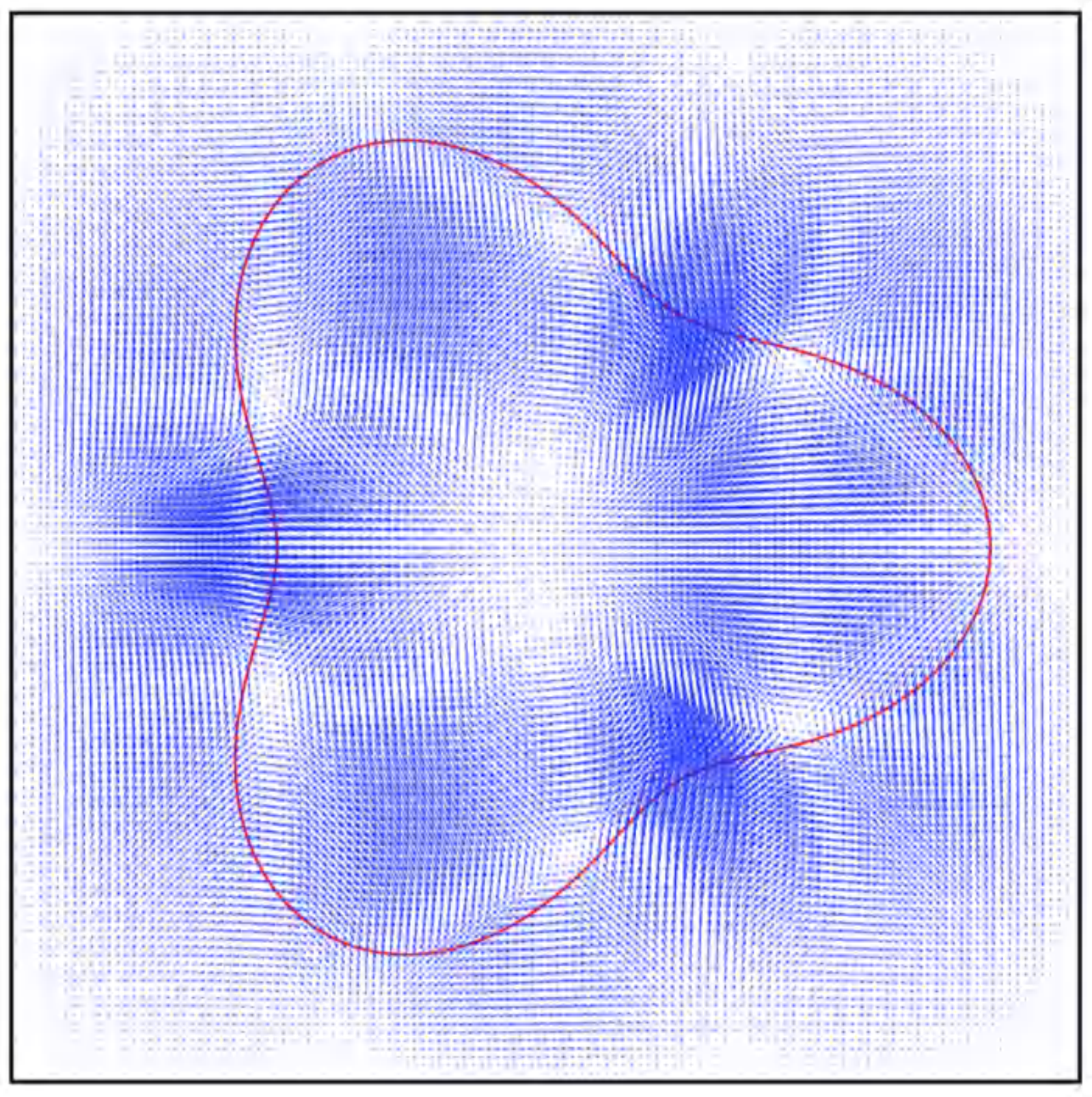}
}
\subfigure[$$]{
\includegraphics[width=0.31\textwidth]{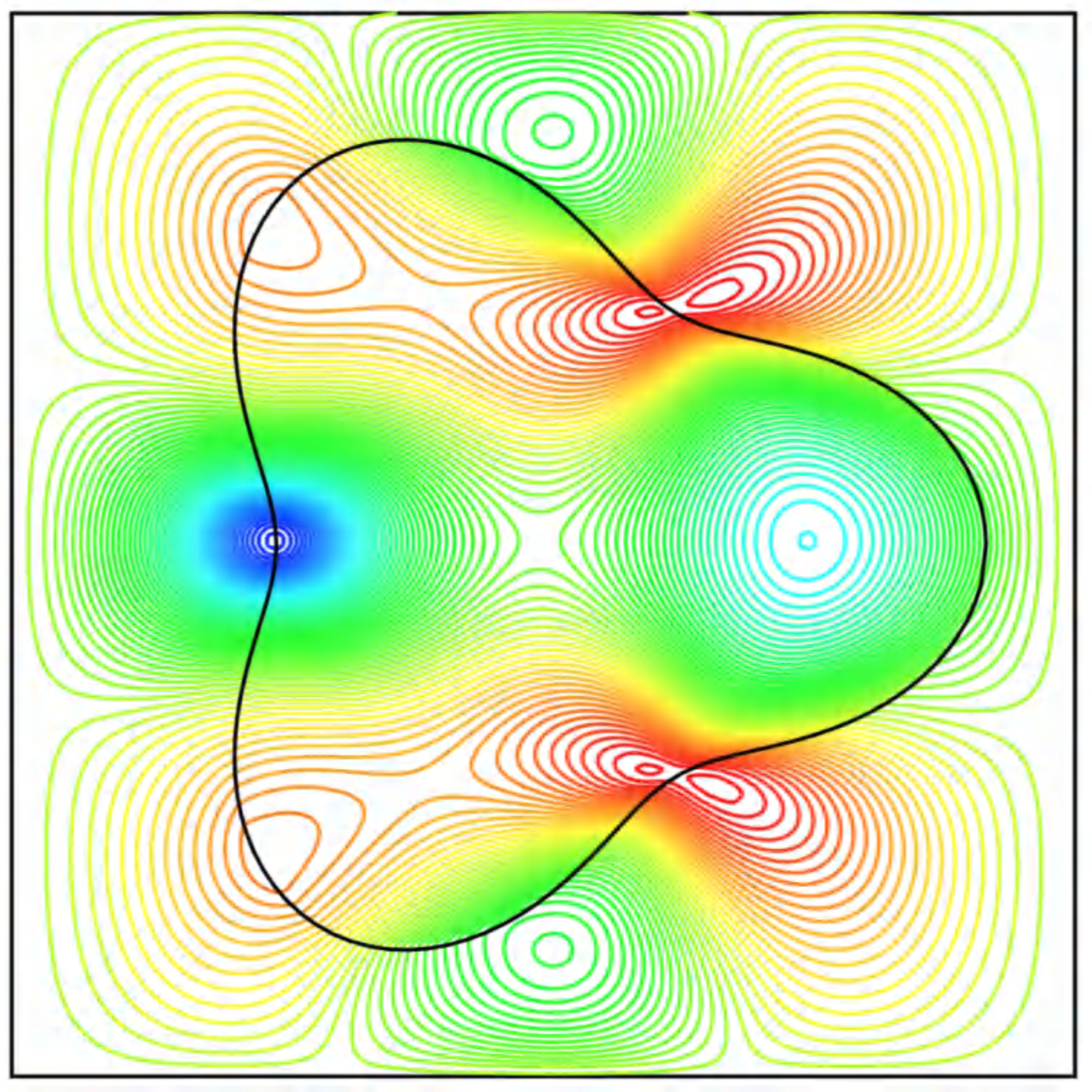}
}
\subfigure[$$]{
\includegraphics[width=0.31\textwidth]{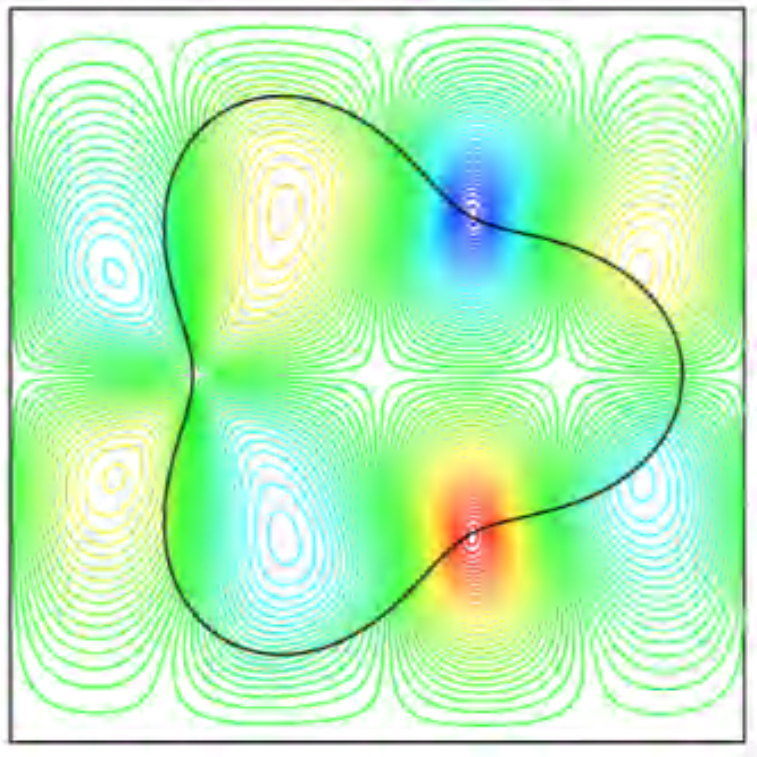}
}
\subfigure[$$]{
\includegraphics[width=0.31\textwidth]{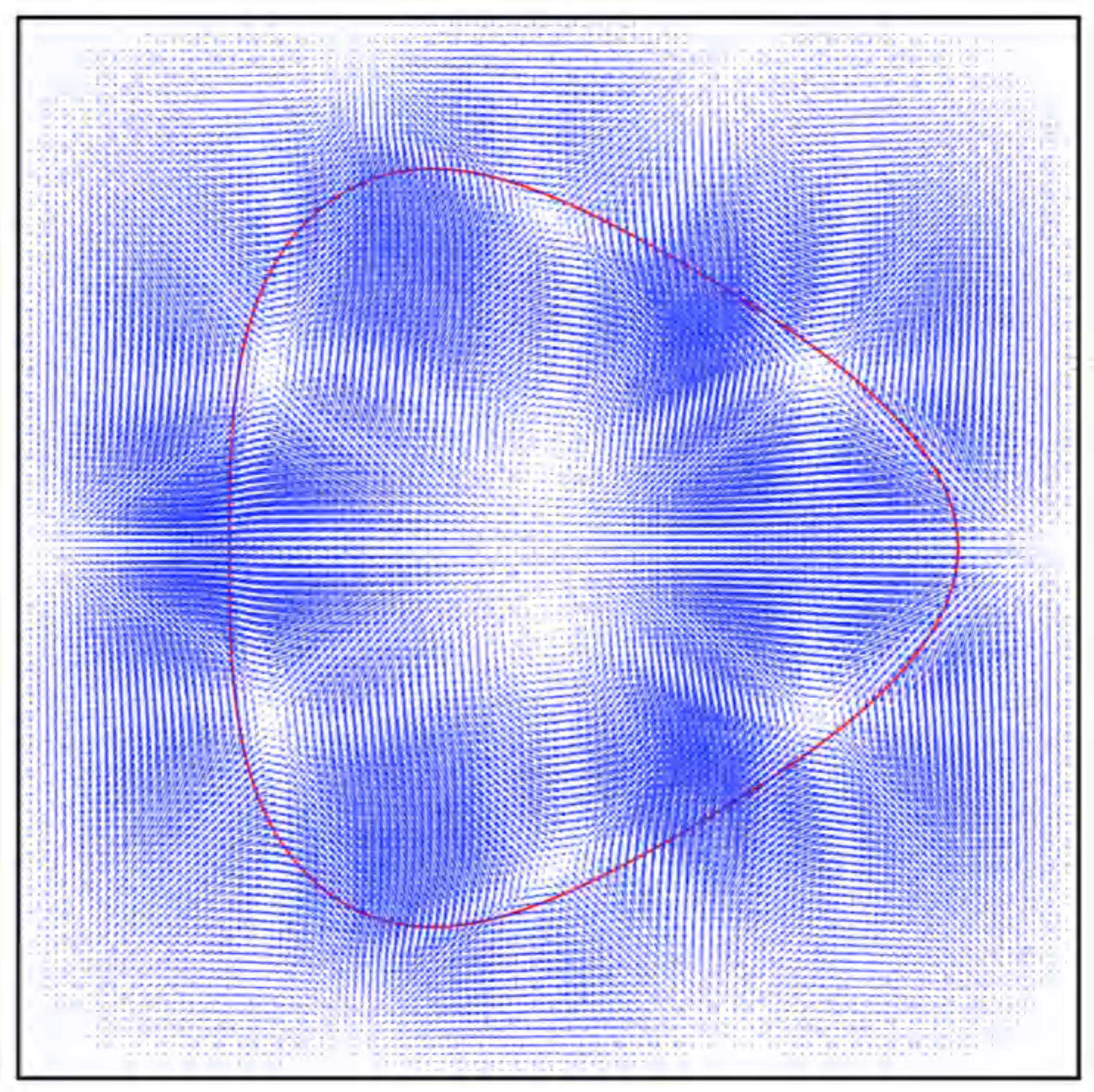}
}
\subfigure[$$]{
\includegraphics[width=0.31\textwidth]{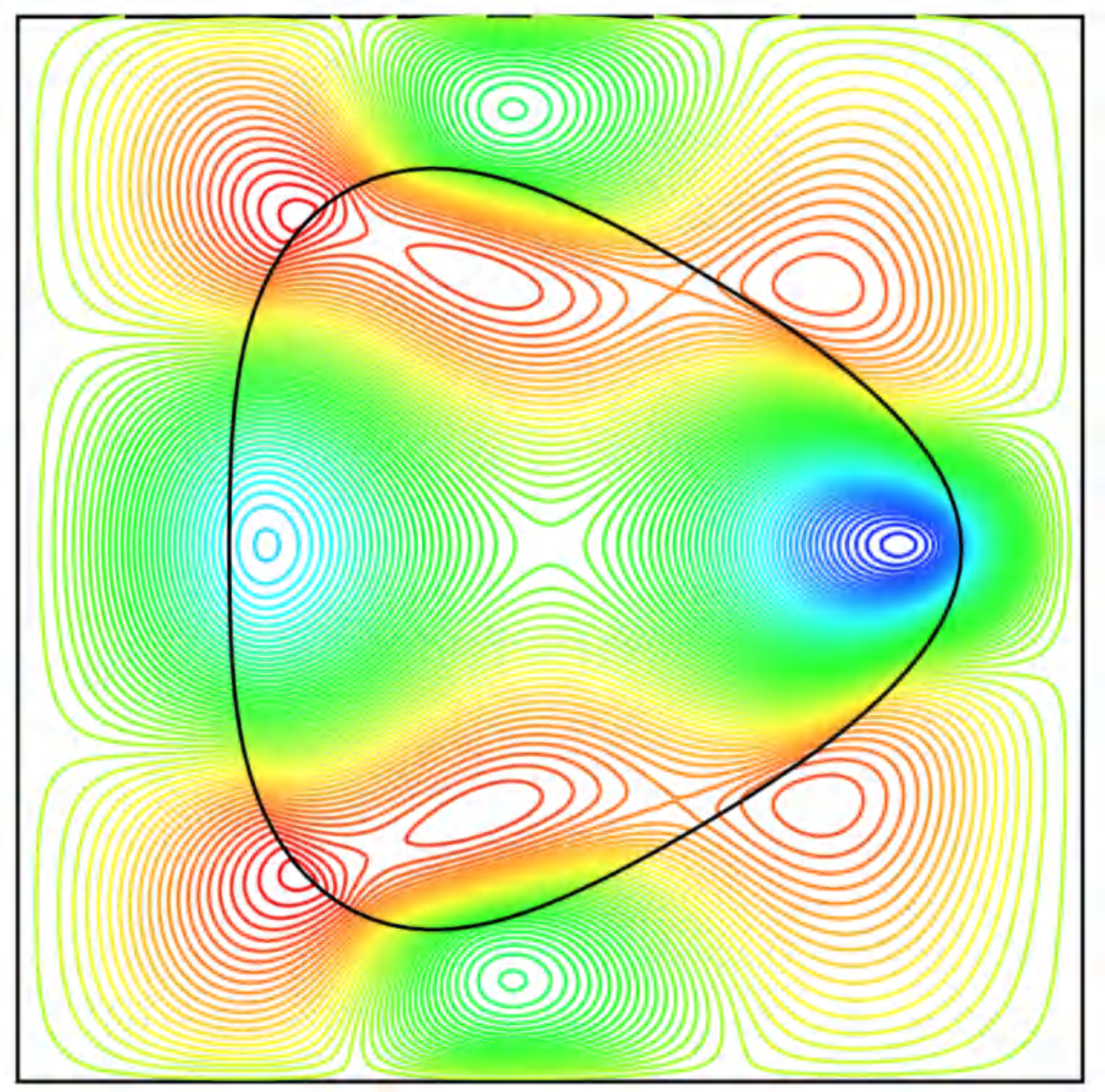}
}
\subfigure[$$]{
\includegraphics[width=0.31\textwidth]{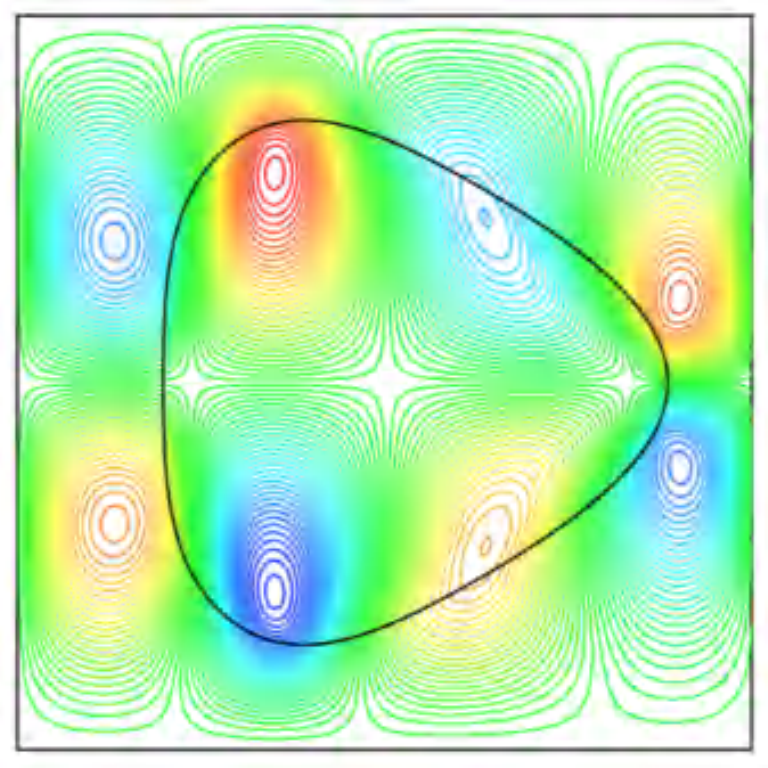}
}
\subfigure[$$]{
\includegraphics[width=0.31\textwidth]{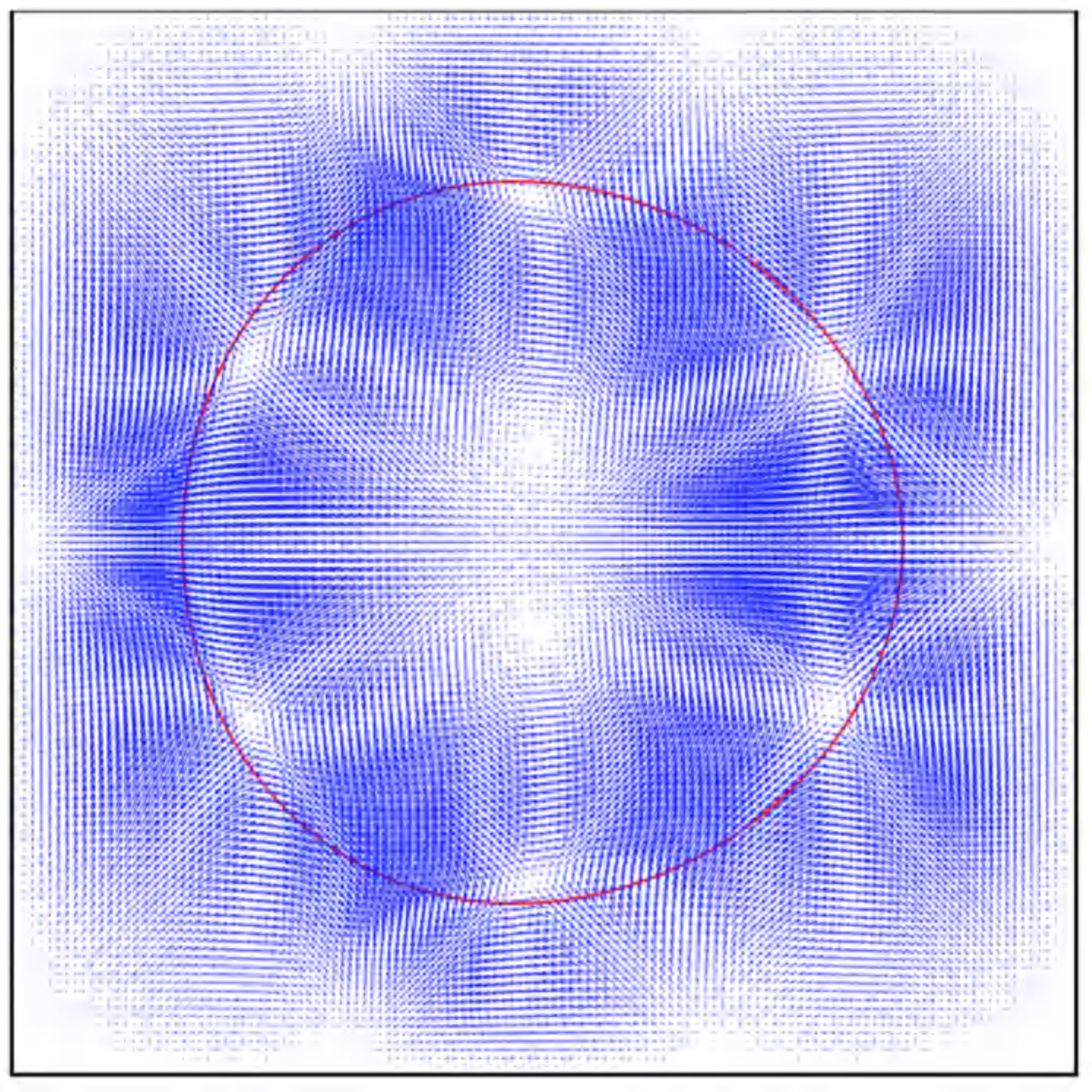}
}
\subfigure[$$]{
\includegraphics[width=0.31\textwidth]{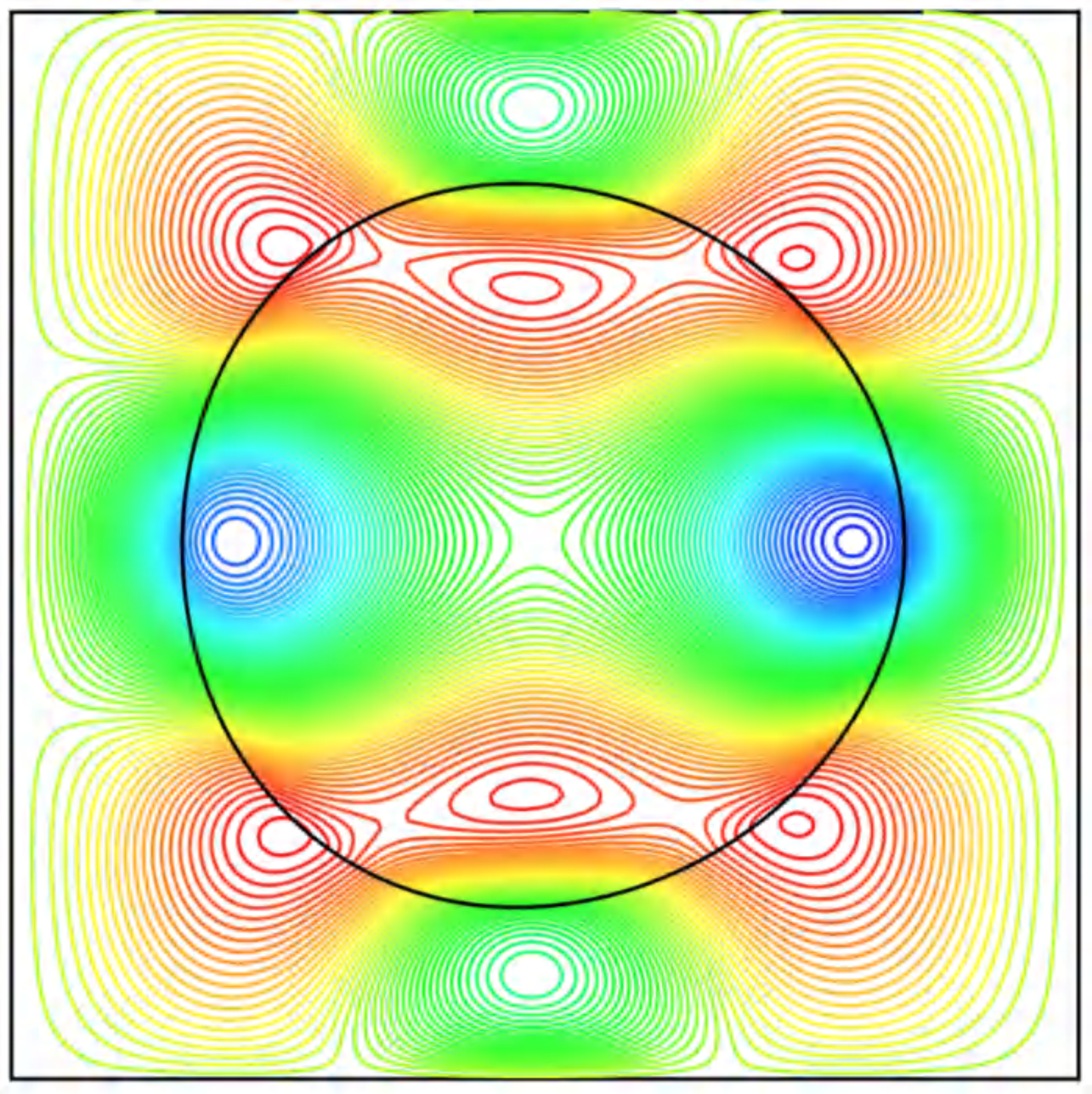}
}
\subfigure[$$]{
\includegraphics[width=0.31\textwidth]{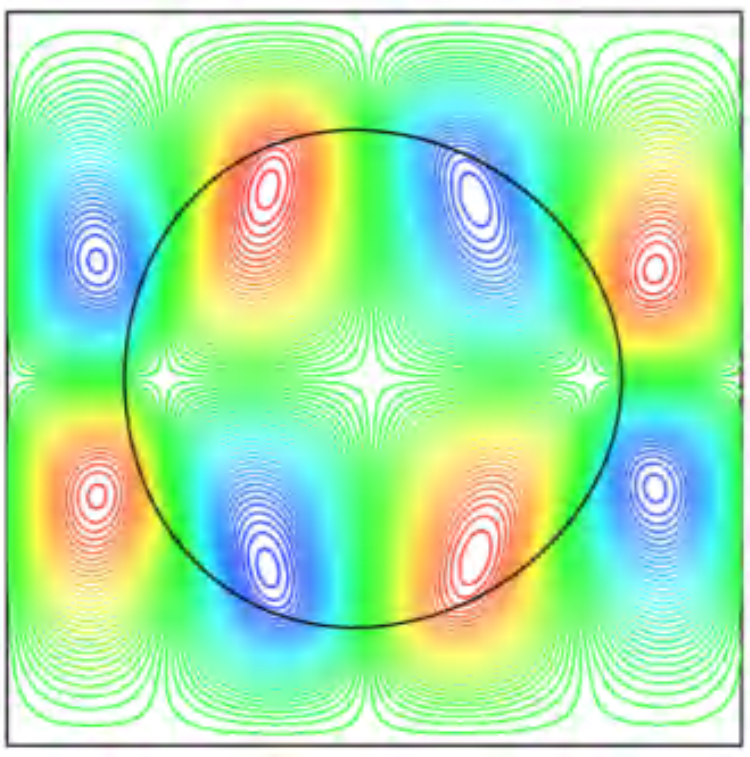}
}
\setlength{\abovecaptionskip}{-0.0cm}
\setlength{\belowcaptionskip}{-0.0cm}
 \caption{Example 3. Evolution of the velocity field $\vect u$ and interface position for three-petaled flower initial  interface problem, computed with a $128\times 128$ grid and $\Delta t = h = 0.0187$. (left: Velocity field $\vect u$; middle: Isolines of the $x$-component $u^{(1)}$; right: Isolines of the $y$-component  $u^{(2)}$.)}
 \label{velocity-e1}
\end{figure}

\begin{figure}[h!]
\centering
\subfigure[$t=0$]{
\includegraphics[width=0.31\textwidth]{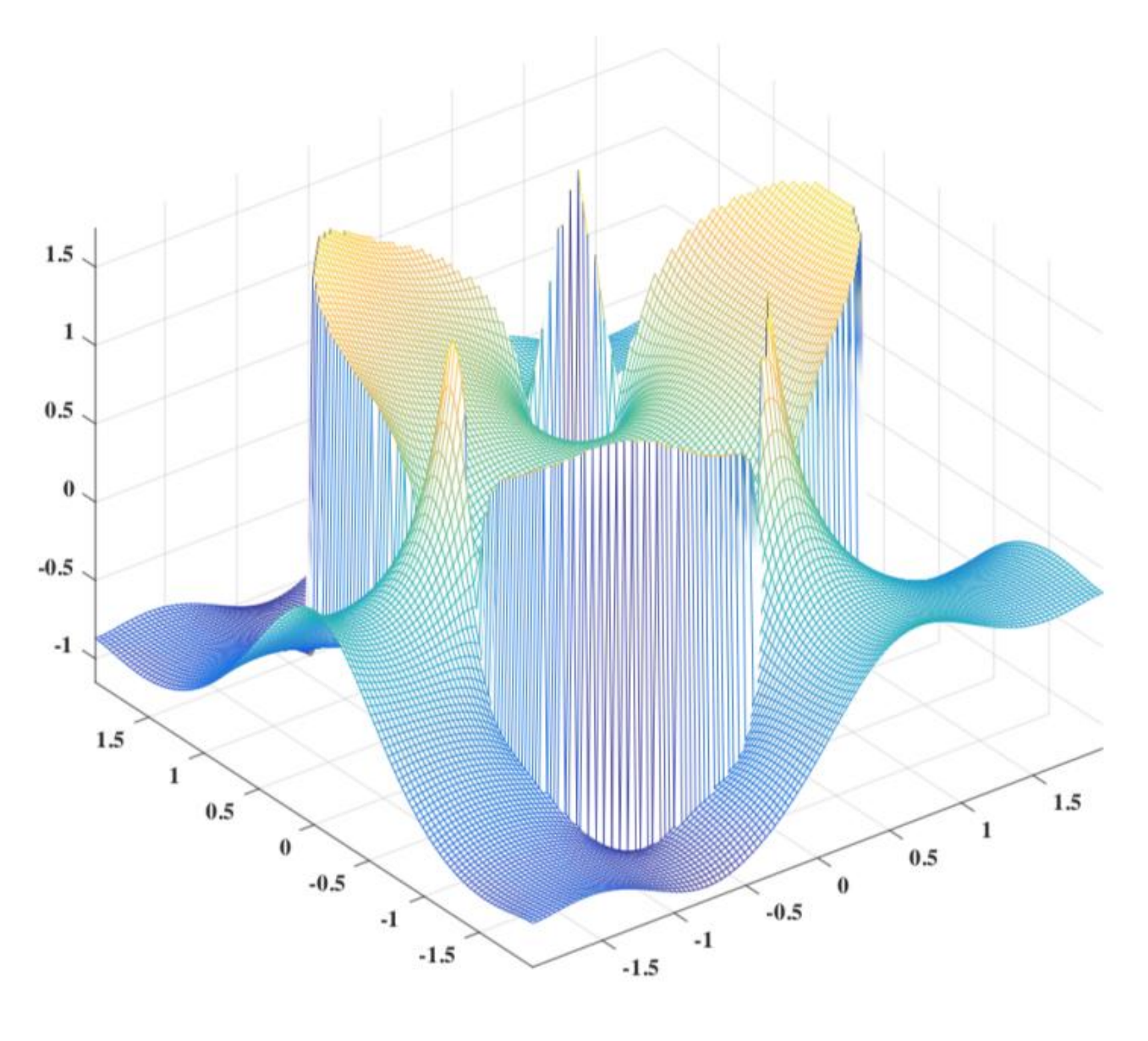}
}
%\subfigure[$t=0.374$]{
%\includegraphics[width=0.23\textwidth]{m_q_t2_e1.pdf}
%}
\subfigure[$t=1.87$]{
\includegraphics[width=0.31\textwidth]{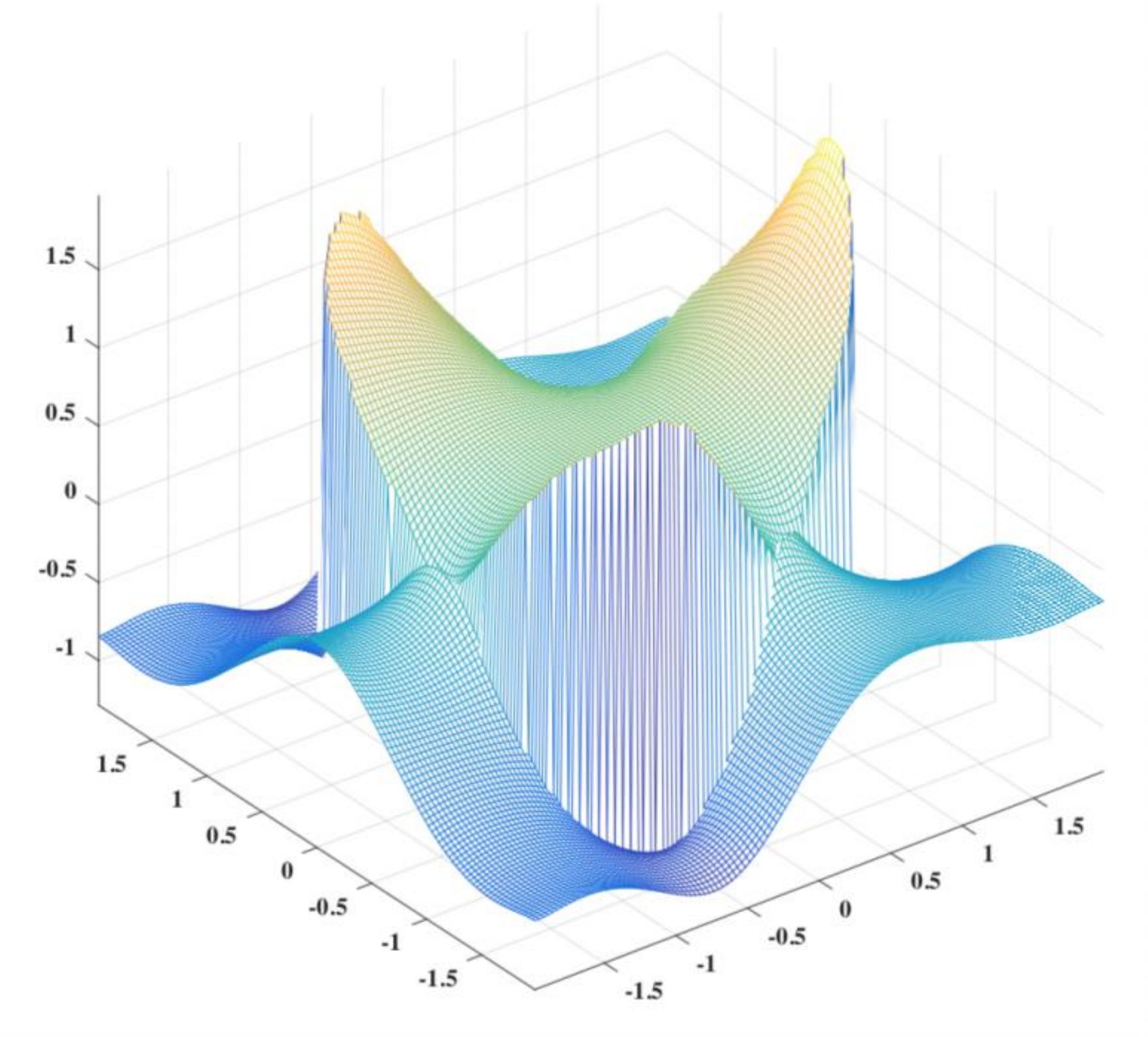}
}
\subfigure[$t=7.48$]{
\includegraphics[width=0.31\textwidth]{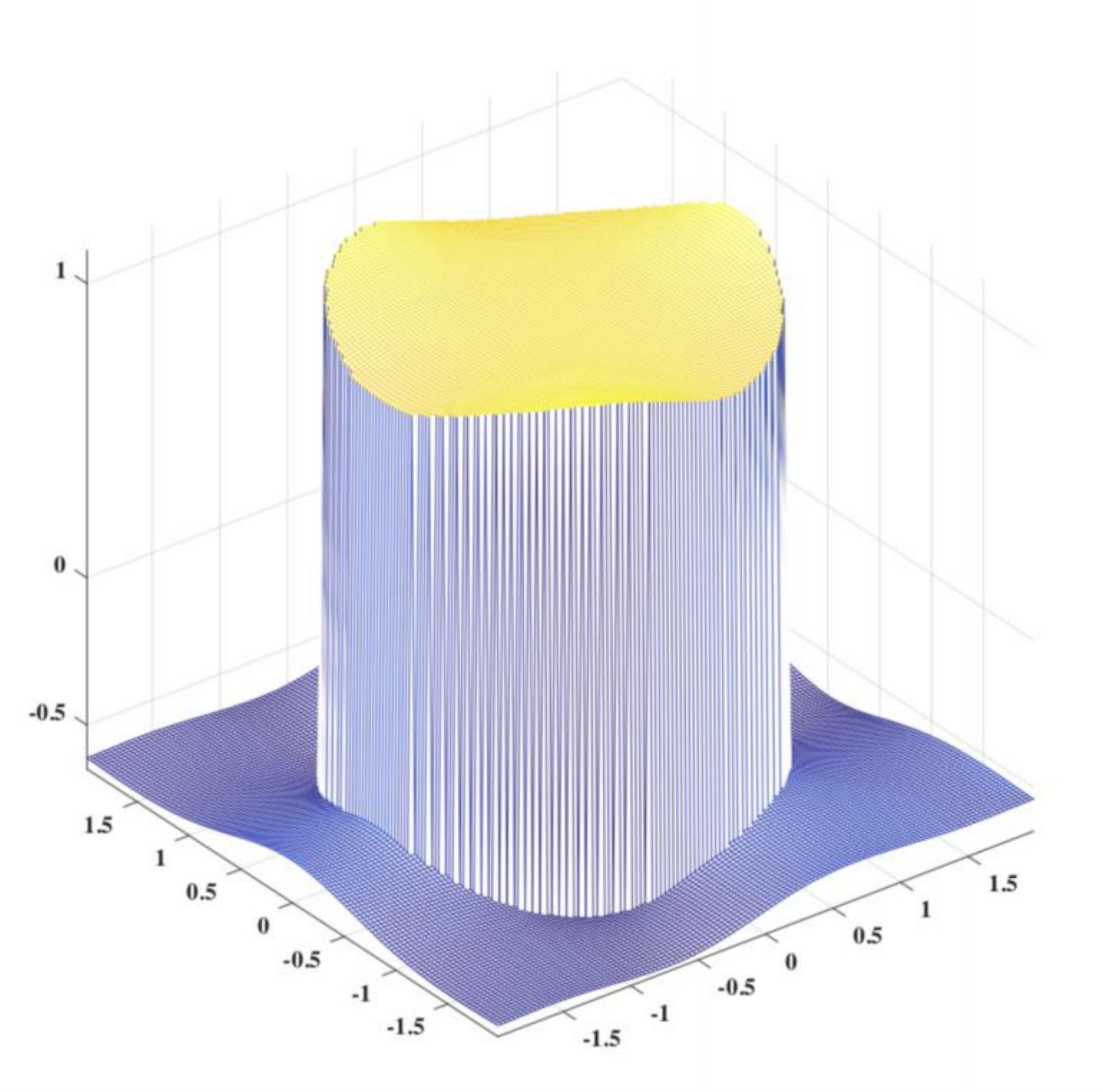}
}
\setlength{\abovecaptionskip}{-0.0cm}
\setlength{\belowcaptionskip}{-0.0cm}
 \caption{Example 3. The pressure distribution at different times.}
 \label{pressure-e1}
\end{figure}

%========================================================================
{\em Example } 3.  In this example, the initial interface is given in polar coordinates by $r=0.8+0.2\sin(3\theta), 0\leq \theta\leq2\pi$.  The interface will relax to its equilibrium,  a circle with radius $r_0 = 0.2$. 
The configurations of the interface at $t=0, t=0.748, t=1.87,t=3.74, t=7.48 $ are shown in Fig. \ref{interface-1} (left). The tension coefficient $T_0$ is set to be 0.5.  The viscosity coefficients outside the interface is taken to be $\mu^- = 1$  and inside the interface is taken to be $\mu^+ = 10$.  The approximation is computed up to a final time $T=8$ with $\Delta t = h = 0.0187$.  A time evolution of the velocity and the interface position are plotted in Fig.\ref{velocity-e1} (left),  isolines of the $x$-component $u^{(1)}$ and isolines of the $y$-component $u^{(2)}$ at different times are presented in Fig. \ref{velocity-e1} (middle) and Fig.\ref{velocity-e1} (right), respectively.  A time evolution of the pressure profile  is shown in Fig. \ref{pressure-e1}. As expected, it can be observed from these figures that  the velocity is continuous but not smooth, while the pressure is discontinuous across the interface.  The sharp jumps in the derivative of the velocity and the pressure are well captured, demonstrating the ability of the scheme to compute the velocity, pressure and interface position in each time step.

\begin{figure}[h!]
\centering
\subfigure[$$]{
\includegraphics[width=0.31\textwidth]{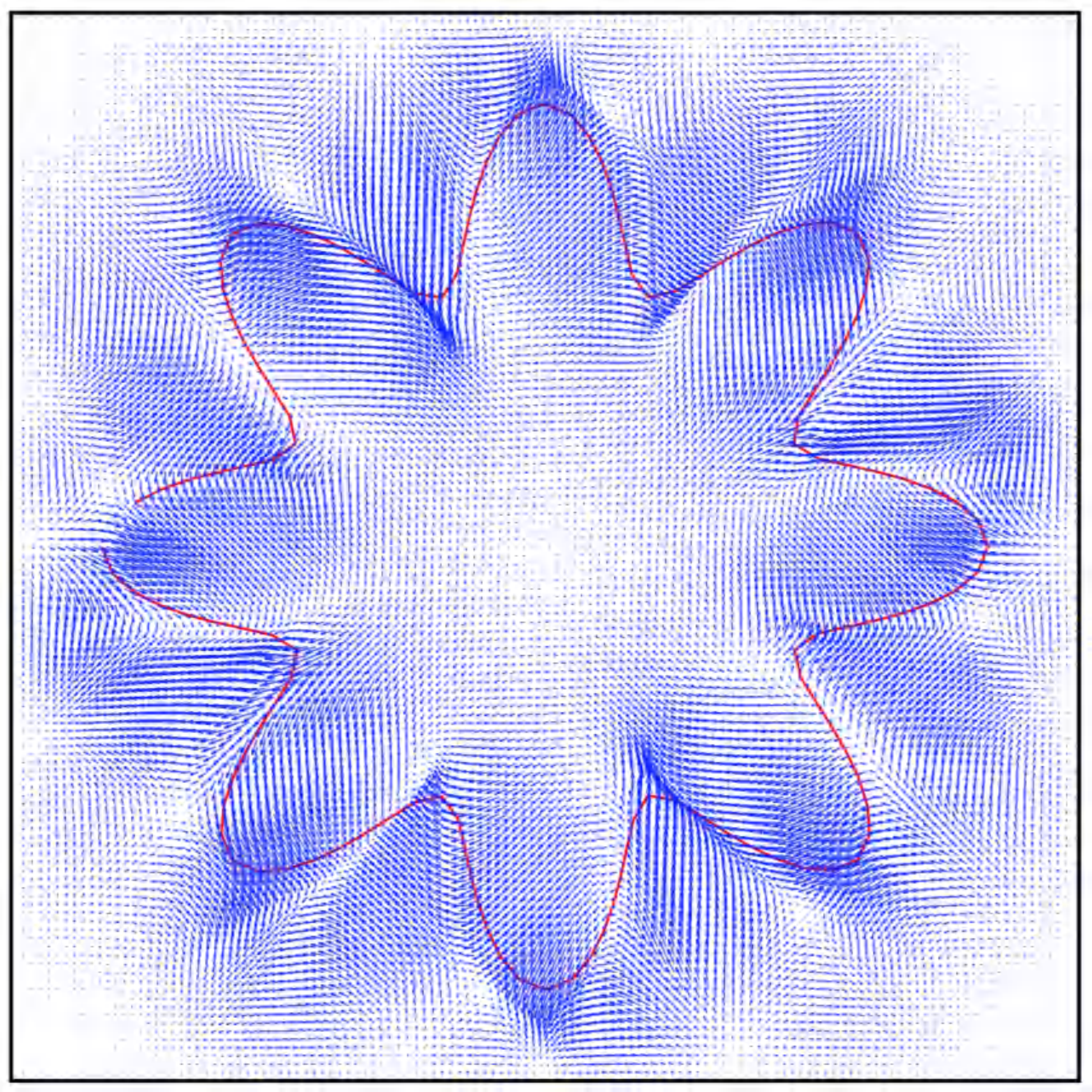}
}
\subfigure[$$]{
\includegraphics[width=0.31\textwidth]{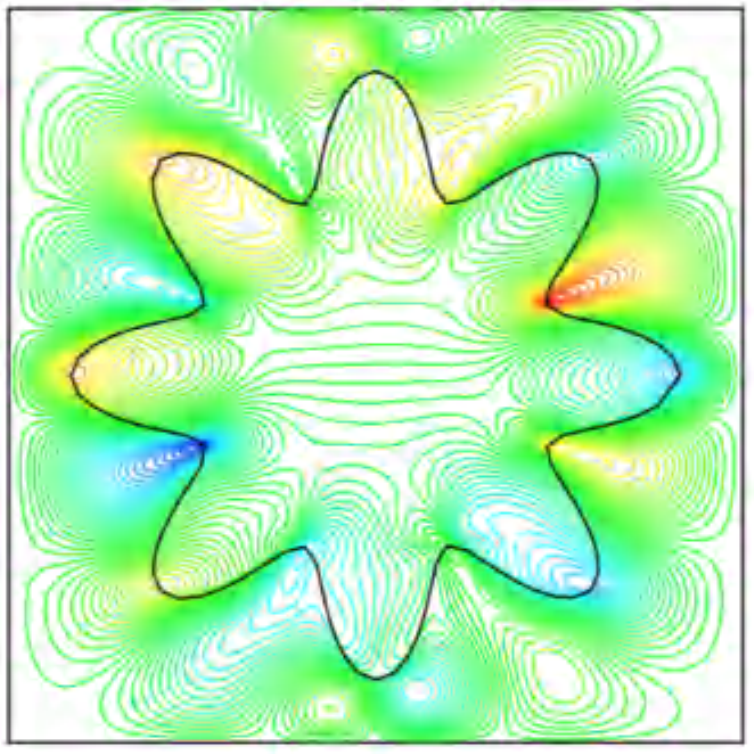}
}
\subfigure[$$]{
\includegraphics[width=0.31\textwidth]{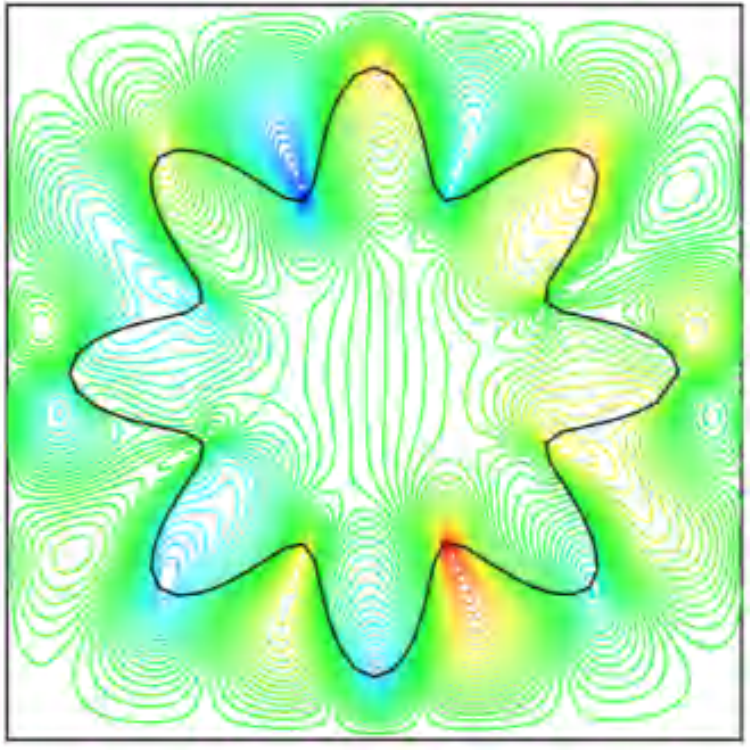}
}
\subfigure[$$]{
\includegraphics[width=0.31\textwidth]{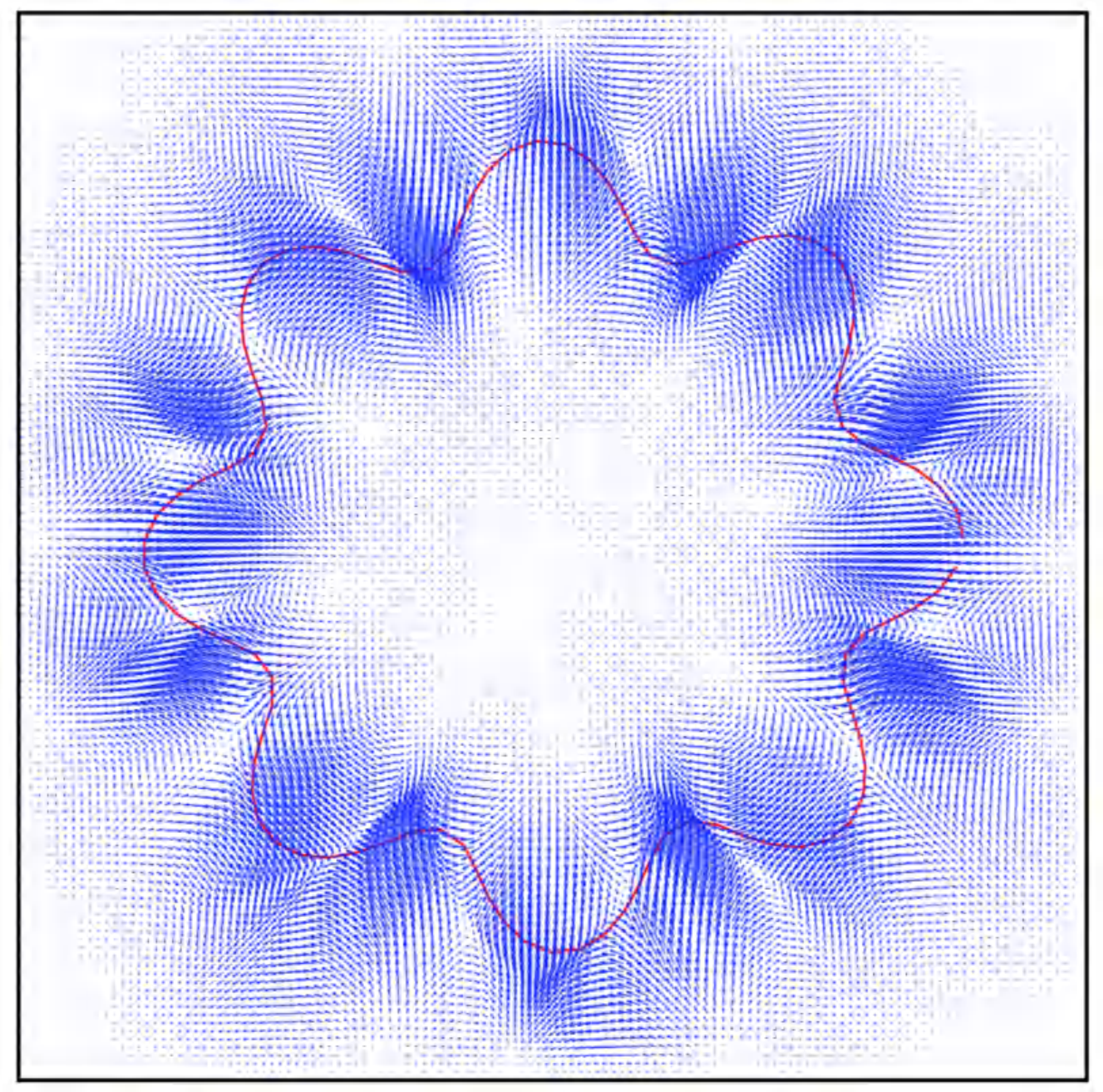}
}
\subfigure[$$]{
\includegraphics[width=0.31\textwidth]{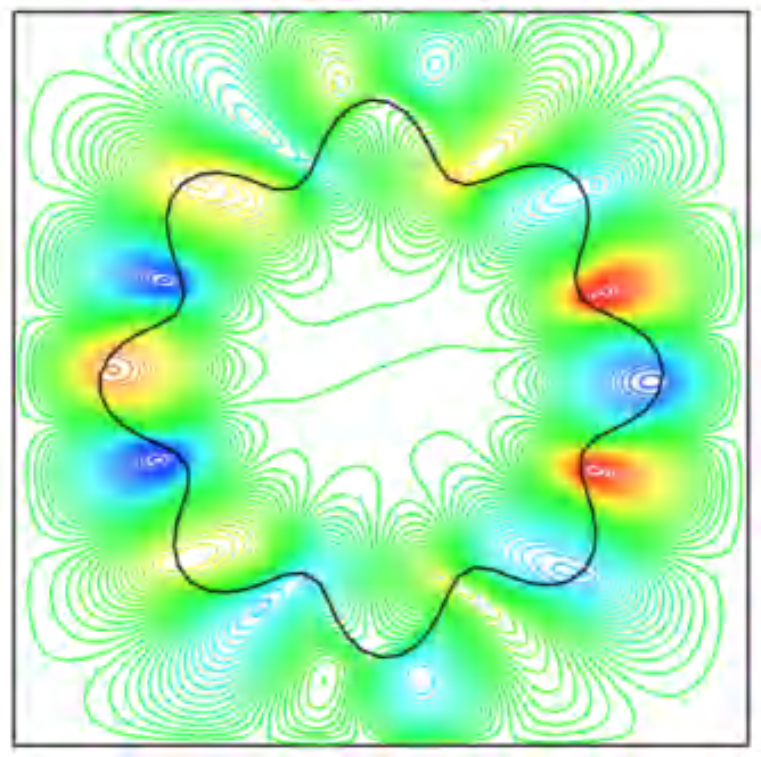}
}
\subfigure[$$]{
\includegraphics[width=0.31\textwidth]{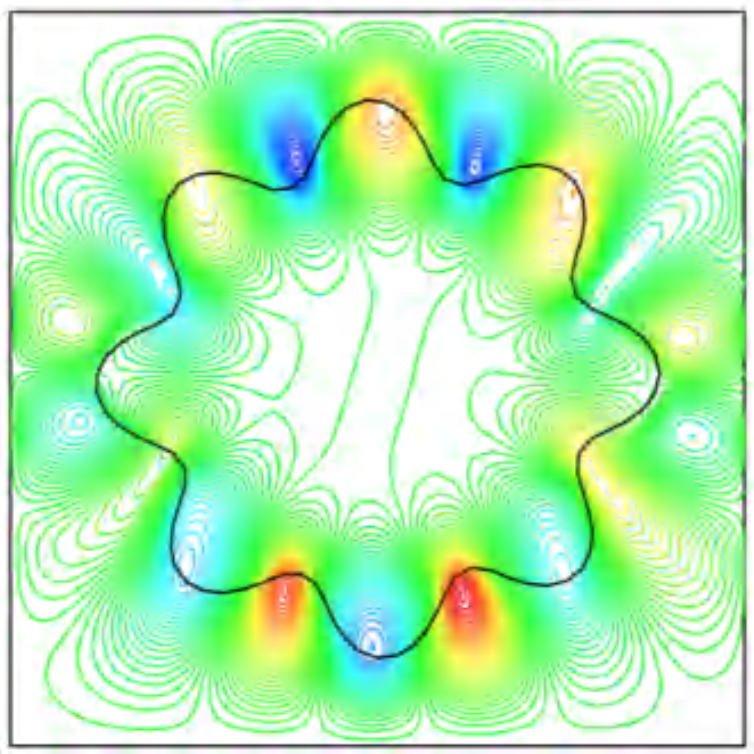}
}
\subfigure[$$]{
\includegraphics[width=0.31\textwidth]{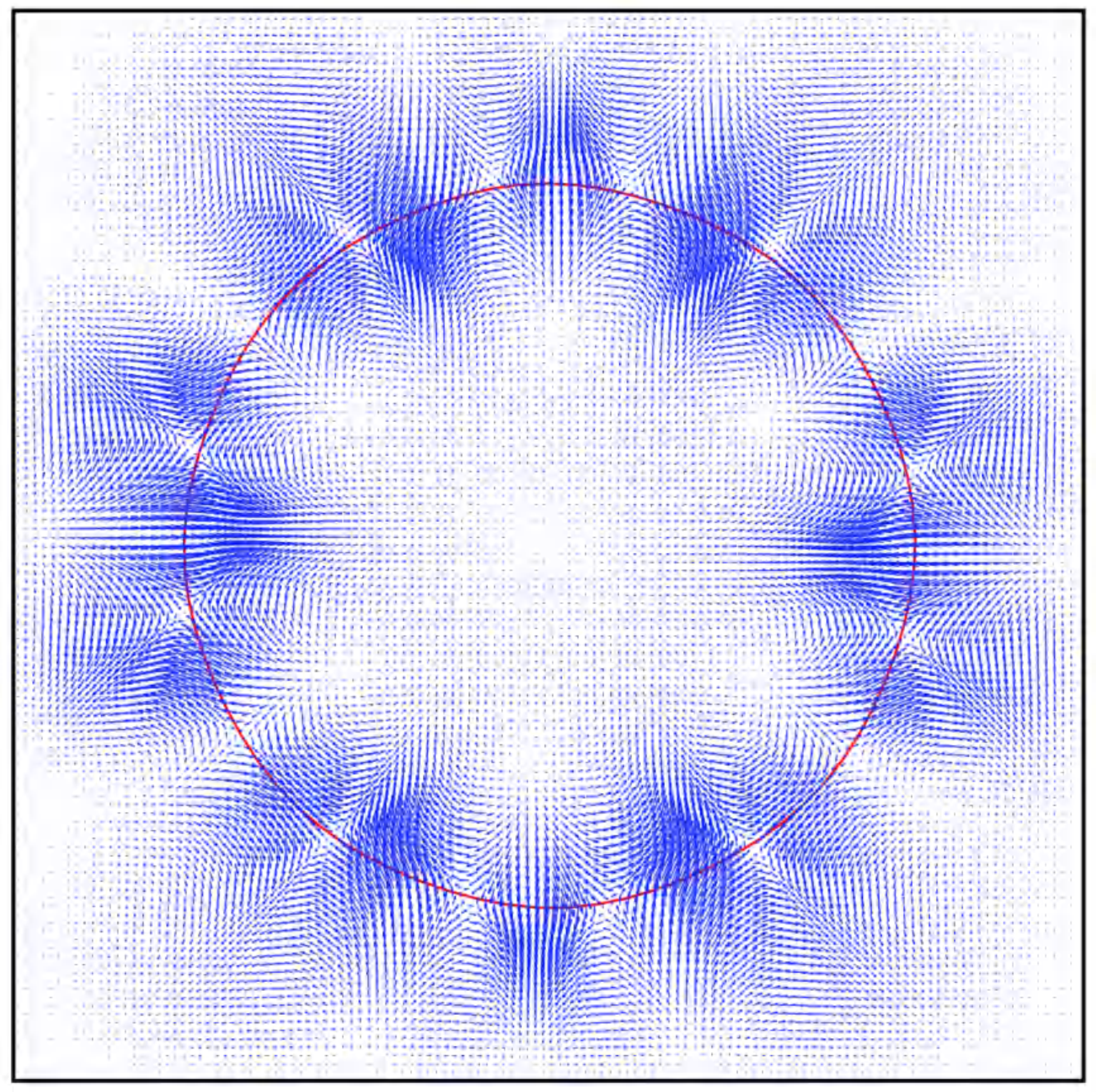}
}
\subfigure[$$]{
\includegraphics[width=0.31\textwidth]{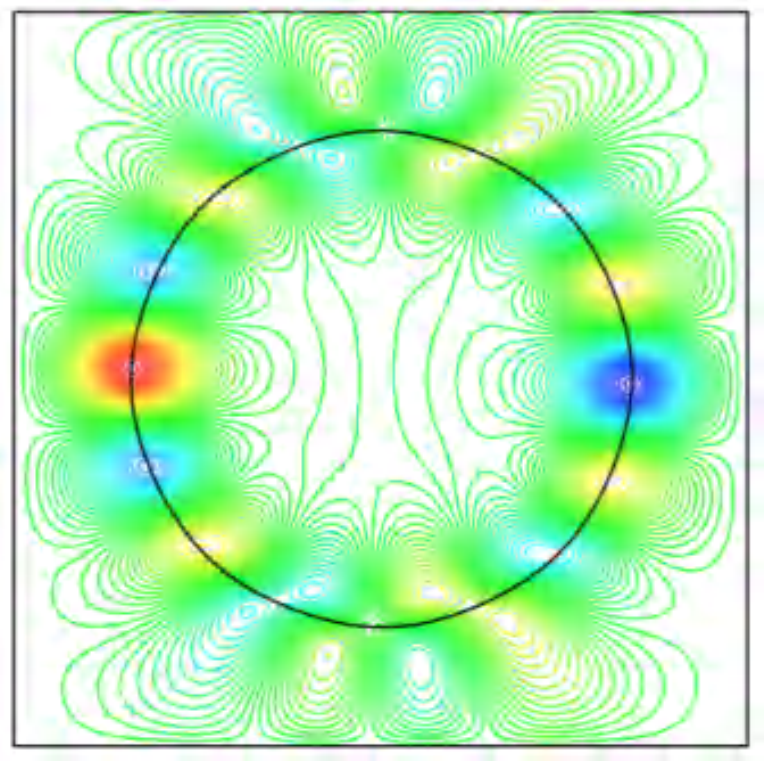}
}
\subfigure[$$]{
\includegraphics[width=0.31\textwidth]{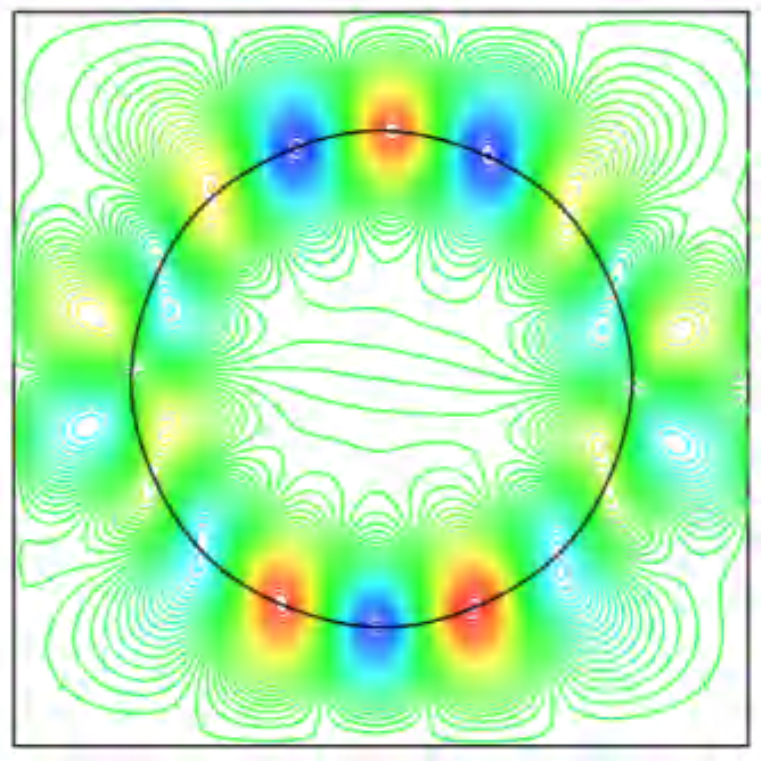}
}
\setlength{\abovecaptionskip}{-0.0cm}
\setlength{\belowcaptionskip}{-0.0cm}
 \caption{Example 4. Evolution of the velocity field $\vect u$ and interface position for five-petaled flower initial  interface problem, computed with a $128\times 128$ grid and $\Delta t = h = 0.0187$. (left: Velocity field $\vect u$; middle: Isolines of the $x$-component $u^{(1)}$; right: Isolines of the $y$-component  $u^{(2)}$.)}
 \label{velocity-e2}
\end{figure}

\begin{figure}[h!]
\centering
\subfigure[$t=0$]{
\includegraphics[width=0.31\textwidth]{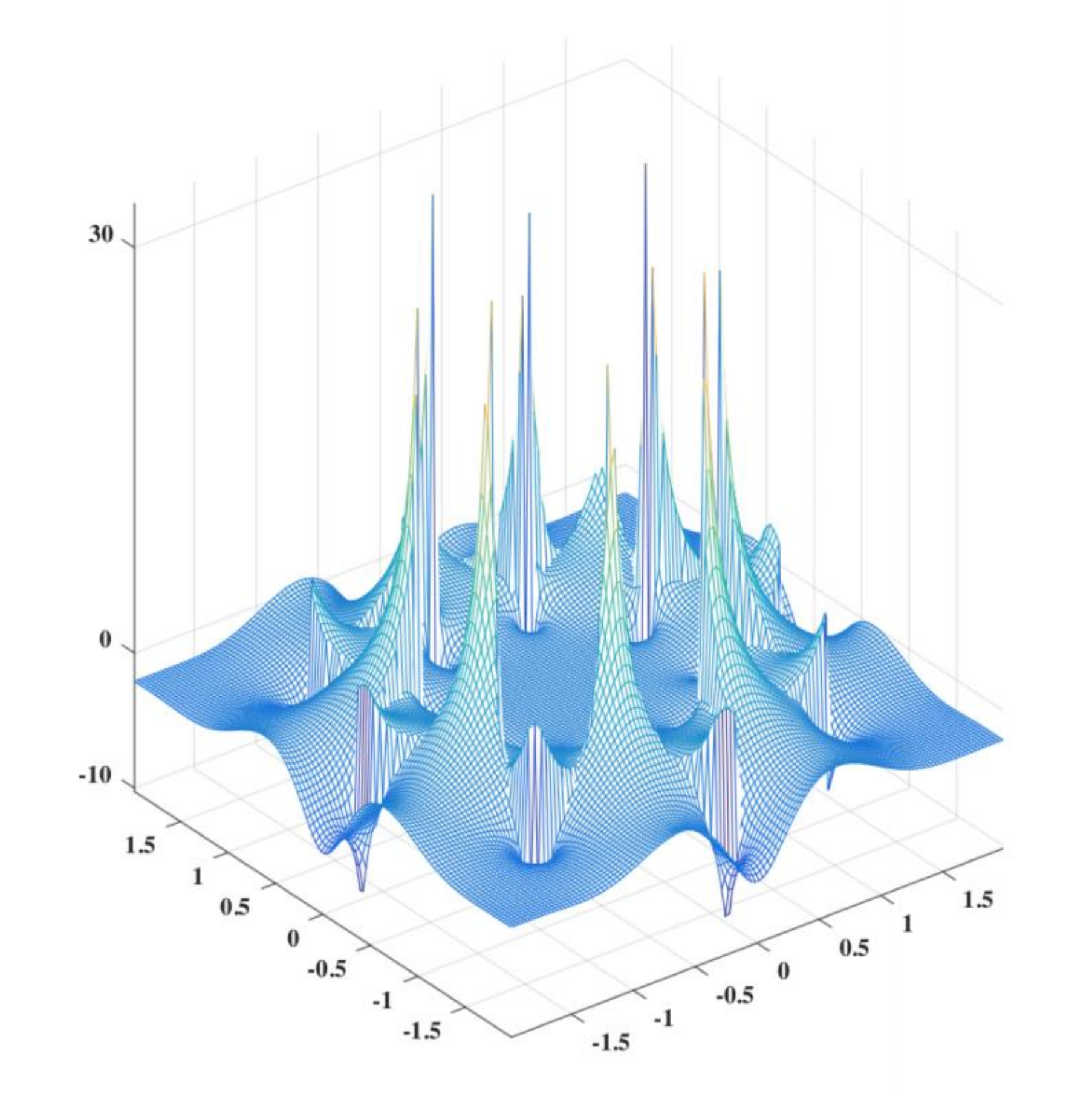}
}
%\subfigure[$t=0.374$]{
%\includegraphics[width=0.23\textwidth]{m_q_t2_e2.pdf}
%}
\subfigure[$t=0.935$]{
\includegraphics[width=0.31\textwidth]{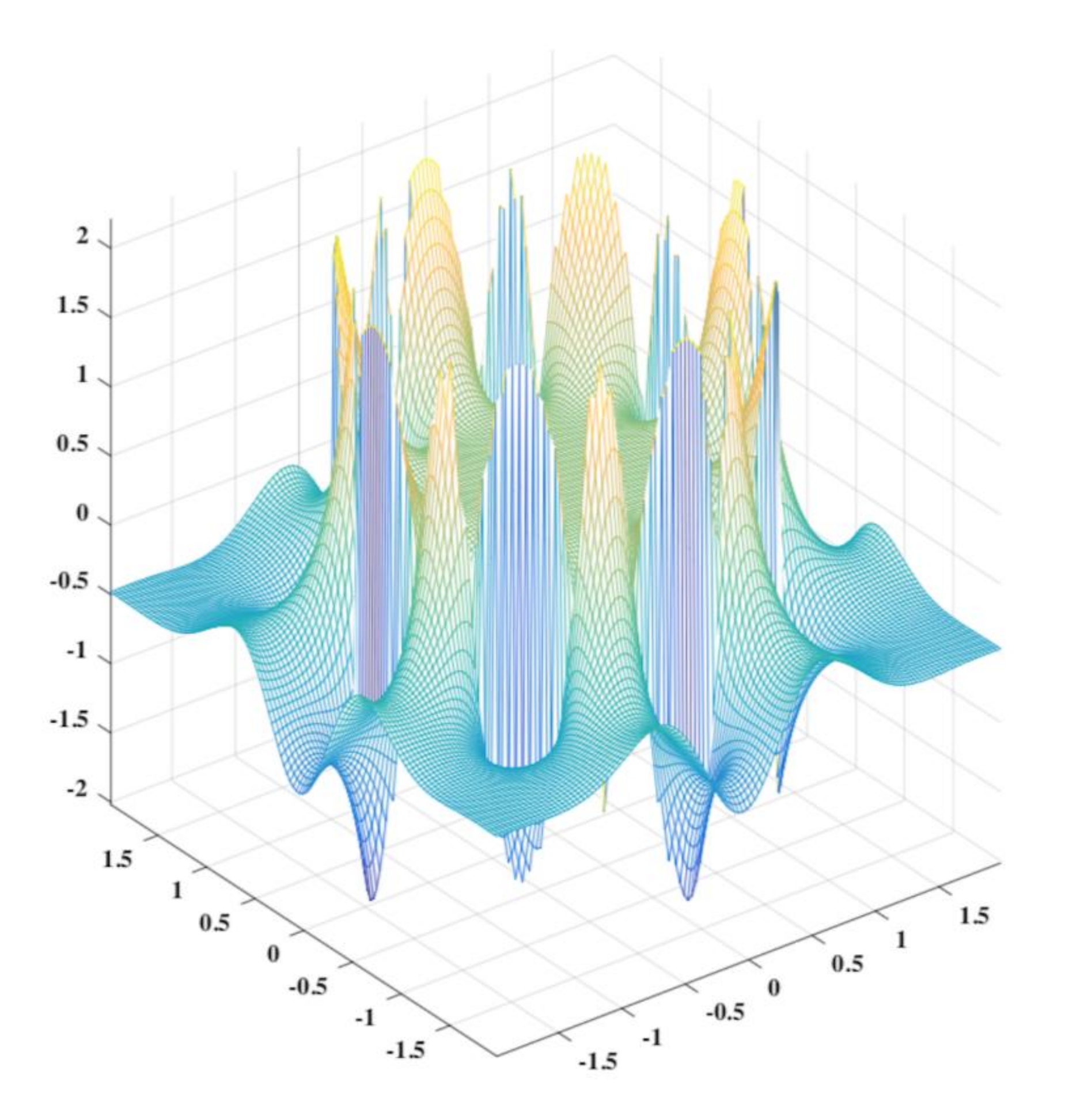}
}
\subfigure[$t=1.87$]{
\includegraphics[width=0.31\textwidth]{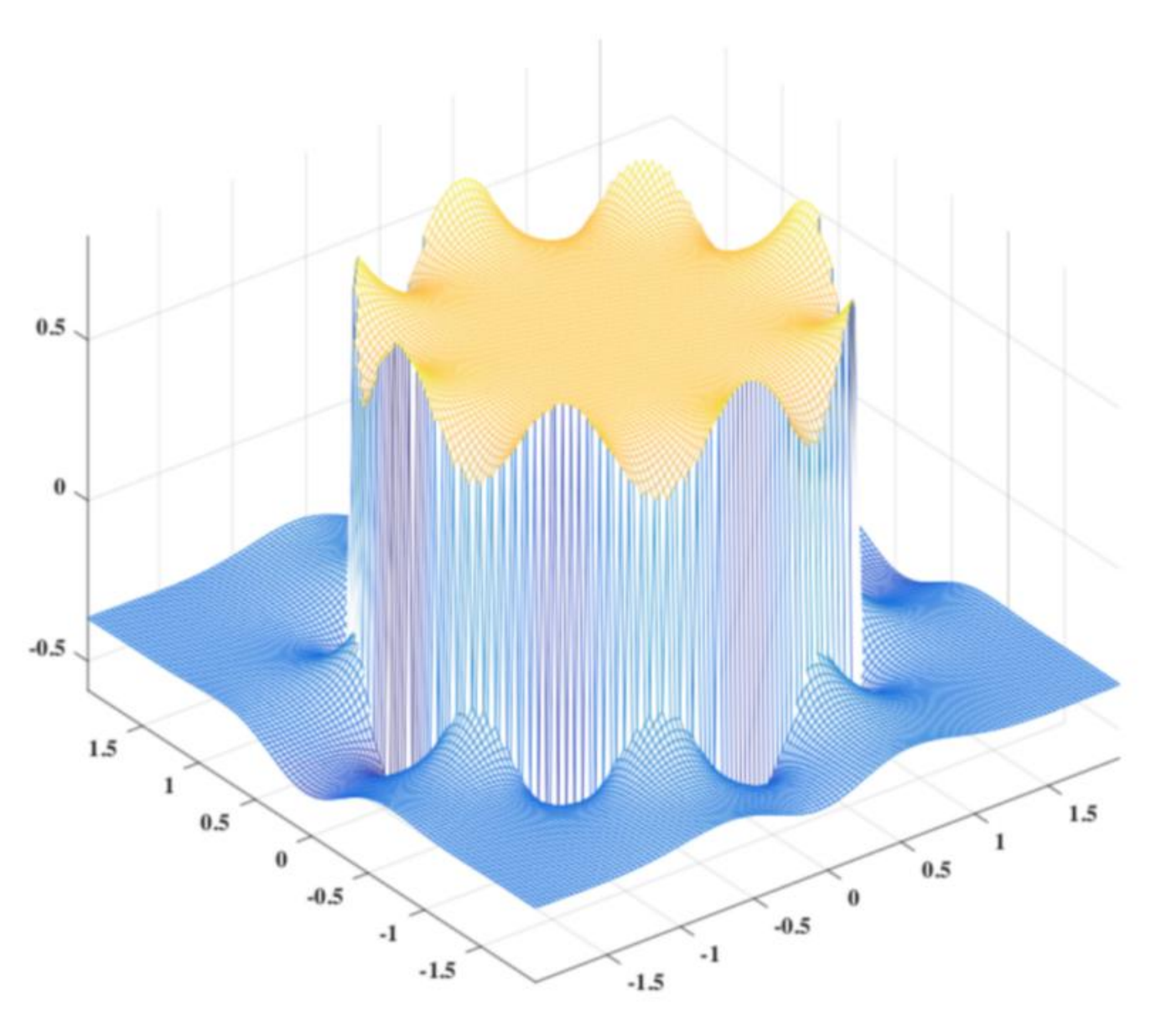}
}
\setlength{\abovecaptionskip}{-0.0cm}
\setlength{\belowcaptionskip}{-0.0cm}
 \caption{Example 4. The pressure distribution at different times.}
 \label{pressure-e2}
\end{figure}

%========================================================================
{\em Example} 4.  This example is to show that the proposed method can handle flows with a more complicated initial interface, which is given by $r=0.8+0.2\sin(8\theta), 0\leq \theta\leq2\pi$ in polar coordinates.  Fig. \ref{interface-1} (right) shows the corresponding  interface configurations at different times $t=0, t=0.187, t=0.374, t=0.748, t=1.87$. The tension coefficient $T_0$ is set to be 0.5.  In this example, the viscosity outside the interface is larger than inside the interface, which is taken as $\mu^+ = 1$ and $\mu^- = 10$.  The approximation is computed up to a final time $T=2$ with $\Delta t = h = 0.0187$.  A time evolution of the velocity is presented in Fig. \ref{velocity-e2} and  a time evolution of the pressure profile  is shown in Fig. \ref{pressure-e2}. From these figures, it is clear  that  the velocity is continuous but not smooth, while the pressure is discontinuous across the interface, suggesting that the proposed method can capture highly discontinuous profile for the pressure in each time step.

\begin{figure}[h!]
\centering
\subfigure{
%\subfigure[$u^{(1)}$]{
\includegraphics[width=0.31\textwidth]{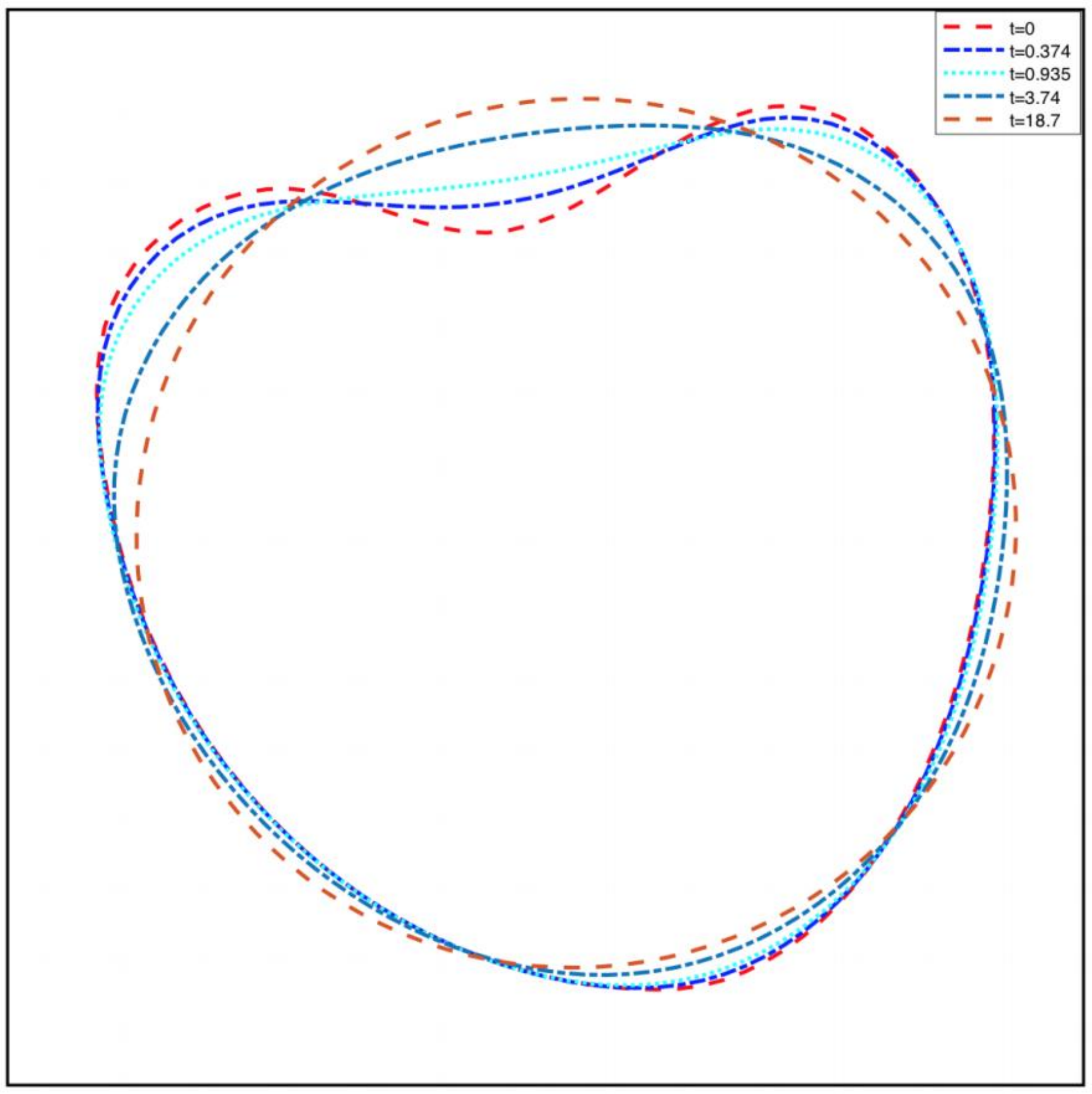}
}
%\subfigure[$u^{(2)}$]{
\subfigure{
\includegraphics[width=0.31\textwidth]{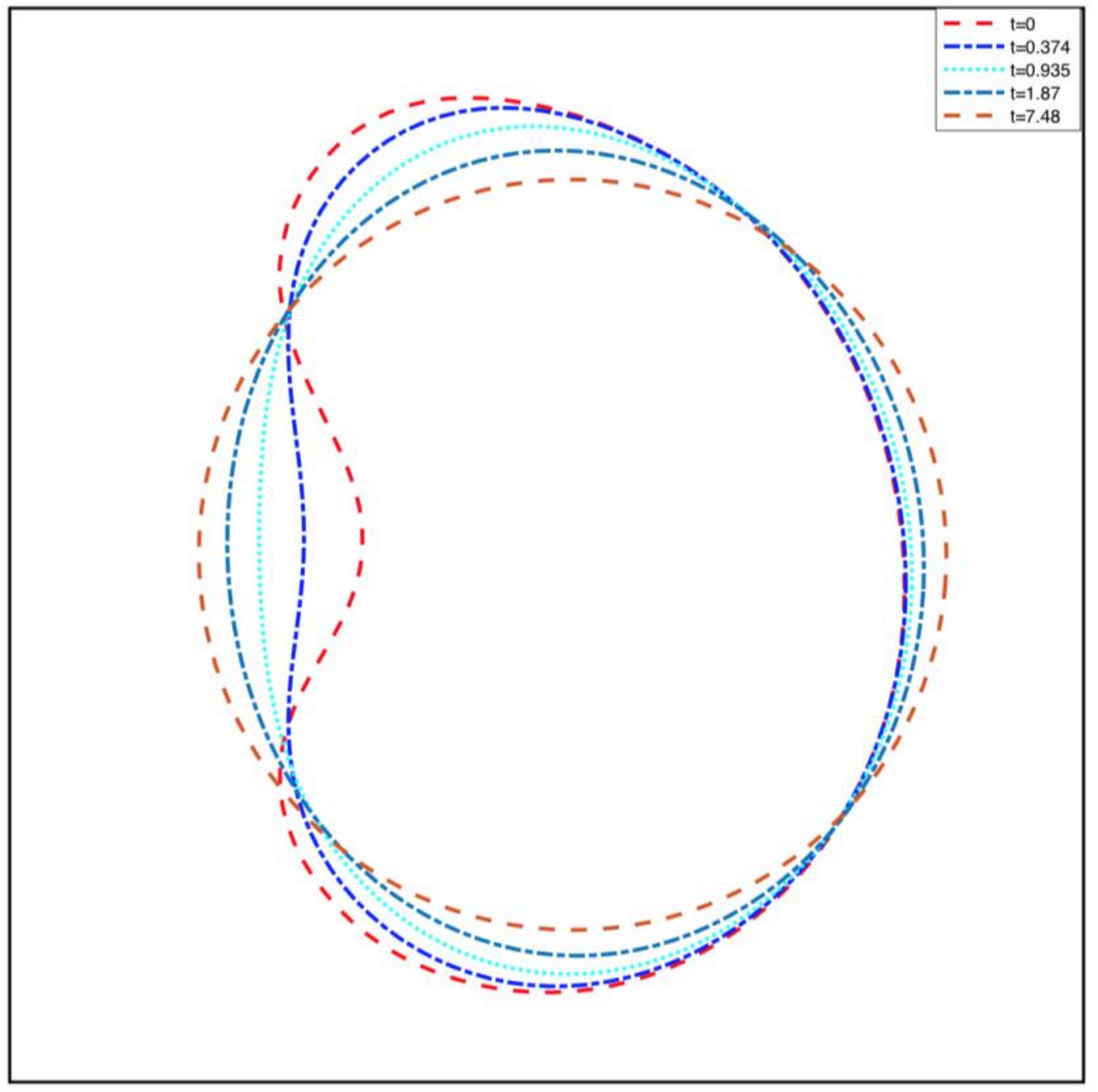}
}
\setlength{\abovecaptionskip}{-0.0cm}
\setlength{\belowcaptionskip}{-0.0cm}
 \caption{The interface configurations at different times in a square domain. (left: Heart-shaped initial  interface; right: Kidney-shaped initial  interface) }
\label{interface-2}
\end{figure}

\begin{figure}[h!]
\centering
\subfigure[$$]{
\includegraphics[width=0.31\textwidth]{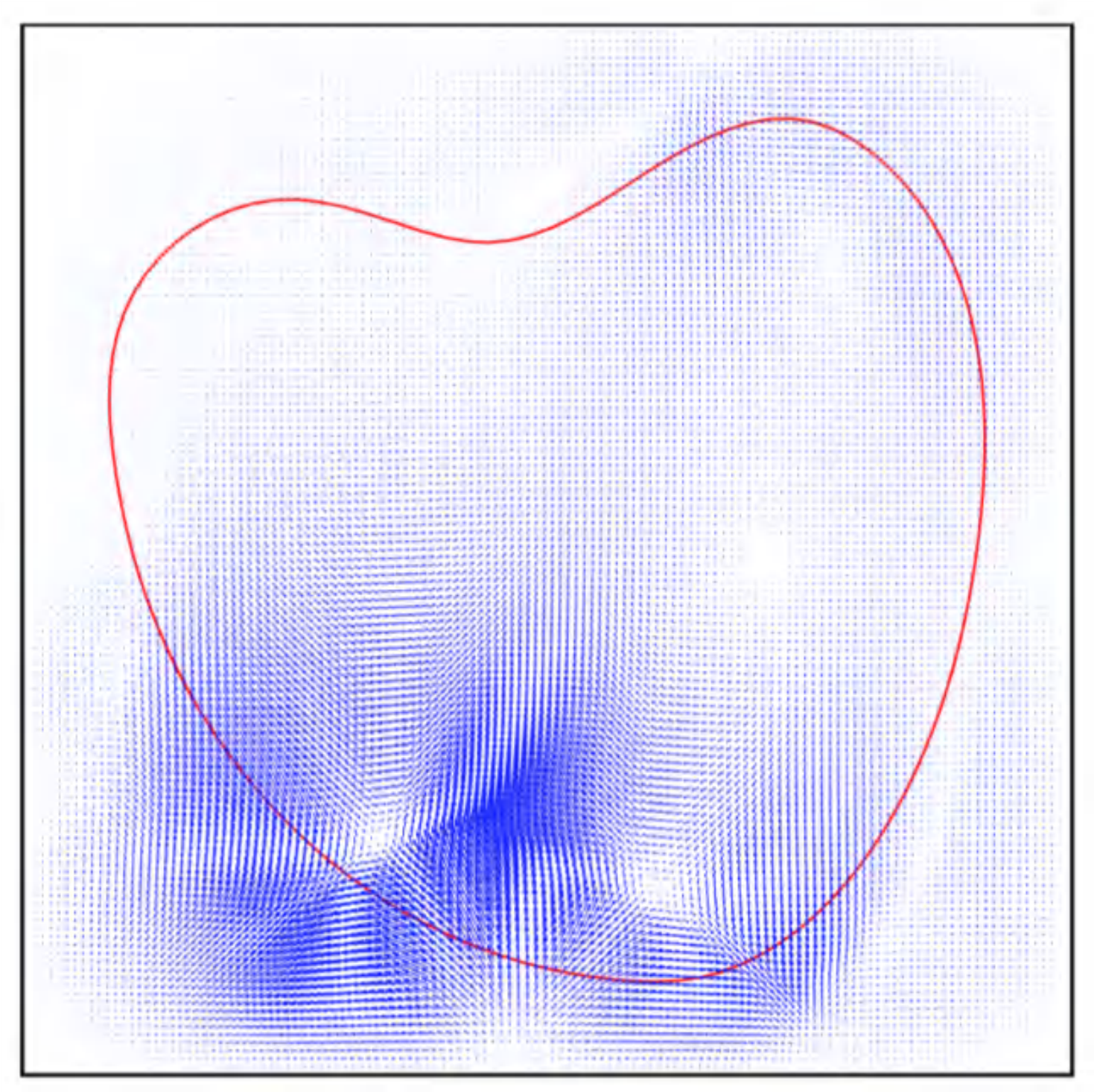}
}
\subfigure[$$]{
\includegraphics[width=0.31\textwidth]{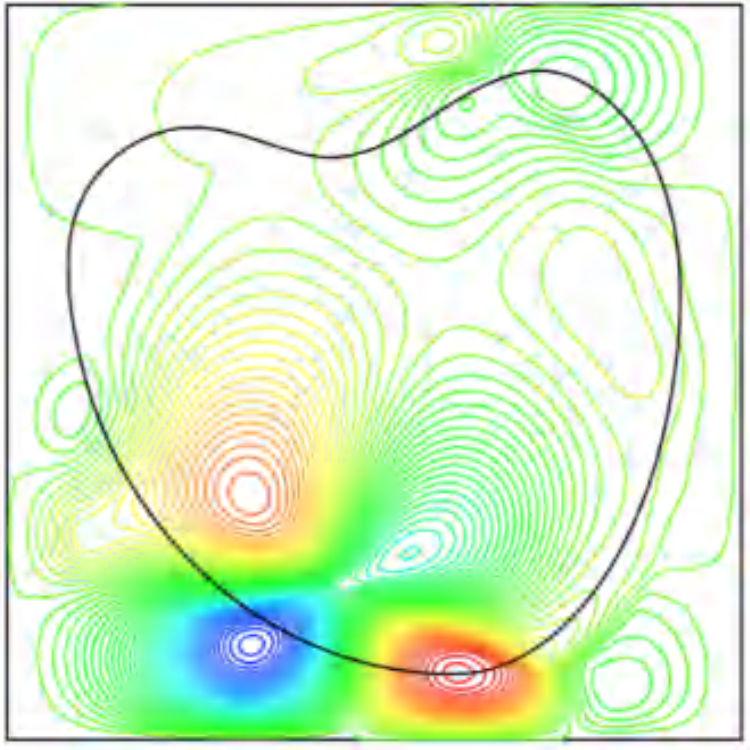}
}
\subfigure[$$]{
\includegraphics[width=0.31\textwidth]{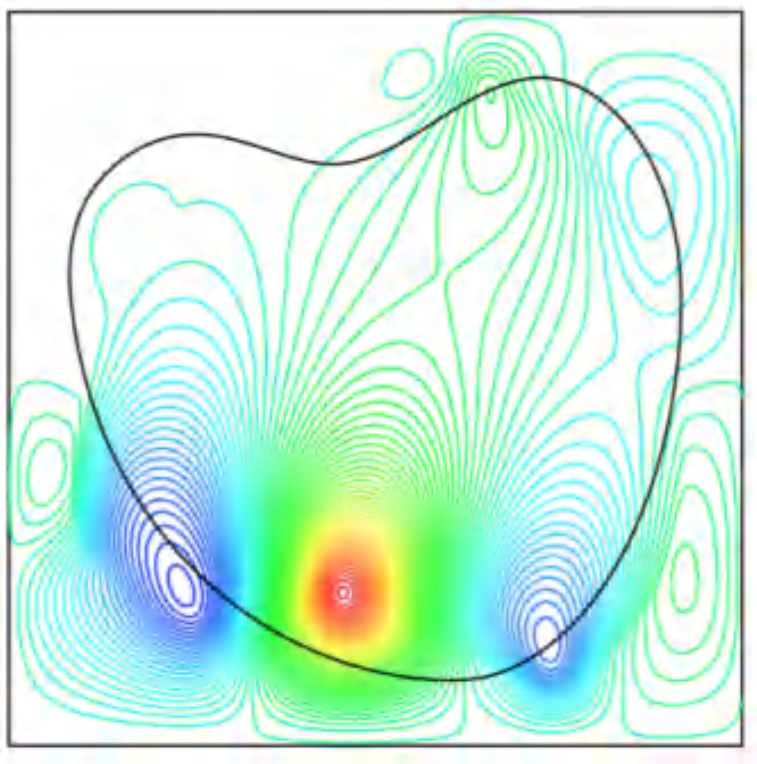}
}
\subfigure[$$]{
\includegraphics[width=0.31\textwidth]{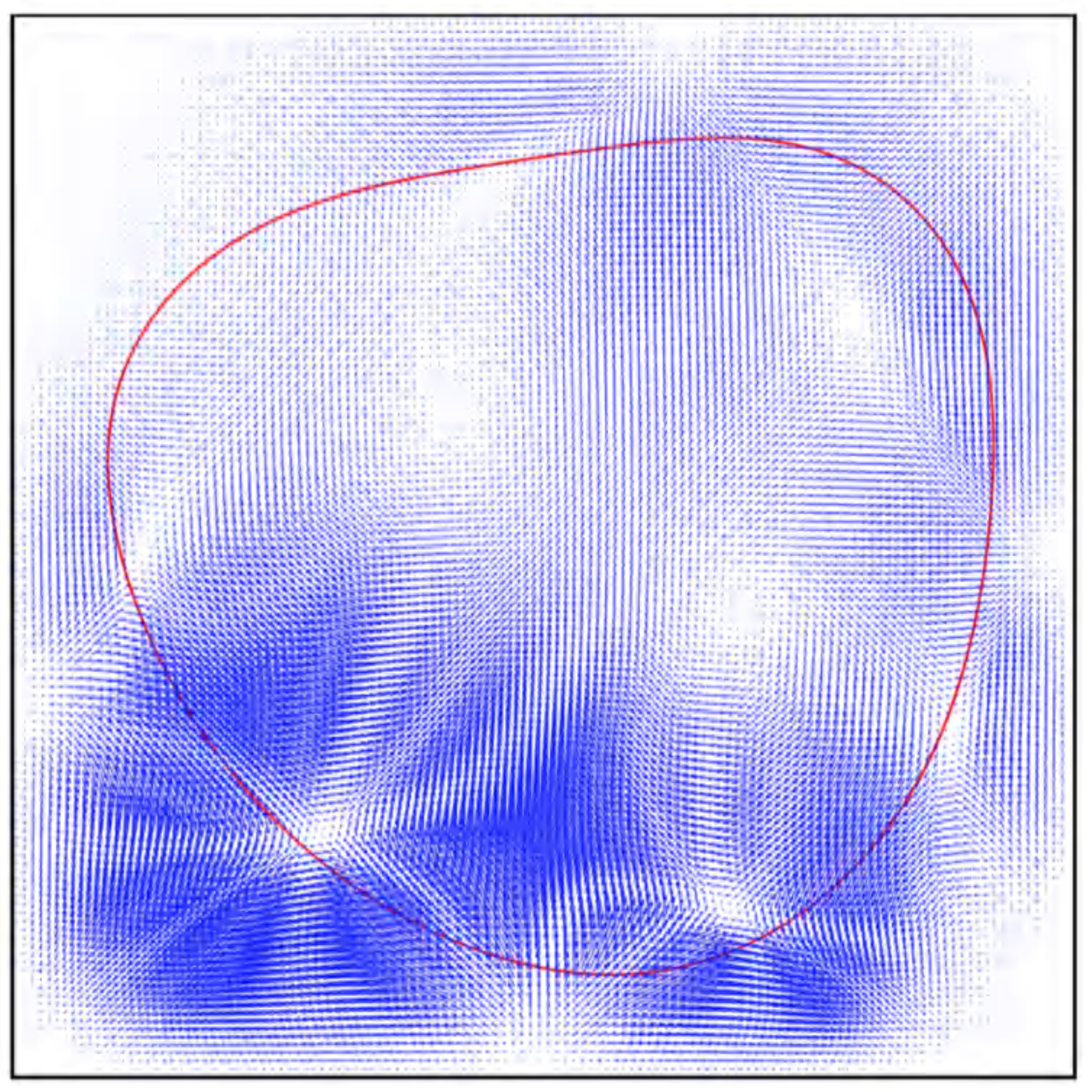}
}
\subfigure[$$]{
\includegraphics[width=0.31\textwidth]{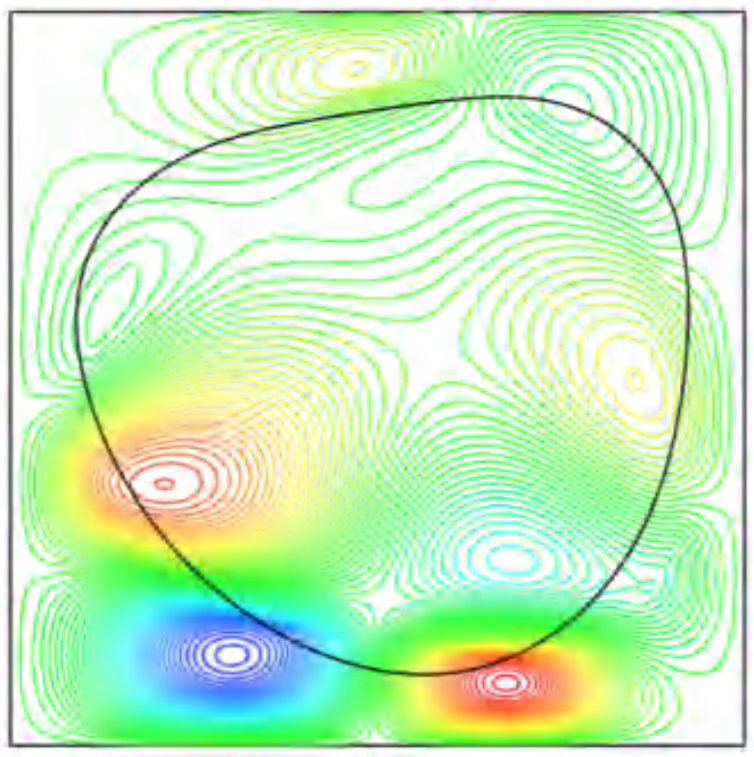}
}
\subfigure[$$]{
\includegraphics[width=0.31\textwidth]{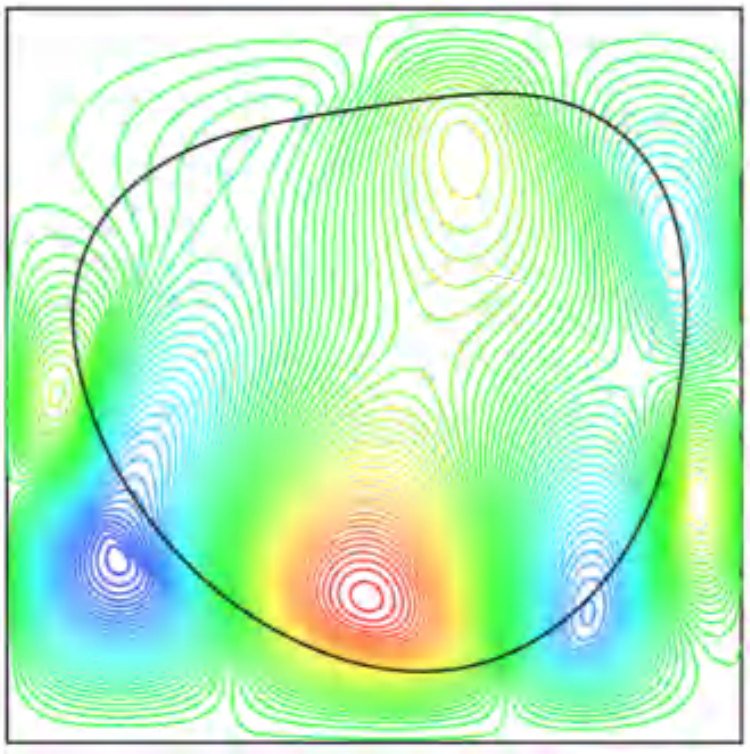}
}
\subfigure[$$]{
\includegraphics[width=0.31\textwidth]{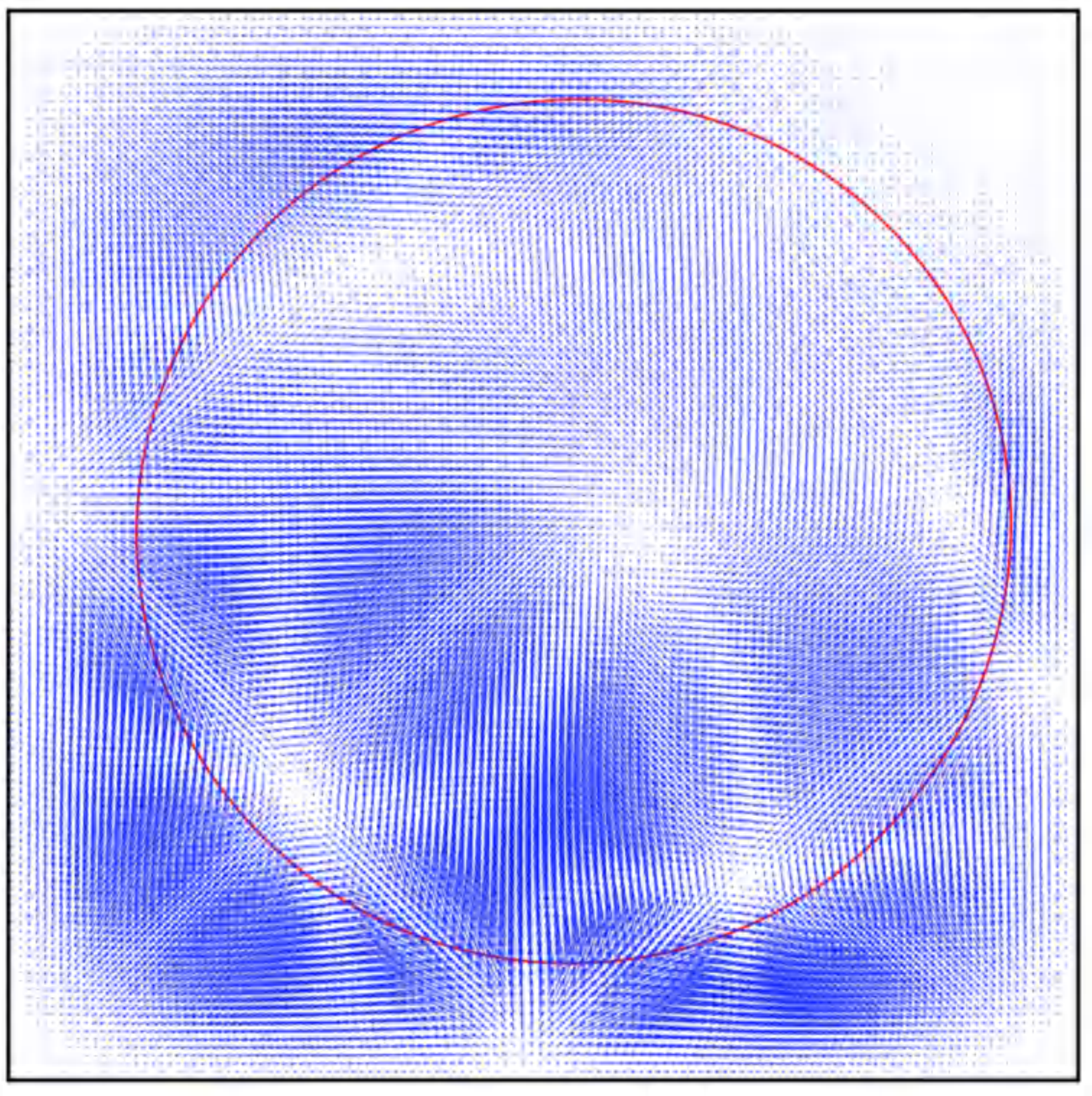}
}
\subfigure[$$]{
\includegraphics[width=0.31\textwidth]{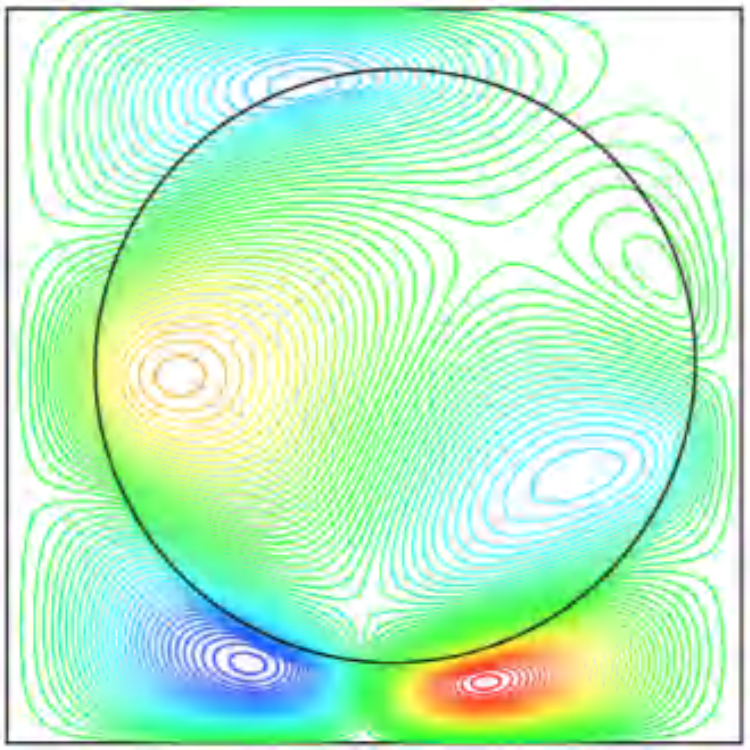}
}
\subfigure[$$]{
\includegraphics[width=0.31\textwidth]{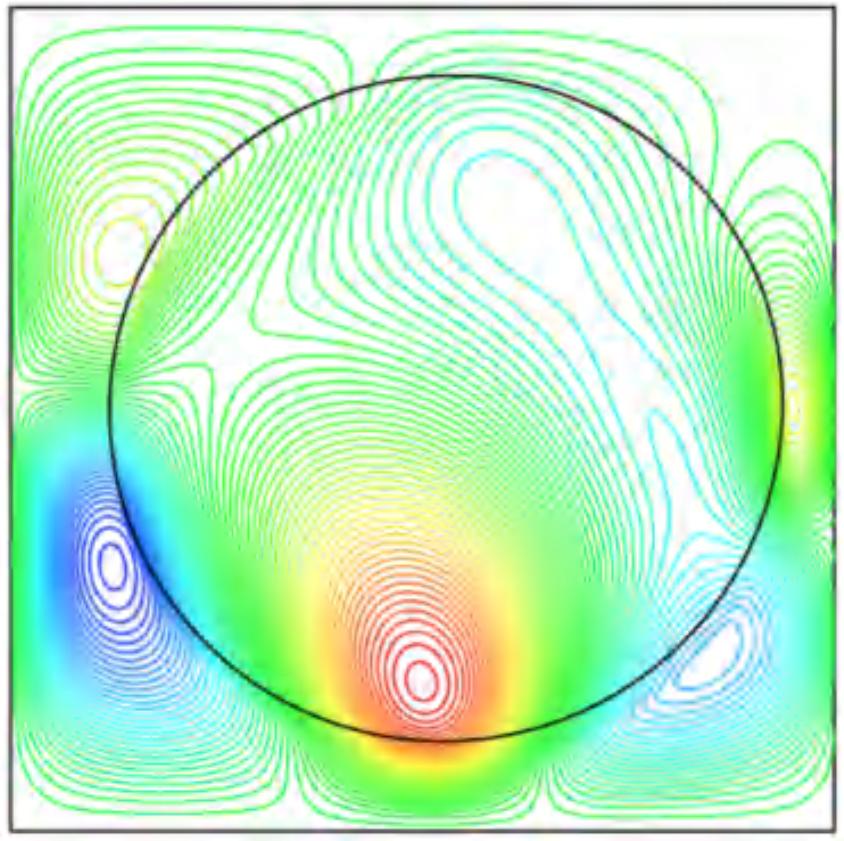}
}
\setlength{\abovecaptionskip}{-0.0cm}
\setlength{\belowcaptionskip}{-0.0cm}
 \caption{Example 5. Evolution of the velocity field $\vect u$ and interface position for heart-shaped initial  interface problem, computed with a $128\times 128$ grid and $\Delta t = h = 0.0187$. (left: Velocity field $\vect u$; middle: Isolines of the $x$-component $u^{(1)}$; right: Isolines of the $y$-component  $u^{(2)}$.)}
 \label{velocity-e3}
\end{figure}

\begin{figure}[h!]
\centering
\subfigure[$t=0$]{
\includegraphics[width=0.31\textwidth]{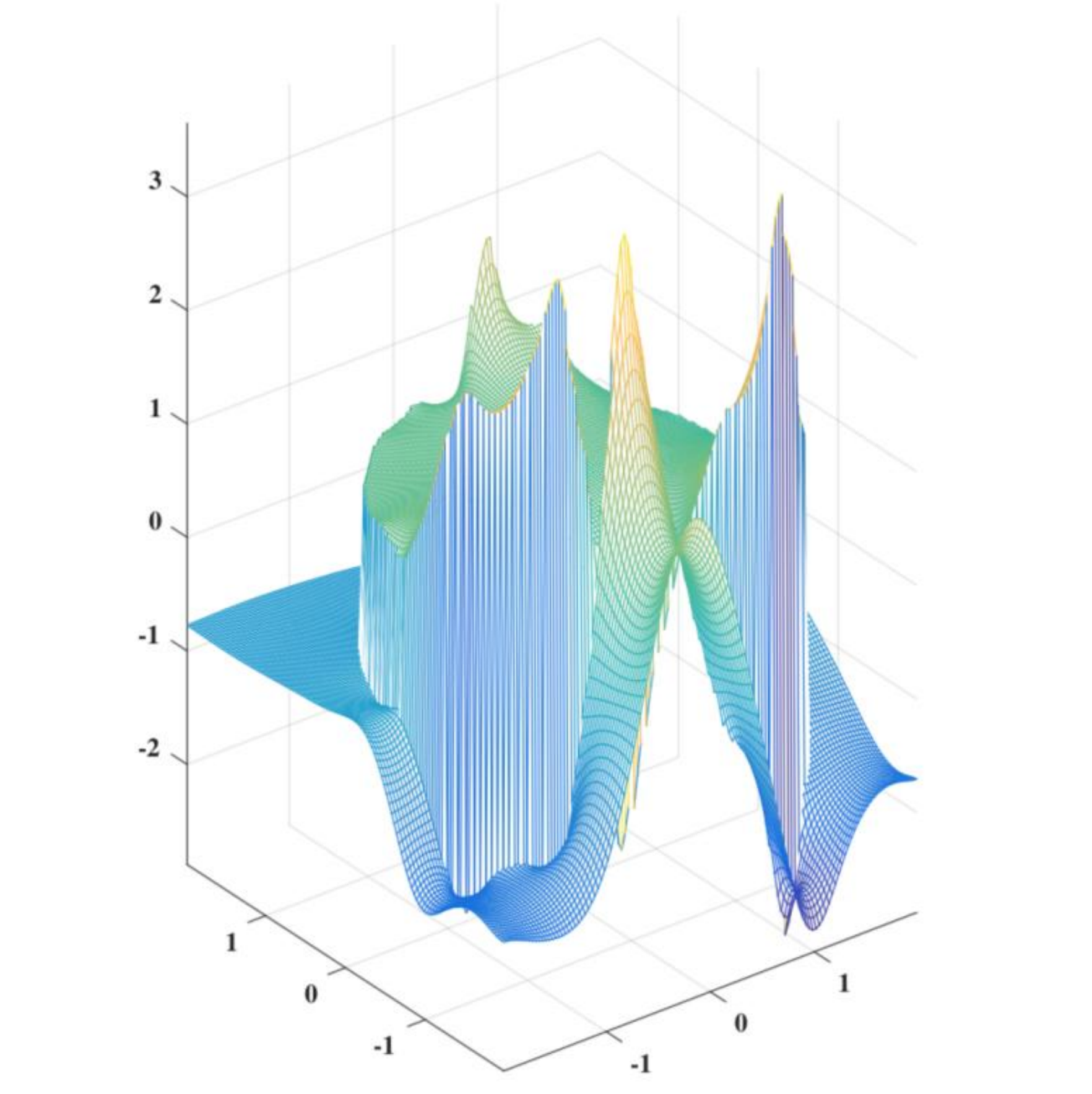}
}
%\subfigure[$t=0.374$]{
%\includegraphics[width=0.23\textwidth]{m_q_t4_e3.pdf}
%}
\subfigure[$t=1.87$]{
\includegraphics[width=0.31\textwidth]{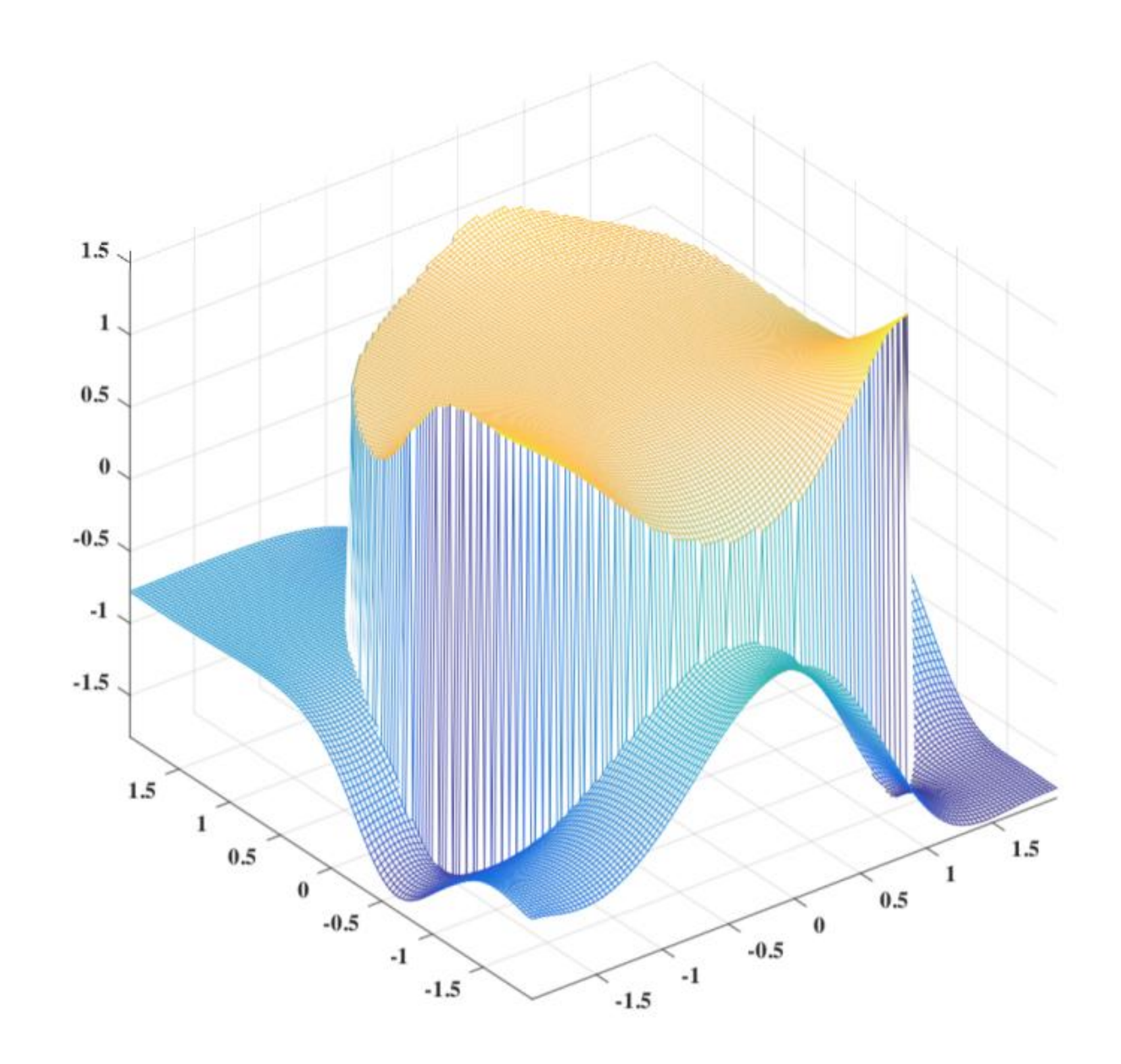}
}
\subfigure[$t=3.74$]{
\includegraphics[width=0.31\textwidth]{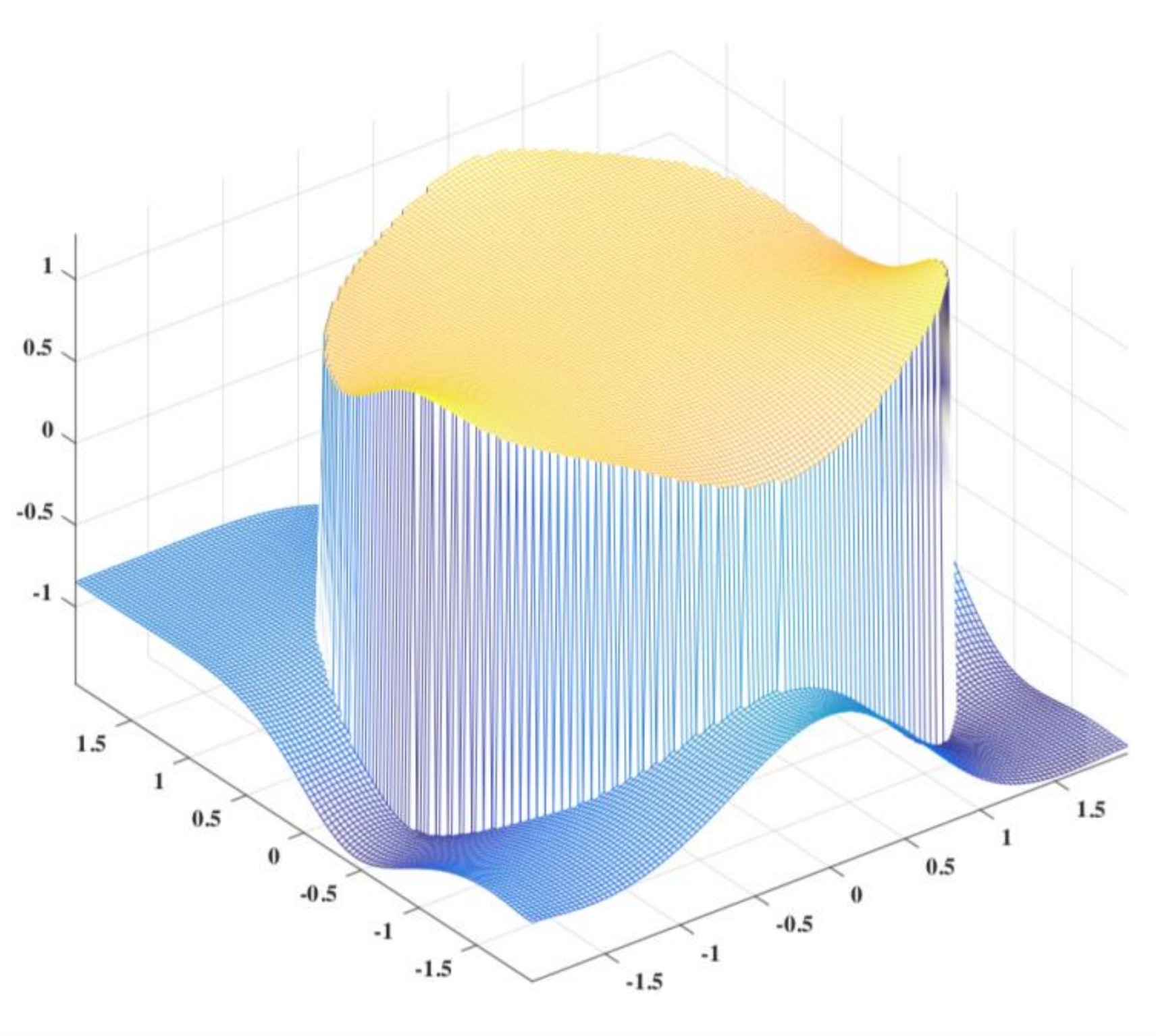}
}
\setlength{\abovecaptionskip}{-0.0cm}
\setlength{\belowcaptionskip}{-0.0cm}
 \caption{Example 5. The pressure distribution at different times.}
 \label{pressure-e3}
\end{figure}

{\em Example }5.  In this example, an initial heart-shaped interface given by $100$ control points is considered.  The interface configurations at different times ($t=0, t=0.374, t=0.935, t=3.74, t=18.7$) are shown in Fig. \ref{interface-2} (left). As we can see, the heart-shaped interface moves according to the velocity field while slowly transforming into a circle.   The tension coefficient $T_0$ is set to be 1.0.  In this example, the viscosity outside the interface is larger than inside the interface, which is taken as $\mu^+ = 1$ and $\mu^- = 10$.  The approximations is computed up to a final time $T=20$ with $\Delta t = h = 0.0187$.  
A time evolution of the velocity is plotted in Fig.\ref{velocity-e3} and a time evolution of the pressure profile is shown in Fig. \ref{pressure-e3}. Similar to the flower-shaped case, the velocity is continuous but not smooth, while the pressure is discontinuous across the interface, demonstrating  that the proposed method can handle more involved interface configurations, which are non-convex in this case.

\begin{figure}[h!]
\centering
\subfigure[$$]{
\includegraphics[width=0.31\textwidth]{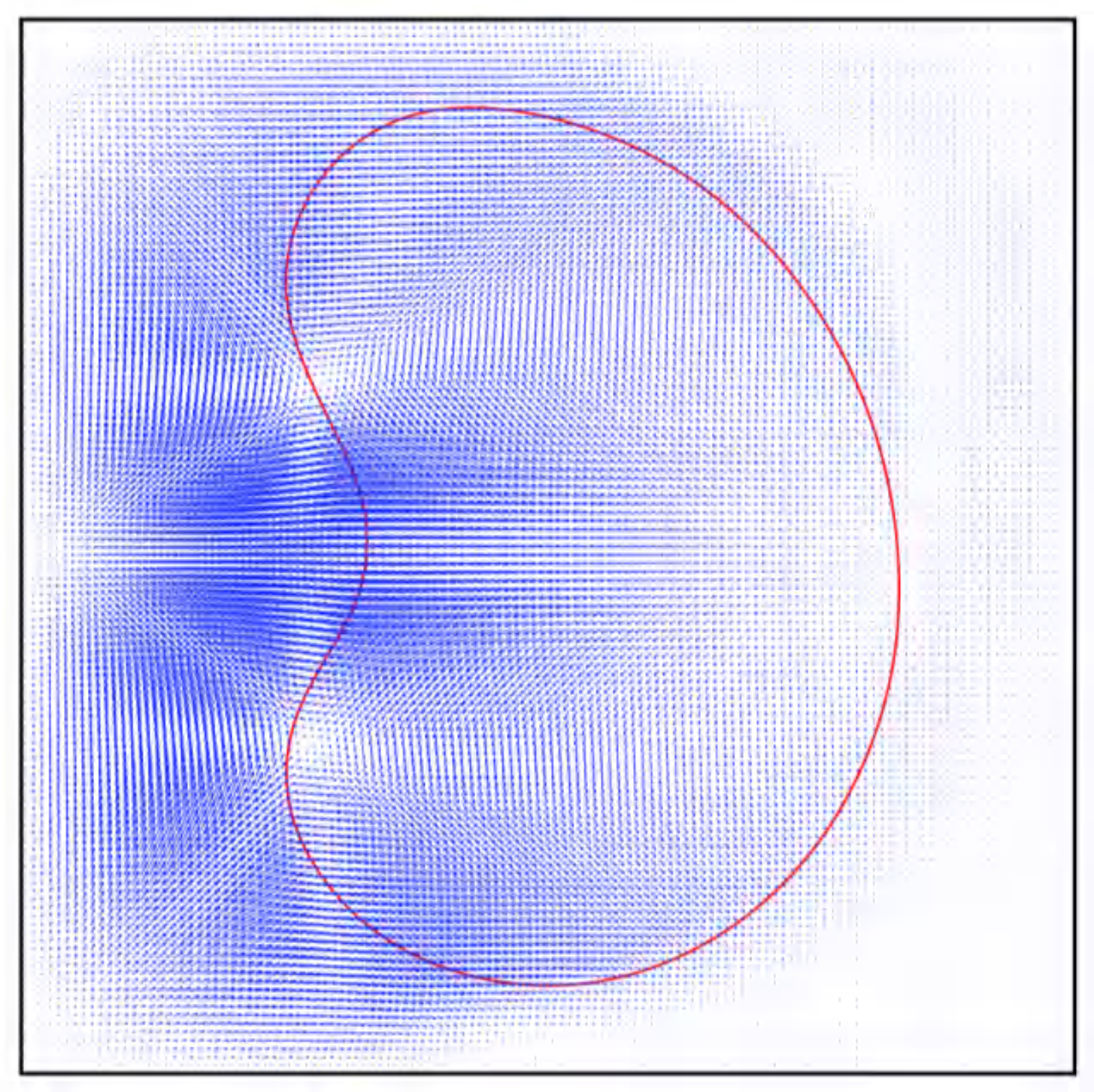}
}
\subfigure[$$]{
\includegraphics[width=0.31\textwidth]{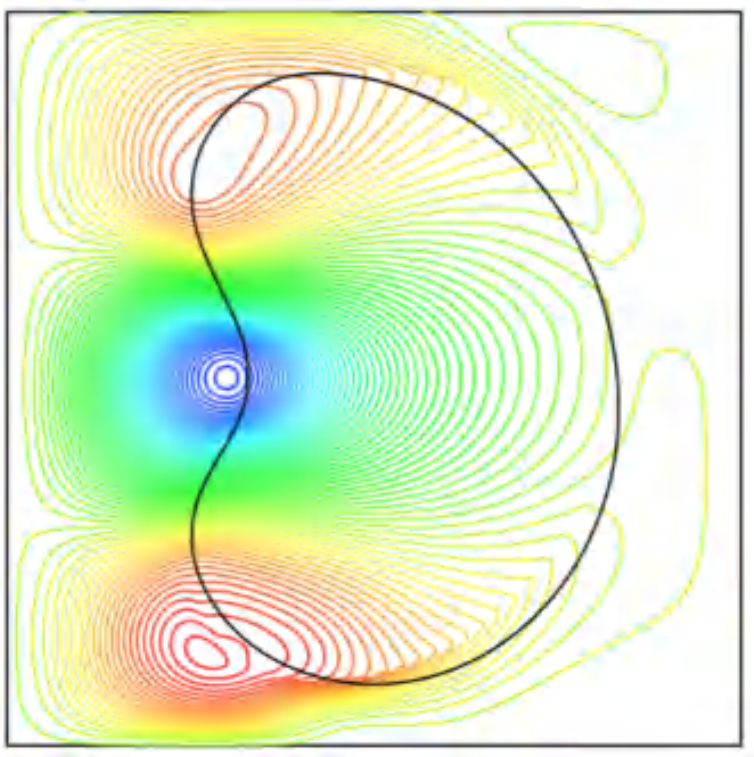}
}
\subfigure[$$]{
\includegraphics[width=0.31\textwidth]{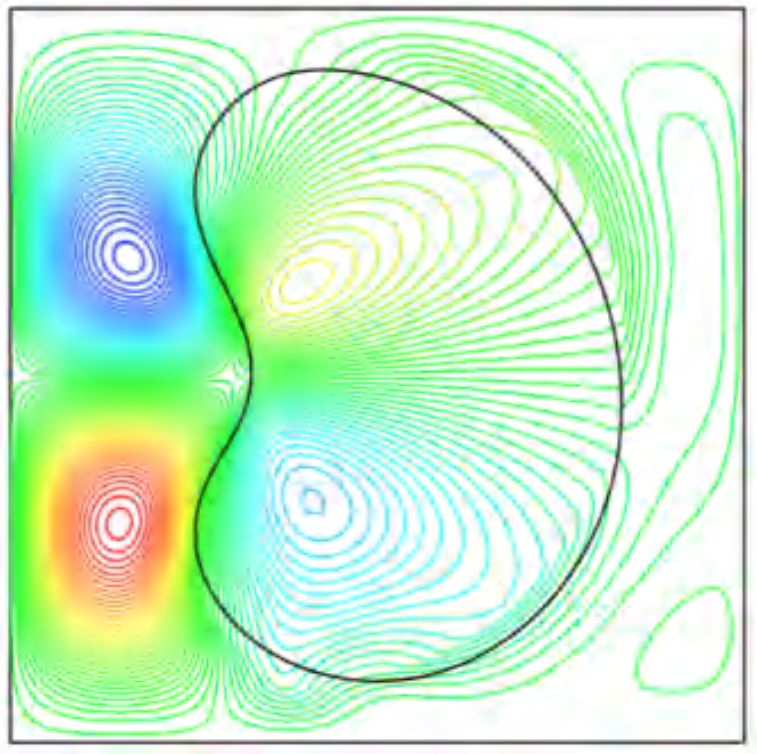}
}
\subfigure[$$]{
\includegraphics[width=0.31\textwidth]{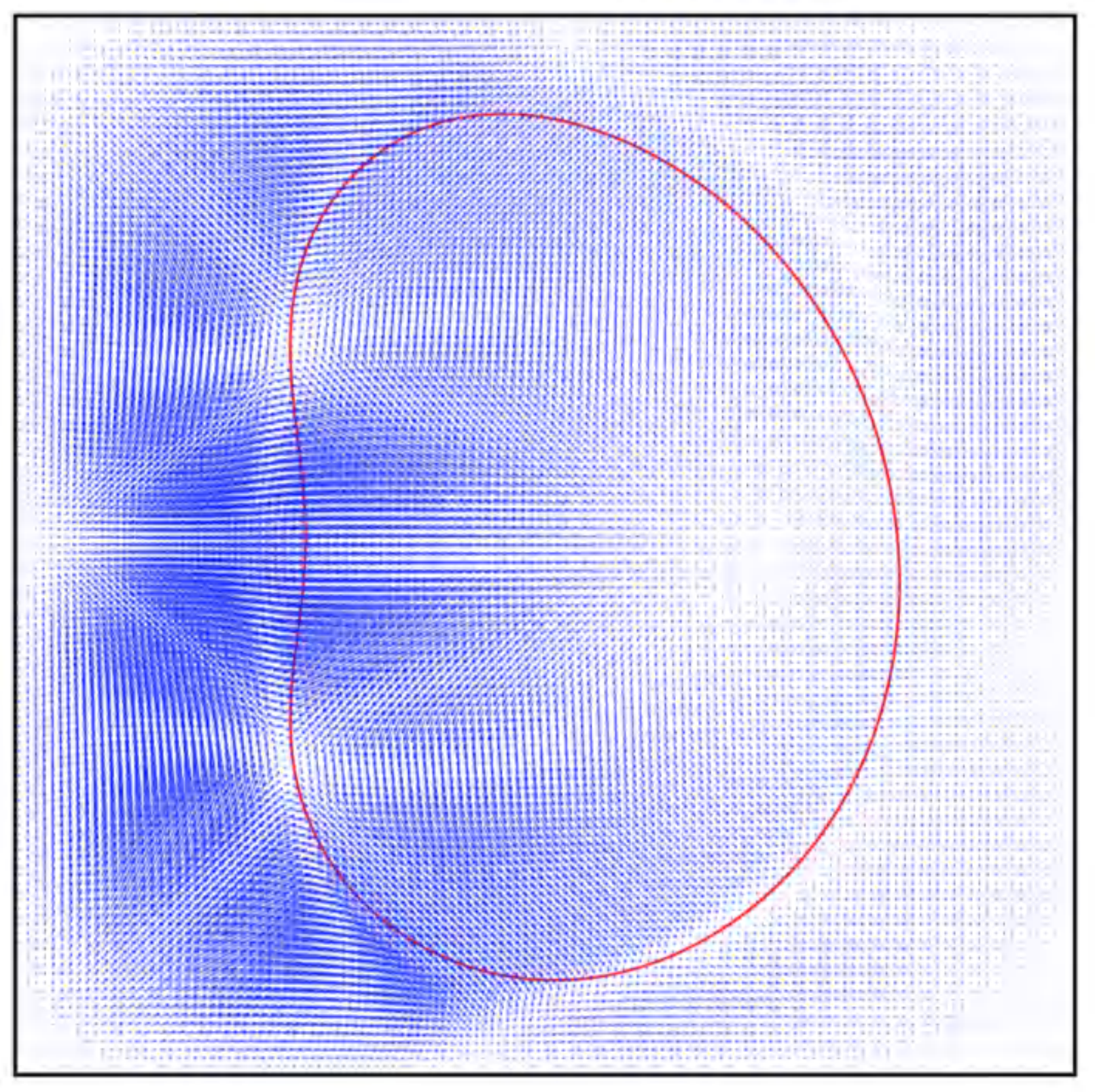}
}
\subfigure[$$]{
\includegraphics[width=0.31\textwidth]{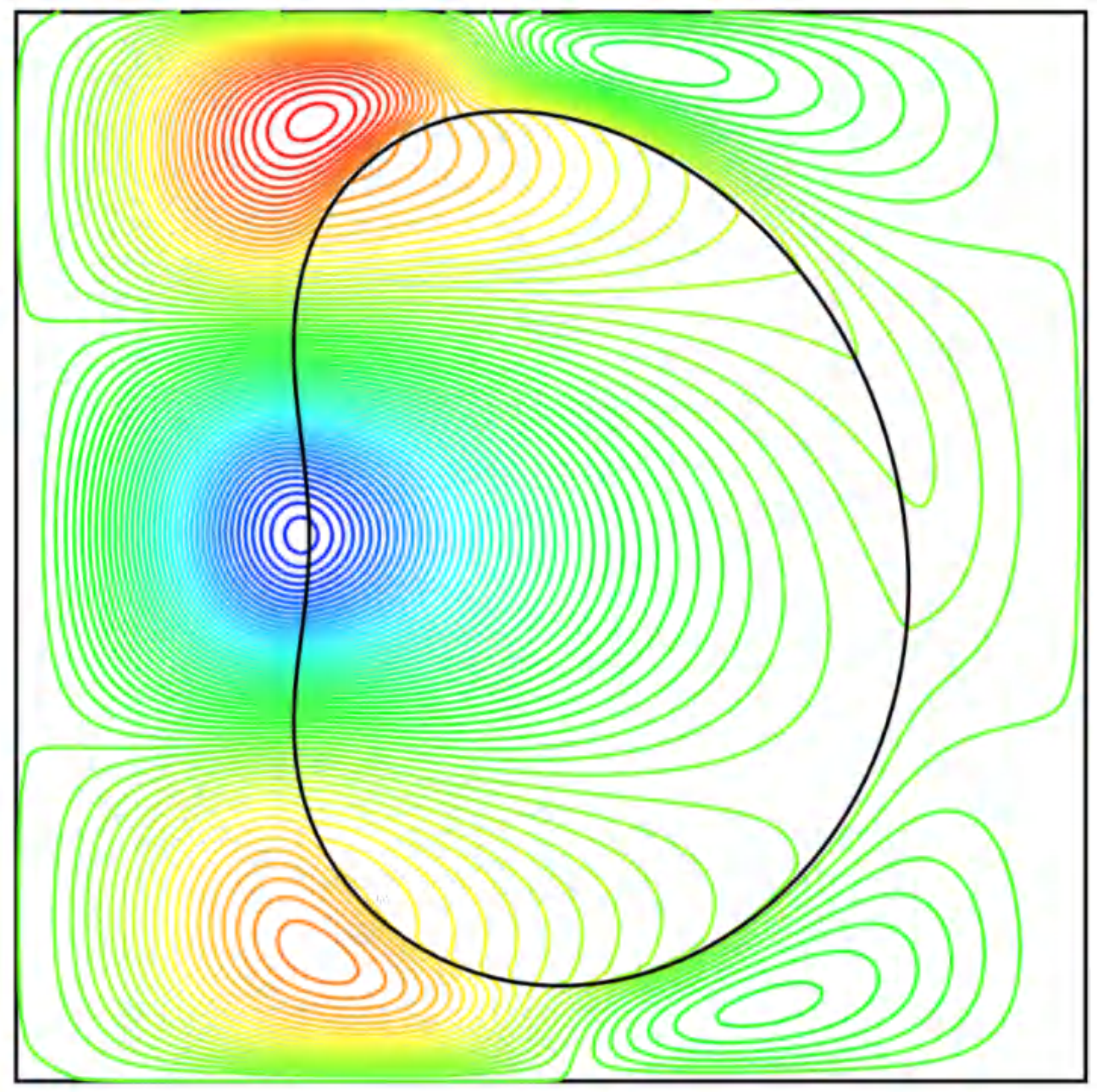}
}
\subfigure[$$]{
\includegraphics[width=0.31\textwidth]{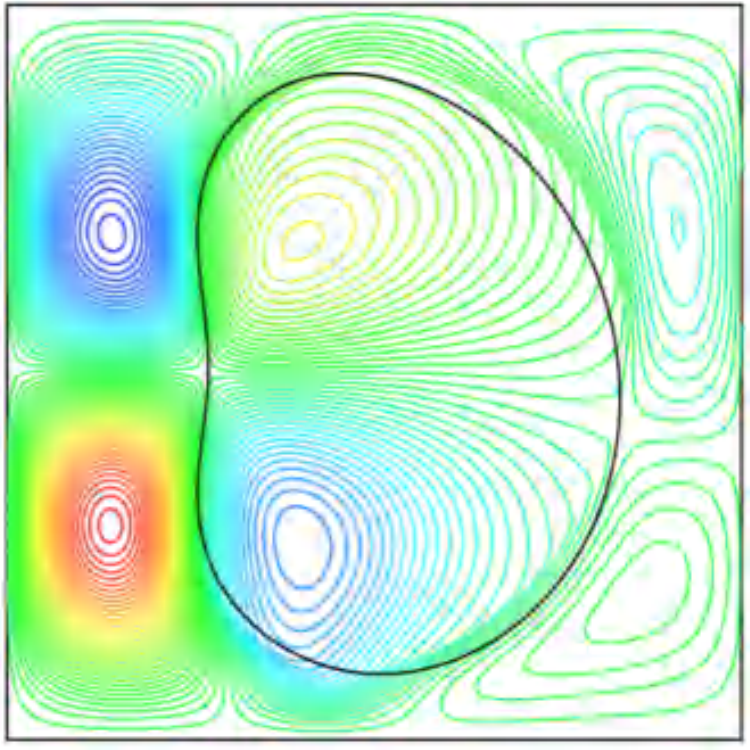}
}
\subfigure[$$]{
\includegraphics[width=0.31\textwidth]{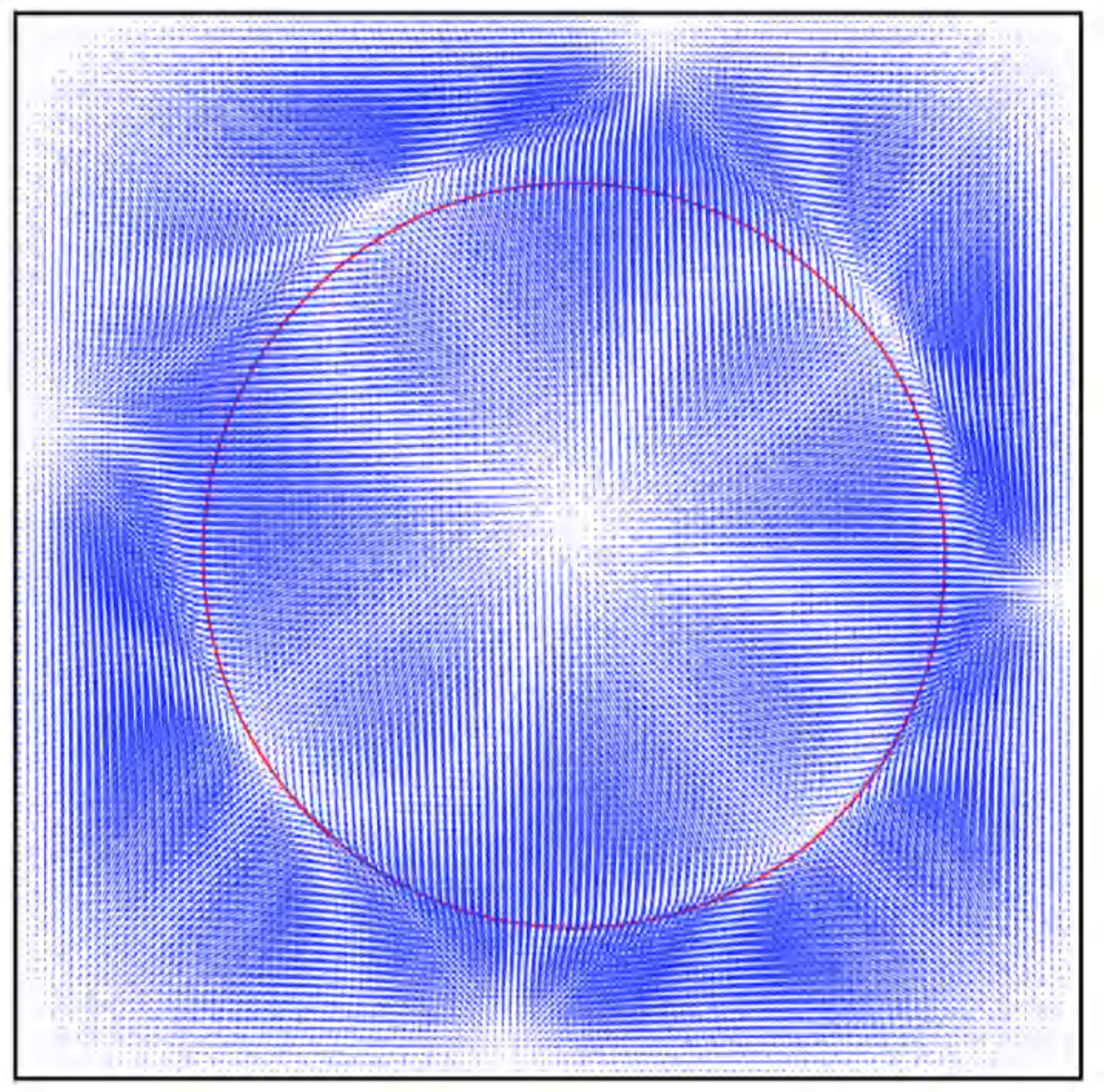}
}
\subfigure[$$]{
\includegraphics[width=0.31\textwidth]{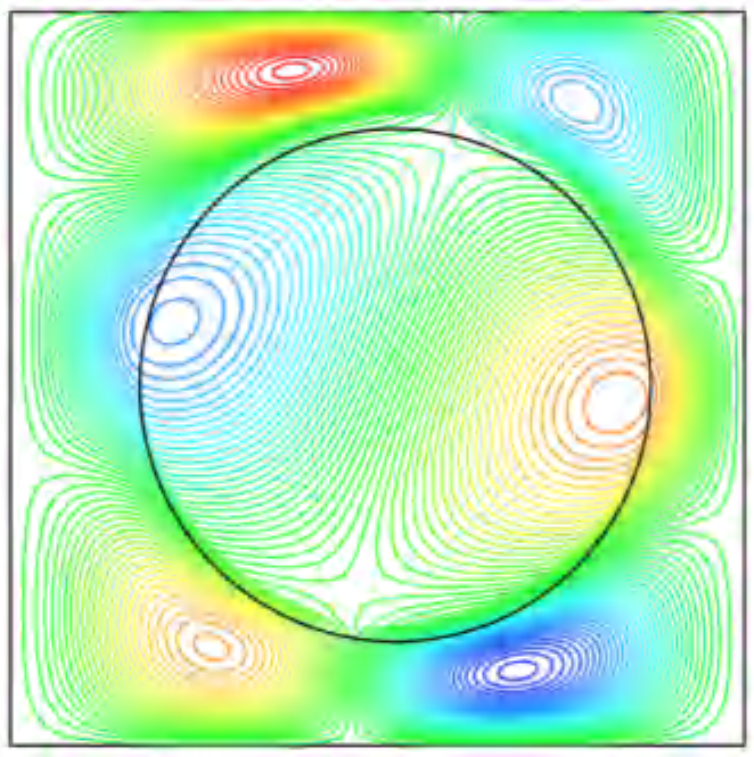}
}
\subfigure[$$]{
\includegraphics[width=0.31\textwidth]{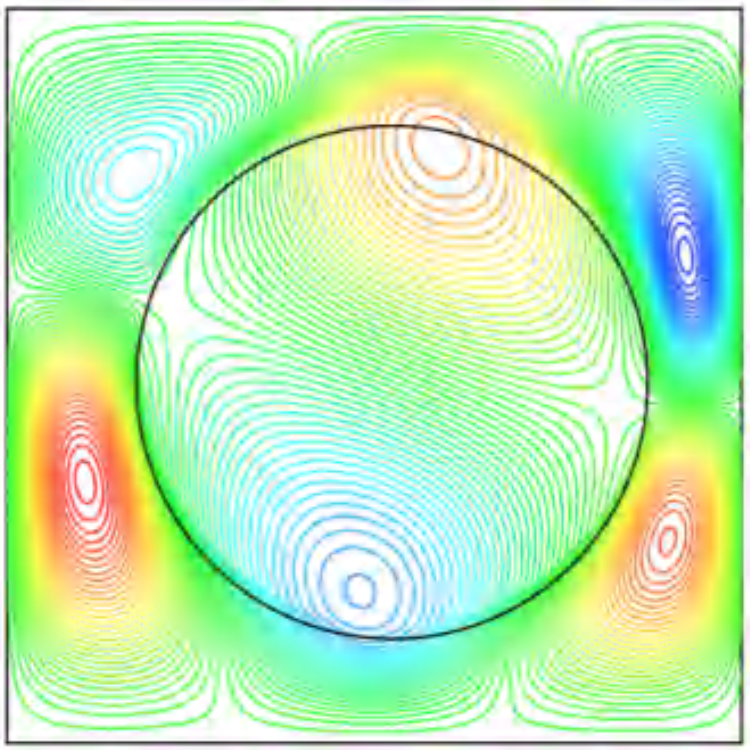}
}
\setlength{\abovecaptionskip}{-0.0cm}
\setlength{\belowcaptionskip}{-0.0cm}
 \caption{Example 6. Evolution of the velocity field $\vect u$ and interface position for kidney-shaped initial  interface problem, computed with a $128\times 128$ grid and $\Delta t = h = 0.0187$. (left: Velocity field $\vect u$; middle: Isolines of the $x$-component $u^{(1)}$; right: Isolines of the $y$-component  $u^{(2)}$.)}
 \label{velocity-e4}
\end{figure}

\begin{figure}[h!]
\centering
\subfigure[$t=0$]{
\includegraphics[width=0.31\textwidth]{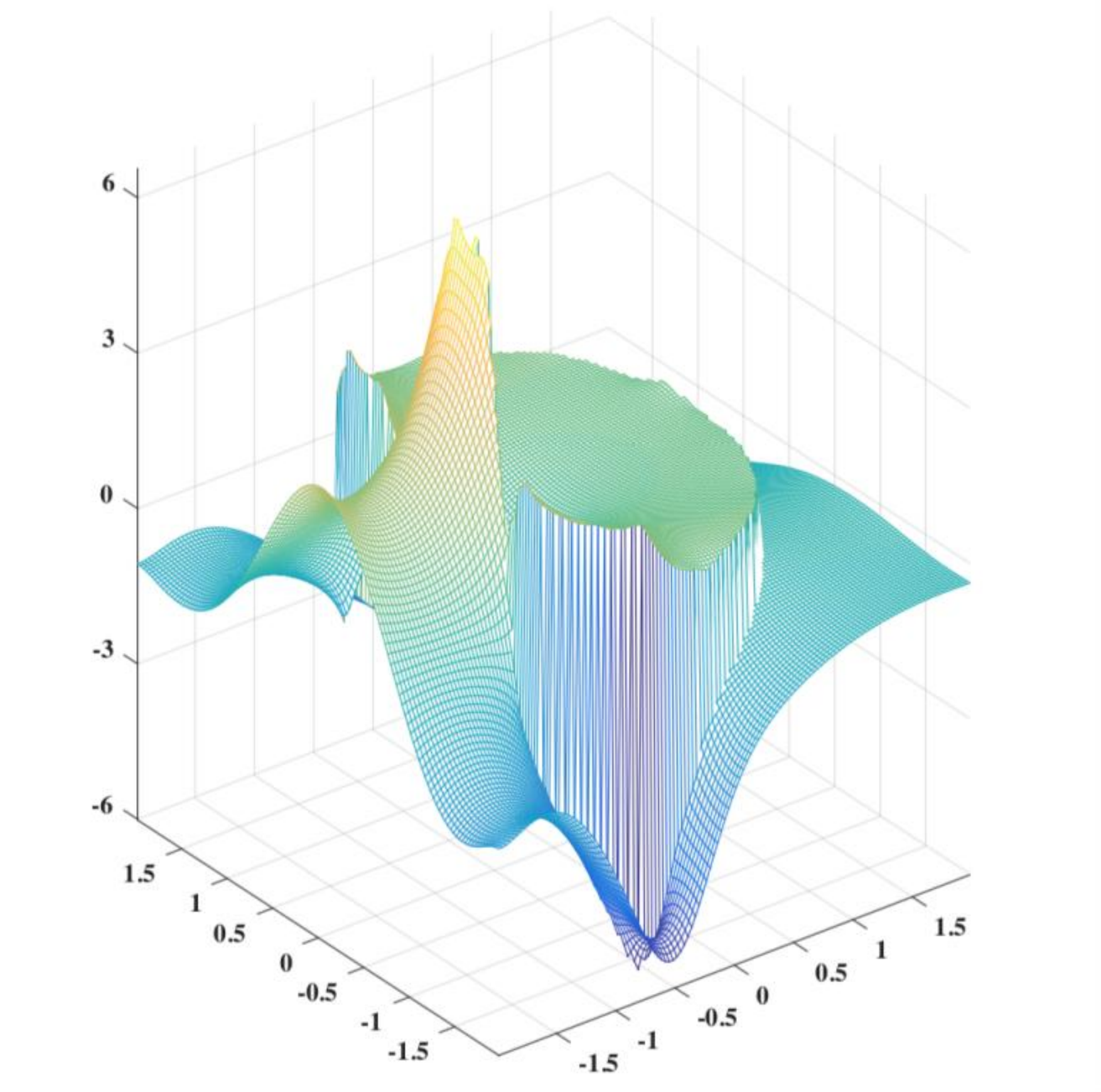}
}
\subfigure[$t=0.748$]{
\includegraphics[width=0.31\textwidth]{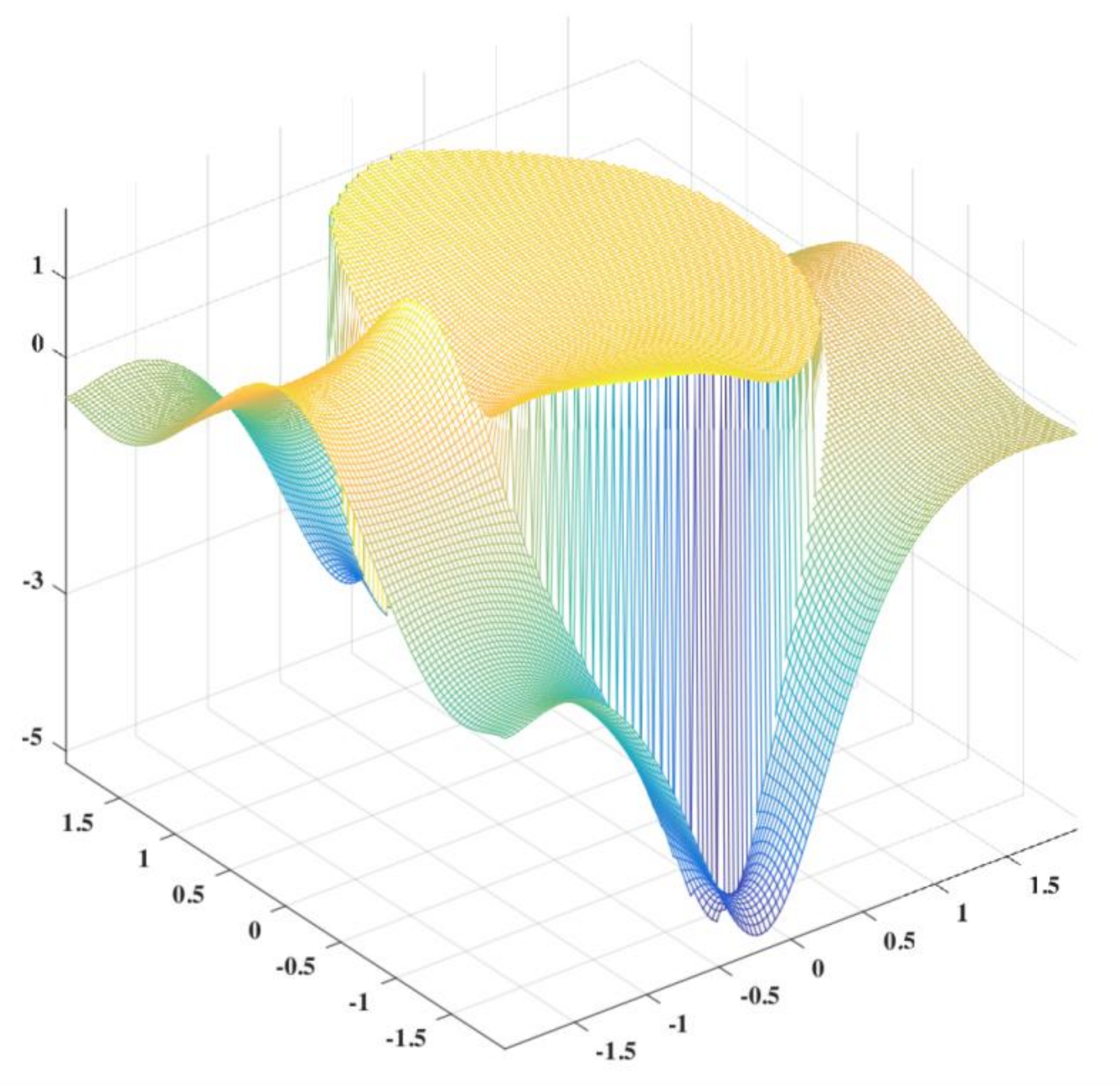}
}
%\subfigure[$t=0.935$]{
%\includegraphics[width=0.23\textwidth]{m_q_t6_e4.pdf}
%}
\subfigure[$t=3.74$]{
\includegraphics[width=0.31\textwidth]{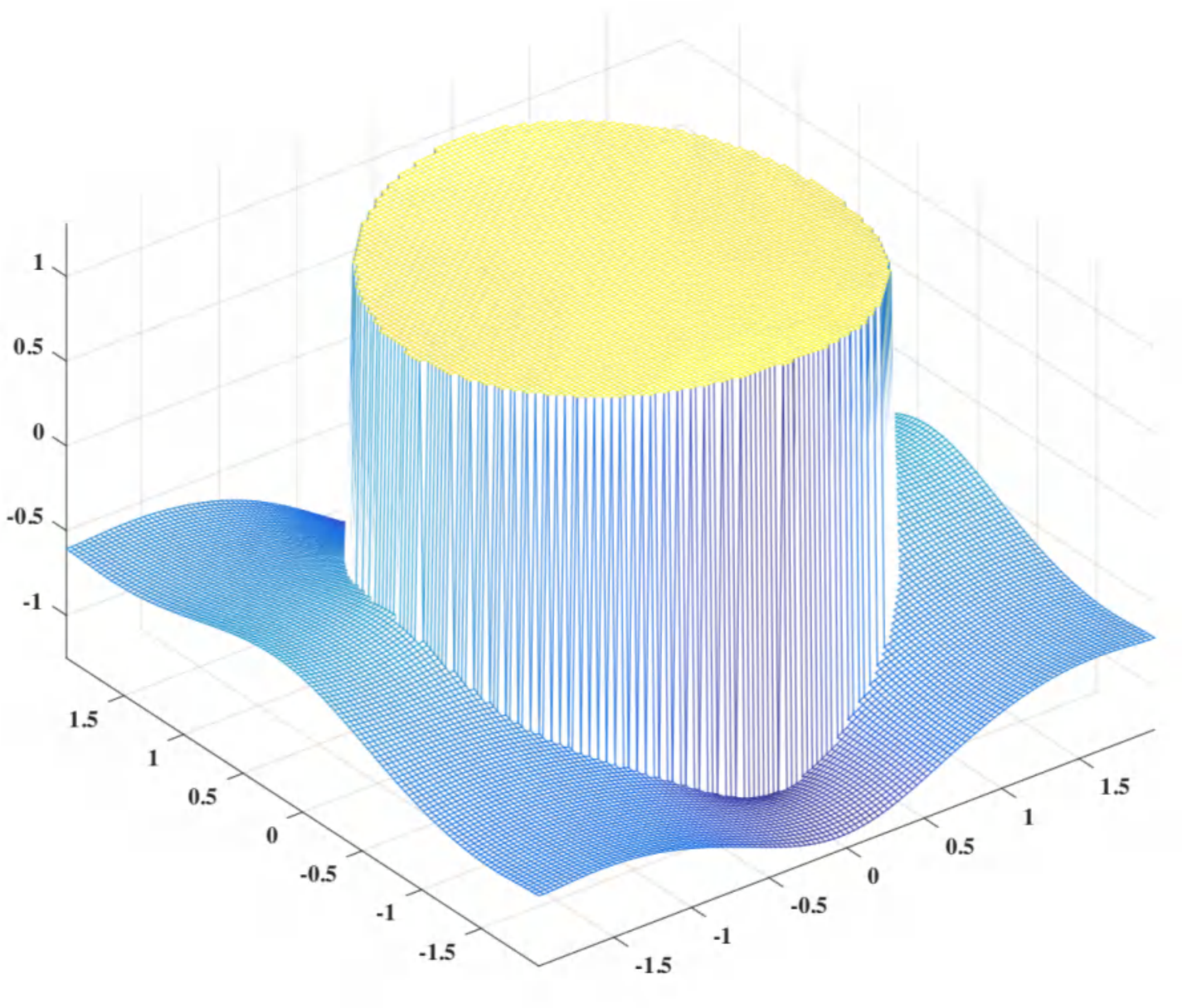}
}
\setlength{\abovecaptionskip}{-0.0cm}
\setlength{\belowcaptionskip}{-0.0cm}
 \caption{Example 6. The pressure distribution at different times.}
 \label{pressure-e4}
\end{figure}

{\em Example }6. In this example, an initial kidney-shaped interface given by $100$ control points is considered.  It is noted that the initial interface is also non-convex. The interface configurations at different moments ($t=0, t=0.374, t=0.935, t=1.87, t=7.48$) are shown in Fig. \ref{interface-2} (right). The tension coefficient $T_0$ is set to be 1.0.  In this example, the viscosity outside the interface is less than inside the interface, which is taken as $\mu^+ = 10$ and $\mu^- = 1$.  The approximation is computed up to a final time $T=8$ with $\Delta t = h = 0.0187$.  A time evolution of the velocity and the interface position are plotted in Fig.\ref{velocity-e4} (left),  isolines of the $x$-component $u^{(1)}$ and isolines of the $y$-component $u^{(2)}$ at different times are presented in Fig. \ref{velocity-e4} (middle) and Fig.\ref{velocity-e4} (right), respectively.  A time evolution of the pressure profile  is shown in Fig. \ref{pressure-e4}. Similar to the flower-shaped case, the velocity is continuous but not smooth, while the pressure is discontinuous across the interface, demonstrating  that the proposed method can handle more involved interface configurations, which are non-convex, too.

%%%%%%%%%%%%%%%%%%%%%%%%%%%%%%%%%%%%%%%%%%%%%%%%%
\section{Conclusion}
This work develops an efficient method to accurately solve Stokes equations involving two immiscible fluids with different viscosities and surface tension, in which the jump conditions for the velocity and the pressure are coupled together. The proposed technique is a combination of the KFBI method and the modified MAC scheme, which can be viewed as further development of the previous KFBI method and a stepping stone to more challenging cases.  By using boundary integral equations, the two-phase Stokes problems are reduced to the single-fluid case, the jump conditions of which can be decoupled easily, thus it can be solved using a modified MAC scheme in the usual way.  
Furthermore, in the iteration for the discrete BIEs, evaluation of the integrals are made indirectly by a Cartesian grid-based method.   
The major advantages of the presented method are that no augmented variable is needed and the matrix of the linear system to obtain the velocity and pressure approximation is unchanged. Thus some regular fast solver is applicable.

In addition, numerical results confirm that the proposed approach achieves second order accuracy for both velocity and pressure as well as the gradient of the velocity, even with a relatively large ratio $\mu^+/\mu^-$. Investigation of this method for moving interfaces is also considered, which suggests that  the proposed method is computationally  efficient and easy to implement. It can be extended to three dimensional problems, even more complex practical application problems, such as Stokes-Darcy problems, solid-fluid interaction, which will be future work.

%%%%%%%%%%%%%%%%%%%%%%%%%%%%%%%%%%%%%%%%%%%%%%%%%%%%%%%%%%%%

\section*{Appendix}
This appendix will illustrate and prove the equivalences between the interface problems and the volume or boundary integrals associated with the Green functions defined in  \eqref{Greendefine}. Here, only the proof for double layer potential is presented since the other two are similar and much easier.  Before proceeding the proof, Green's second identity is given by 
\begin{equation}
\label{G}
\int_{\Omega^+}(-\Delta \vect  u +\nabla p)\vect  vd\vect y- \int_{\Omega^+} (-\Delta \vect  v+\nabla q)\vect u d\vect y
=-\int_{\Gamma}T(\vect  u,p)\vect  vds_{\vect  y}  + \int_{\Gamma}T(\vect  v,q)\vect  uds_{\vect  y}.
\end{equation}

\begin{proof}
For the continuous function $\pmb\varphi(\vect x)$ defined on $\Gamma$, choose functions  $\vect  v_1(\vect  x)$  and $q_1(\vect  x)$ defined in $\Omega^+$ such that
\begin{equation}
\label{v1}
\begin{split}
-\Delta \vect  v_1+\nabla q_1&=\vect  0,  \;\,\;{\rm in}\; \Omega^+,\\
\nabla\cdot \vect  v_1&=0,\;\,\; {\rm in}\; \Omega^+,\\
\vect  v_1&=\pmb\varphi, \;\; {\rm on}\; \Gamma.
\end{split}
\end{equation}
Using the generalized Green second identity, one obtains 
\begin{equation*}
\begin{split}
\int_{\Omega^+}(-\Delta \vect  v_1+\nabla q_1)\vect  G_{\vect  v}d\vect  y - \int_{\Omega^+}(-\Delta \vect  G_{\vect  v}+\nabla G_q)\vect  v_1d\vect  y 
=-\int_{\Gamma}T(\vect  v_1, q_1)\vect  G_{\vect  v}ds_{\vect  y} +\int_{\Gamma} T(\vect  G_{\vect  v},G_q)\vect  v_1ds_{\vect y}.
\end{split}
\end{equation*}
By the first equation in \eqref{v1} and the definition of Green's pairs $(\vect G_{\vect v}, G_q)$ in \eqref{Greendefine}, one get
\begin{equation}
\label{G1}
\begin{split}
\int_{\Gamma} T(\vect  G_{\vect  v},G_q)\vect  v_1ds_{\vect y} = 
 \int_{\Gamma}T(\vect  v_1, q_1)\vect  G_{\vect  v}ds_{\vect  y}-\begin{cases}
\vect  v_1, \;{\rm if}\; \vect  x\in \Omega^+,\\
\vect 0, \;\;\;{\rm if}\; \vect  x\in\Omega^-,
\end{cases}
\end{split}
\end{equation}
Thus
\begin{equation}
\label{A1}
\begin{split}
\mathcal{M}_{\vect v}\pmb\varphi(\vect  x)&=\int_{\Gamma} T(\vect  G_{\vect  v},G_q)\pmb\varphi ds_{\vect  y} 
= \int_{\Gamma} T(\vect  G_{\vect  v},G_q)\vect  v_1ds_{\vect  y}
=\int_{\Gamma} T(\vect  v_1,q_1)\vect  G_{\vect  v}ds_{\vect  y}
-\begin{cases}
\vect  v_1, \;{\rm if}\; \vect  x\in \Omega^+,\\
0, \;\;\;{\rm if}\; \vect  x\in\Omega^-.
\end{cases}
\end{split}
\end{equation}
Recalling the fact  that $\vect  G_{\vect  v} = 0$ if $\vect  x\in\partial\Omega$, then 
$\int_{\Gamma} T(\vect  v_1,q_1)\vect  G_{\vect  v}ds_{\vect  y}=0,$ if  $\vect  x\in\partial\Omega$.
Thus, $\mathcal{M}_{\vect v}\pmb\varphi$  satisfies the boundary condition in \eqref{d-i}.

Substituting the relation $\nabla q_1 =\Delta \vect v_1$ into \eqref{G1}, we get
\begin{equation}
\label{G2}
\begin{split}
2\int_{\Gamma} \dfrac{\partial G_q}{\partial \vect  n_{\vect  y}}\cdot\vect  v_1ds_{\vect  y} = \int_{\Gamma}T(\vect  v_1,q_1)\cdot G_qds_{\vect  y}
-\begin{cases}
q_1, \;\,\;{\rm if}\; \vect  x\in \Omega^+,\\
0, \;\;\;\;{\rm if}\; \vect  x\in\Omega^-,
\end{cases}
\end{split}
\end{equation}
where identity \eqref{dd} has been used. Thus
\begin{equation}
\label{A2}
\begin{split}
\mathcal{M}_q\pmb\varphi(\vect  x)=2\int_{\Gamma} \dfrac{\partial G_q}{\partial \vect  n_{\vect  y}}\cdot\pmb\varphi ds_{\vect  y}
=2\int_{\Gamma} \dfrac{\partial G_q}{\partial \vect  n_{\vect  y}}\cdot\vect  v_1ds_{\vect  y} 
= \int_{\Gamma} T(\vect  v_1,q_1) \cdot  G_qds_{\vect  y}
-\begin{cases}
q_1, \;\,\;{\rm if}\; \vect  x\in \Omega^+,\\
0, \;\;\;\;{\rm if}\; \vect  x\in\Omega^-.
\end{cases}
\end{split}
\end{equation}
Applying the elliptic operator $-\Delta$ to identity \eqref{A1} and $\nabla$ to \eqref{A2}, then adding them together,  one can derive
\begin{equation*}
-\Delta(\mathcal{M}_{\vect v}\pmb\varphi)+\nabla (\mathcal{M}_q\pmb\varphi)=\int_{\Gamma} T(\vect  v_1,q_1)(-\Delta\vect  G_{\vect  v} + \nabla G_q)d\vect  y
-\begin{cases}
-\Delta \vect  v_1+\nabla q_1, \;{\rm if}\; \vect  x\in \Omega^+,\\
0, \qquad\qquad\;\quad{\rm if}\; \vect  x\in\Omega^-.
\end{cases}
\end{equation*}
That is, the double layer boundary integrals $\mathcal{M}_{\vect v}\pmb\varphi$ and $\mathcal{M}_q\pmb\varphi$ satisfy the  first equation with homogeneous right term in \eqref{d-i}.  Moreover,  by applying the elliptic operator $\nabla\cdot$ to \eqref{A1}, it is easy to see that the double layer boundary integral $\mathcal{M}_{\vect v}\pmb\varphi$ is incompressible.

Next, the discontinuity of the function $\mathcal{M}_{\vect v}\pmb\varphi$ and the continuity of the traction $T(\mathcal{M}_{\vect v}\pmb\varphi, \mathcal{M}_q\pmb\varphi)$ across $\Gamma$ will be illustrated.  First consider a constant density $\pmb\varphi=\vect  c$. By setting $\vect  u=\vect  c, p=0, \vect  v=\vect  G_{\vect  v}, q=G_q$ in equation \eqref{A1}, one readily sees that 
\begin{equation}
\label{C}
\mathcal{M}_{\vect v}\pmb\varphi(\vect  x)=\begin{cases}
-\vect  c,\;\;\vect  x\in\Omega^+,\\
\vect  0,\quad\;\vect  x\in\Omega^-.
\end{cases}
\end{equation}
A further application of Green's second identity \eqref{G} on $\Omega^+_{\epsilon}$, which is the part of $\Omega^+$ remaining after deleting $K(\vect  x,\epsilon)$ with radius $\epsilon$ and center $\vect  x\in \Gamma$,
shows that 
\begin{equation*}
\int_{S_{\epsilon}+C_{\epsilon}} \Big(T(\vect  G_{\vect  v},G_q)\vect  cds_{\vect  y} - T(\vect  c, 0)\vect  G_{\vect  v}\Big)ds_{\vect  y} = 0.
\end{equation*}
Here $C_{\epsilon}$ denotes the part of the surface of the circle $K(\vect  x, \epsilon)$ that is contained in $\Omega^+$, and $S_{\epsilon}$ is the piece of the interface $\Gamma$ remaining after deleting $K(\vect  x, \epsilon)$.  Since $T(\vect  c,0)=0$, the following identity holds
\begin{equation*}
\vect  v(\vect  x) = - \mathcal{M}_{\vect v}\pmb\varphi(\vect x)= -\lim\limits_{\epsilon\rightarrow 0}\int_{S_{\epsilon}} T(\vect  G_{\vect  v},G_q)\vect  cds_{\vect  y} = \lim\limits_{\epsilon\rightarrow 0}\int_{C_{\epsilon}} T(\vect  G_{\vect  v},G_q)\vect  cds_{\vect  y}.
\end{equation*}
As done before, if the integral is carried over the whole surface of the circle $K(\vect  x, \epsilon)$, which is defined as 
$C_{(x, \epsilon)}$, one can derive the following identity 
\begin{equation*}
\int_{C_{(x, \epsilon)}} T(\vect  G_{\vect  v},G_q)\vect  cds_{\vect  y} = -\int_{K(\vect  x, \epsilon)}(-\Delta \vect G_{\vect v}+\nabla G_q)\vect cd\vect y +\int_{C_{(x, \epsilon)}} T(\vect c, 0)\vect G_{\vect v}ds_{\vect y} = -\vect  c,
\end{equation*}
for any $\vect x\in K(\vect  x, \epsilon)$.  Since $T(\vect  G_{\vect  v},G_q)$ is symmetric, the  above integral over a half of $C_{(x, \epsilon)}$ equals $-\frac{1}{2}\vect  c$. While as $\epsilon\rightarrow 0$, the integral over $C_{\epsilon}$ approaches the integral over the semi-circle. Thus,  
\begin{equation}
\label{A3}
\vect  v(\vect  x) = - \mathcal{M}_{\vect v}\pmb\varphi(\vect x)=\begin{cases}
\vect  c,\quad\vect  x\in\Omega^+,\\[4pt]
\dfrac{1}{2}\vect  c,\;\vect  x\in\Gamma,\\[4pt]
\vect  0,\quad\vect  x\in\Omega^-.
\end{cases}
\end{equation}
In order to determine values of the jumps of $\vect  v$ on $\Gamma$ for a continuous density $\pmb\varphi$, one needs to show the continuity of the following function 
\begin{equation}
\label{A4}
\widetilde{\vect  v}(\vect  x) = - \int_{\Gamma} T(\vect  G_{\vect  v},G_q)\pmb\varphi(\vect  y) ds_{\vect  y} + \pmb\varphi(\vect  y_0)\int_{\Gamma} T(\vect  G_{\vect  v},G_q)ds_{\vect  y},  \quad\vect y_0\in \Gamma.
\end{equation}
Actually,  provided that 
\begin{equation*}
 \int_{\Gamma} |T(\vect  G_{\vect  v},G_q)| d\vect  y\leq {\rm Const}, 
\end{equation*}
for $\vect x \in \Omega$, $\widetilde{\vect v}(\vect x)$ is continuous. 
It is noted that this assumption is reasonable because that the difference between Green's function defined in \eqref{Greendefine} and the fundamental solution of the Stokes system in the free space  is a harmonic function defined on the domain $\Omega$ and the integral of the fundamental solution on $\Gamma$ is bounded \cite{ladyzhenskaya1969mathematical}.  

Taking the limiting values on $\Gamma$ from inside of $\Omega$ and outside of $\Omega$ about the equation \eqref{A4}, one can obtain 
the following relations
\begin{equation}
\label{A5}
\begin{split}
\widetilde{\vect v}^+(\vect x) &= -(\mathcal{M}_{\vect v}\pmb\varphi)^+(\vect x) -  \pmb\varphi(\vect x) = \vect v^+(\vect x) -  \pmb\varphi(\vect x),\\
\widetilde{\vect v}^-(\vect x) &= -(\mathcal{M}_{\vect v}\pmb\varphi)^-(\vect x) -0 = \vect v^-(\vect x).
\end{split}
\end{equation}
The jump relation $[\![ \vect v(\vect x)]\!] = \pmb\varphi(\vect x)$ holds because the continuity of function $\widetilde{\vect v}$.
Moreover, 
\begin{equation}
\label{A6}
\widetilde{\vect  v}(\vect x)\Big|_{\Gamma} = -(\mathcal{M}_{\vect v}\pmb\varphi)(\vect x) -\dfrac{1}{2} \pmb\varphi(\vect x).
\end{equation}
Combining \eqref{A5} and \eqref{A6}, 
one can derive
\begin{equation*}
\begin{split}
\vect v^+(\vect x) &=  \frac{1}{2} \pmb\varphi(\vect x) -(\mathcal{M}_{\vect v}\pmb\varphi)(\vect x),\\[4pt]
\vect v^-(\vect x) &=  -\frac{1}{2} \pmb\varphi (\vect x)-(\mathcal{M}_{\vect v}\pmb\varphi)(\vect x).
\end{split}
\end{equation*}

In addition, note that the normal flux $\vect  \nabla\widetilde{\vect  v}$ is continuous across the interface $\Gamma$,  and the flux $\vect  n\nabla(\int_{\Gamma} T(\vect  G_{\vect  v},G_q)\pmb\varphi(\vect  x_0) d\vect  y)$  is also continuous across the interface by \eqref{A3}.
Thus, the normal flux $\nabla (\mathcal{M}_{\vect v}\pmb\varphi)$ is continuous across the interface $\Gamma$.  

Furthermore, following the line in proving the discontinuity of the function $\mathcal{M}_{\vect v}\varphi(\vect x)$,   one can show the continuity of the double layer integral $\mathcal{M}_q\varphi(\vect x)$. 
This gives the continuity of traction $T(\mathcal{M}_{\vect v}, \mathcal{M}_q)$, which ends the proof.

\end{proof}

\section*{Acknowledgement}
Haixia Dong is partially supported by NSFC under Grant NO. 12001193, the Scientific Research Fund of Hunan Provincial Education Department (No.20B376), Changsha Municipal Natural Science Foundation (No. kq2014073). Wenjun Ying is partially supported by the Strategic Priority Research Program of Chinese Academy of Sciences (Grant No. XDA25010405), the National Natural Science Foundation of China (Grant No. DMS-11771290) and the Science Challenge Project of China (Grant No. TZ2016002). 

%Jiwei Zhang is partially supported by NSFC under grant No. 12171376, 2020-JCJQ- ZD-029 and NSAF U1930402.

%============================================================================== 
%%****************************************************************************** 
%%\bibliographystyle{plain}
%\bibliographystyle{plain}      % basic style, author-year citations
%%\bibliographystyle{spmpsci}      % mathematics and physical sciences
%%\bibliographystyle{spphys}       % APS-like style for physics
%\bibliography{references}  

\begin{thebibliography}{10}

\bibitem{adjerid2015immersed}
Slimane Adjerid, Nabil Chaabane, and Tao Lin.
\newblock An immersed discontinuous finite element method for Stokes interface
  problems.
\newblock {\em Computer Methods in Applied Mechanics and Engineering},
  293:170--190, 2015.

\bibitem{adjerid2019immersed}
Slimane Adjerid, Nabil Chaabane, Tao Lin, and Pengtao Yue.
\newblock An immersed discontinuous finite element method for the Stokes
  problem with a moving interface.
\newblock {\em Journal of Computational and Applied Mathematics}, 362:540--559,
  2019.

\bibitem{Beale2004grid}
J.~Thomas Beale.
\newblock A grid-based boundary integral method for elliptic problems in three
  dimensions.
\newblock {\em SIAM Journal on Numerical Analysis}, 42(2):599--620, 2004.

\bibitem{chang1996level}
Yu-Chung Chang, TY~Hou, B~Merriman, and Stanley Osher.
\newblock A level set formulation of Eulerian interface capturing methods for
  incompressible fluid flows.
\newblock {\em Journal of computational Physics}, 124(2):449--464, 1996.

\bibitem{chen2018direct}
Xiaohong Chen, Zhilin Li, and Juan~Ruiz {\'A}lvarez.
\newblock A direct IIM approach for two-phase Stokes equations with
  discontinuous viscosity on staggered grids.
\newblock {\em Computers \& Fluids}, 172:549--563, 2018.

\bibitem{chen2021p2}
Yuan Chen and Xu~Zhang.
\newblock A p2-p1 partially penalized immersed finite element method for Stokes
  interface problems.
\newblock {\em International journal of numerical analysis and modeling},
  18(1), 2021.

\bibitem{chessa2003extended}
Jack Chessa and Ted Belytschko.
\newblock An extended finite element method for two-phase fluids.
\newblock {\em J. Appl. Mech.}, 70(1):10--17, 2003.

\bibitem{cogan2005modeling}
NG~Cogan, Ricardo Cortez, and Lisa Fauci.
\newblock Modeling physiological resistance in bacterial biofilms.
\newblock {\em Bulletin of mathematical biology}, 67(4):831--853, 2005.

\bibitem{cortez2001method}
Ricardo Cortez.
\newblock The method of regularized Stokeslets.
\newblock {\em SIAM Journal on Scientific Computing}, 23(4):1204--1225, 2001.

\bibitem{dong2018hybridizable}
Haixia Dong, Wenjun Ying, and Jiwei Zhang.
\newblock A hybridizable discontinuous Galerkin method for elliptic interface
  problems in the formulation of boundary integral equations.
\newblock {\em Journal of Computational and Applied Mathematics}, 344:624--639,
  2018.

\bibitem{dong2022second}
Haixia Dong, Zhongshu Zhao, Shuwang Li, Wenjun Ying, and Jiwei Zhang.
\newblock Second order convergence of a modified MAC scheme for Stokes
  interface problem.
\newblock {\em Preprint}.

\bibitem{fogelson1986numerical}
AL~Fogelson and CS~Peskin.
\newblock Numerical solution of the three-dimensional Stokes' equations in the
  presence of suspended particles.
\newblock In {\em Unknown Host Publication Title}. Soc. Ind. \& Appl. Math,
  1986.

\bibitem{gross2007extended}
Sven Gro{\ss} and Arnold Reusken.
\newblock An extended pressure finite element space for two-phase
  incompressible flows with surface tension.
\newblock {\em Journal of Computational Physics}, 224(1):40--58, 2007.

\bibitem{gross2007finite}
Sven Gross and Arnold Reusken.
\newblock Finite element discretization error analysis of a surface tension
  force in two-phase incompressible flows.
\newblock {\em SIAM journal on numerical analysis}, 45(4):1679--1700, 2007.

\bibitem{gross2011numerical}
Sven Gross and Arnold Reusken.
\newblock {\em Numerical methods for two-phase incompressible flows},
  volume~40.
\newblock Springer Science \& Business Media, 2011.

\bibitem{hansbo2014cut}
Peter Hansbo, Mats~G Larson, and Sara Zahedi.
\newblock A cut finite element method for a Stokes interface problem.
\newblock {\em Applied Numerical Mathematics}, 85:90--114, 2014.

\bibitem{he2019stabilized}
Xiaoxiao He, Fei Song, and Weibing Deng.
\newblock A stabilized nonconforming Nitsche's extended finite element method
  for Stokes interface problems.
\newblock {\em arXiv preprint arXiv:1905.04844}, 2019.

\bibitem{hou2012numerical}
Gene Hou, Jin Wang, and Anita Layton.
\newblock Numerical methods for fluid-structure interaction—a review.
\newblock {\em Communications in Computational Physics}, 12(2):337--377, 2012.

\bibitem{ji2022immersed}
Haifeng Ji, Feng Wang, Jinru Chen, and Zhilin Li.
\newblock An immersed CR-p0 element for Stokes interface problems and the
  optimal convergence analysis.
\newblock {\em Computer Methods in Applied Mechanics and Engineering},
  399:115306, 2022.

\bibitem{jones2021class}
Derrick Jones and Xu~Zhang.
\newblock A class of nonconforming immersed finite element methods for Stokes
  interface problems.
\newblock {\em Journal of Computational and Applied Mathematics}, 392:113493,
  2021.

\bibitem{kim2019immersed}
Woojin Kim and Haecheon Choi.
\newblock Immersed boundary methods for fluid-structure interaction: A review.
\newblock {\em International Journal of Heat and Fluid Flow}, 75:301--309,
  2019.

\bibitem{kirchhart2016analysis}
Matthias Kirchhart, Sven Gross, and Arnold Reusken.
\newblock Analysis of an XFEM discretization for Stokes interface problems.
\newblock {\em SIAM Journal on Scientific Computing}, 38(2):A1019--A1043, 2016.

\bibitem{kress1989linear}
Rainer Kress, V~Maz'ya, and V~Kozlov.
\newblock {\em Linear integral equations}, volume~17.
\newblock Springer, 1989.

\bibitem{ladyzhenskaya1969mathematical}
Olga~A Ladyzhenskaya and Richard~A Silverman.
\newblock {\em The mathematical theory of viscous incompressible flow},
  volume~12.
\newblock Gordon \& Breach New York, 1969.

\bibitem{laymuns2022corrected}
Genaro Laymuns and Manuel~A S{\'a}nchez.
\newblock Corrected finite element methods on unfitted meshes for Stokes moving
  interface problem.
\newblock {\em Computers \& Mathematics with Applications}, 108:159--174, 2022.

\bibitem{layton2008efficient}
Anita~T Layton.
\newblock An efficient numerical method for the two-fluid Stokes equations with
  a moving immersed boundary.
\newblock {\em Computer Methods in Applied Mechanics and Engineering},
  197(25):2147--2155, 2008.

\bibitem{lee2003immersed}
Long Lee and Randall~J LeVeque.
\newblock An immersed interface method for incompressible Navier--Stokes
  equations.
\newblock {\em SIAM Journal on Scientific Computing}, 25(3):832--856, 2003.

\bibitem{lehrenfeld2012nitsche}
Christoph Lehrenfeld and Arnold Reusken.
\newblock Nitsche-XFEM with streamline diffusion stabilization for a two-phase
  mass transport problem.
\newblock {\em SIAM journal on scientific computing}, 34(5):A2740--A2759, 2012.

\bibitem{leveque1994immersed}
Randall~J Leveque and Zhilin Li.
\newblock The immersed interface method for elliptic equations with
  discontinuous coefficients and singular sources.
\newblock {\em SIAM Journal on Numerical Analysis}, 31(4):1019--1044, 1994.

\bibitem{leveque1997immersed}
Randall~J LeVeque and Zhilin Li.
\newblock Immersed interface methods for Stokes flow with elastic boundaries or
  surface tension.
\newblock {\em SIAM Journal on Scientific Computing}, 18(3):709--735, 1997.

\bibitem{li2001maximum}
Zhilin Li and Kazufumi Ito.
\newblock Maximum principle preserving schemes for interface problems with
  discontinuous coefficients.
\newblock {\em SIAM Journal on Scientific Computing}, 23(1):339--361, 2001.

\bibitem{li2006immersed}
Zhilin Li and Kazufumi Ito.
\newblock {\em The immersed interface method: numerical solutions of PDEs
  involving interfaces and irregular domains}, volume~33.
\newblock Siam, 2006.

\bibitem{li2007augmented}
Zhilin Li, Kazufumi Ito, and Ming-Chih Lai.
\newblock An augmented approach for Stokes equations with a discontinuous
  viscosity and singular forces.
\newblock {\em Computers \&amp; Fluids}, 36(3):622--635, 2007.

\bibitem{li2017accurate}
Zhilin Li, Haifeng Ji, and Xiaohong Chen.
\newblock Accurate solution and gradient computation for elliptic interface
  problems with variable coefficients.
\newblock {\em SIAM journal on numerical analysis}, 55(2):570--597, 2017.

\bibitem{li2001immersed}
Zhilin Li and Ming-Chih Lai.
\newblock The immersed interface method for the Navier--Stokes equations with
  singular forces.
\newblock {\em Journal of Computational Physics}, 171(2):822--842, 2001.

\bibitem{lundberg2019distributed}
Andrew Lundberg, Pengtao Sun, and Cheng Wang.
\newblock Distributed Lagrange multiplier-fictitious domain finite element
  method for Stokes interface problems.
\newblock {\em Int. J. Numer. Anal. Model}, 16(6):939--963, 2019.

\bibitem{mayo1992implicit}
Anita~A Mayo and Charles~S Peskin.
\newblock An implicit numerical method for fluid dynamics problems with
  immersed elastic boundaries.
\newblock {\em Contemporary Mathematics}, 141:261--261, 1992.

\bibitem{mokbel2018phase}
Dominic Mokbel, Helmut Abels, and Sebastian Aland.
\newblock A phase-field model for fluid--structure interaction.
\newblock {\em Journal of computational physics}, 372:823--840, 2018.

\bibitem{peskin1977numerical}
Charles~S Peskin.
\newblock Numerical analysis of blood flow in the heart.
\newblock {\em Journal of computational physics}, 25(3):220--252, 1977.

\bibitem{peskin2002immersed}
Charles~S Peskin.
\newblock The immersed boundary method.
\newblock {\em Acta numerica}, 11:479--517, 2002.

\bibitem{Saad1993GMRES}
Youcef Saad.
\newblock A flexible inner-outer preconditioned GMRES algorithm.
\newblock {\em SIAM Journal on Scientific Computing}, 14(2):461--469, 1993.

\bibitem{saad1986gmres}
Youcef Saad and Martin~H Schultz.
\newblock Gmres: A generalized minimal residual algorithm for solving
  nonsymmetric linear systems.
\newblock {\em SIAM Journal on scientific and statistical computing},
  7(3):856--869, 1986.

\bibitem{sun2019fictitious}
Pengtao Sun.
\newblock Fictitious domain finite element method for Stokes/elliptic interface
  problems with jump coefficients.
\newblock {\em Journal of Computational and Applied Mathematics}, 356:81--97,
  2019.

\bibitem{sun2020distributed}
Pengtao Sun and Cheng Wang.
\newblock Distributed Lagrange multiplier/fictitious domain finite element
  method for Stokes/parabolic interface problems with jump coefficients.
\newblock {\em Applied Numerical Mathematics}, 152:199--220, 2020.

\bibitem{tan2009immersed}
Zhijun Tan, Duc-Vinh Le, KM~Lim, and BC~Khoo.
\newblock An immersed interface method for the incompressible Navier--Stokes
  equations with discontinuous viscosity across the interface.
\newblock {\em SIAM Journal on Scientific Computing}, 31(3):1798--1819, 2009.

\bibitem{tan2011implementation}
Zhijun Tan, KM~Lim, and BC~Khoo.
\newblock An implementation of MAC grid-based IIM-Stokes solver for
  incompressible two-phase flows.
\newblock {\em Communications in Computational Physics}, 10(5):1333--1362,
  2011.

\bibitem{tu1992stability}
Cheng Tu and Charles~S Peskin.
\newblock Stability and instability in the computation of flows with moving
  immersed boundaries: a comparison of three methods.
\newblock {\em SIAM Journal on Scientific and Statistical Computing},
  13(6):1361--1376, 1992.

\bibitem{wang2013hybridizable}
Bo~Wang and BC~Khoo.
\newblock Hybridizable discontinuous Galerkin method (HDG) for Stokes interface
  flow.
\newblock {\em Journal of Computational Physics}, 247:262--278, 2013.

\bibitem{wang2019nonconforming}
Nan Wang and Jinru Chen.
\newblock A nonconforming Nitsche’s extended finite element method for Stokes
  interface problems.
\newblock {\em Journal of Scientific Computing}, 81(1):342--374, 2019.

\bibitem{wang2015new}
Qiuliang Wang and Jinru Chen.
\newblock A new unfitted stabilized Nitsche’s finite element method for
  Stokes interface problems.
\newblock {\em Computers \& Mathematics with Applications}, 70(5):820--834,
  2015.

\bibitem{xie2019fourth}
Yaning Xie and Wenjun Ying.
\newblock A fourth-order kernel-free boundary integral method for the modified
  Helmholtz equation.
\newblock {\em Journal of Scientific Computing}, 78(3):1632--1658, 2019.

\bibitem{xie2019high}
Yaning Xie, Wenjun Ying, and Wei-Cheng Wang.
\newblock A high-order kernel-free boundary integral method for the biharmonic
  equation on irregular domains.
\newblock {\em Journal of Scientific Computing}, 80(3):1681--1699, 2019.

\bibitem{xu20083d}
Sheng Xu and Z~Jane Wang.
\newblock A 3d immersed interface method for fluid--solid interaction.
\newblock {\em Computer Methods in Applied Mechanics and Engineering},
  197(25):2068--2086, 2008.

\bibitem{ying2013fast}
Wenjun Ying and J~Thomas Beale.
\newblock A fast accurate boundary integral method for potentials on closely
  packed cells.
\newblock {\em Communications in Computational Physics}, 14(04):1073--1093,
  2013.

\bibitem{ying2007kernel}
Wenjun Ying and Craig~S Henriquez.
\newblock A kernel-free boundary integral method for elliptic boundary value
  problems.
\newblock {\em Journal of computational physics}, 227(2):1046--1074, 2007.

\bibitem{ying2013kernel}
Wenjun Ying and Wei-Cheng Wang.
\newblock A kernel-free boundary integral method for implicitly defined
  surfaces.
\newblock {\em Journal of Computational Physics}, 252:606--624, 2013.

\bibitem{ying2014kernel}
Wenjun Ying and Wei-Cheng Wang.
\newblock A kernel-free boundary integral method for variable coefficients
  elliptic PDEs.
\newblock {\em Communications in Computational Physics}, 15(04):1108--1140,
  2014.

\bibitem{zhang1997immersed}
Chaoming Zhang and Randall~J LeVeque.
\newblock The immersed interface method for acoustic wave equations with
  discontinuous coefficients.
\newblock {\em Wave motion}, 25(3):237--263, 1997.

\end{thebibliography}
%============================================================================== 

%

%============================================================================== 
\end{document}